\documentclass{amsart}[12pt]
\usepackage{amsmath,amssymb,cite,geometry,bm,graphicx,pgfplots}
\usepackage{amsmath,amssymb,amscd,cite,geometry,pgfplots,bm,easybmat,enumerate,graphicx}
\usepackage[colorlinks,citecolor=blue,linkcolor=blue,urlcolor=blue]{hyperref}
\numberwithin{equation}{section}
\allowdisplaybreaks[4]

\usepackage{mathrsfs}

\begin{document}

\title{Remarks on paper ``Two-term spectral asymptotics in linear elasticity''}

\author{Genqian Liu}
\address{School of Mathematics and Statistics, Beijing Institute of Technology, Beijing 100081, China}
\email{liugqz@bit.edu.cn, or liuswgq@163.com}

\subjclass[2020]{58J50 (Primary) 53C21, 35P20 (Secondary)} 
\keywords{Spectral problems,  Riemannian geometry,  Asymptotic distribution of eigenvalues, Linear elasticity}

\begin{abstract}
\   In this note, we shall point out that all ``numerically calculations'' and figures  in \cite{CaFrLeVa-23} are wrong because these calculations are based on some incorrect formulas. Furthermore, by pointing out several serious errors in \cite{CaFrLeVa-23} and especially by Section 7,  Proposition 7.1, Remarks 7.2--7.3, and Section 8  (a result of A. Pierzchalski and B. {\O}rsted)  we show that the conclusions published by Matteo Capoferri, Leonid Friedlander, Michael Levitin and Dmitri Vassiliev (J Geom Anal (2023)33:242) as well as the main ``algorithm'' theory of the book \cite{SaVa-97} are completely wrong. Finally, we explain the correctness of proof of Theorem 1.1 in our paper \cite{Liu-21} by giving some remarks and putting the whole proof in Appendix (see also \cite{Liu-22b} and \cite{Liu-22c}).
\end{abstract}

\maketitle
\vskip 0.29 true cm
\maketitle
\vskip 0.19 true cm
\section{Introduction}

 \vskip 0.49 true cm

Let $(\Omega,g)$ be a compact smooth Riemannian $n$-manifold with smooth boundary $\partial \Omega$. Let $P_g$ be the Navier--Lam\'{e} operator (see \cite{Liu-19} and \cite{Liu-21}):
\begin{eqnarray} \label {1-1} P_g\mathbf{u}:=\mu \nabla^* \nabla \mathbf{u} -(\mu +\lambda) \,\mbox{grad}\; \mbox{div}\, \mathbf{u} -\mu \, \mbox{Ric} (\mathbf{u}), \;  \;\; \mathbf{u}=(u^1, \cdots, u^n),\end{eqnarray}
where $\mu$ and $\lambda$ are the Lam\'{e} parameters satisfying $\mu>0$ and $\mu+\lambda\ge  0$,  $(\nabla^* \nabla \mathbf{u})^k:=-\sum_{j,k=1}^n \nabla_j \nabla^j u^k $ is the Bochner Laplacian (see (2.11) of \cite{Liu-21}), $\mbox{div}$ and $\mbox{grad}$ are the usual divergence and gradient operators, and  \begin{eqnarray} \label{18/12/22} \mbox{Ric} (\mathbf{u})= \big(\sum\limits_{l=1}^nR^{\,1}_{l} u^l,  \sum\limits_{l=1}^n R^{\,\,2}_{\,l} u^l, \cdots, \sum\limits_{l=1}^n R^{\,\,n}_{\,l} u^l\big)\end{eqnarray} denotes the action of Ricci tensor $\mbox{R}_l^{\;j}:=\sum_{k=1}^n R^{k\,\,j}_{\,lk}$ on $\mathbf{u}$.
We denote by $P_g^-$ and $P_g^+$ the Navier--Lam\'{e} operators with zero Dirichlet and zero Neumann boundary conditions, respectively. The Dirichlet boundary condition is $\mathbf{u}\big|_{\partial \Omega}$, and the Neumann (i.e., free) boundary condition is \begin{eqnarray*}2\mu  (\mbox{Def}\, \mathbf{u})^\# \boldsymbol{\nu} + \lambda (\mbox{div}\, \mathbf{u} )\boldsymbol{\nu}\;\;\mbox{on}\;\; \partial \Omega,\end{eqnarray*} where $\mbox{Def}\, \mathbf{u}= \frac{1}{2} (\nabla \mathbf{u} +(\nabla \mathbf{u})^T)$, $\,(\nabla \mathbf{u})^T$ is the transpose of $\nabla \mathbf{u}$. $ \;\#$ is the sharp operator (for a tensor) by raising index, and $\boldsymbol{\nu}$ is the unit outer normal to $\partial \Omega$.
Since $P_g^-$ (respectively, $P_g^+$) is an unbounded, self-adjoint and positive (respectively,  nonnegative) operator in $[H^1_0(\Omega)]^n$ (respectively, $[H^1(\Omega)]^n$) with discrete spectrum $0< \tau_1^- < \tau_2^- \le \cdots \le \tau_k^- \le \cdots \to +\infty$ (respectively, $0\le \tau_1^+ < \tau_2^+ \le \cdots \le \tau_k^+ \le \cdots \to +\infty$), one has
\begin{eqnarray} \label{1-4} P_g^\mp {\mathbf{u}}_k^\mp =\tau_k^\mp {\mathbf{u}}_k^\mp,\end{eqnarray} where ${\mathbf{u}}_k^-\in  [H^1_0(\Omega)]^n$ (respectively, ${\mathbf{u}}_k^+\in  [H^1(\Omega)]^n$) is the eigenvector corresponding to elastic eigenvalue $\tau_k^{-}$ (respectively, $\tau_k^{+}$).

\vskip 0.10 true cm
We introduce the {\it partition function}, or  the {\it trace of the heat semigroup} for the Lam\'{e} operator with zero Dirichlet (respectively, zero Neumann) boundary condition, by
$\mathcal{Z}^- (t) := \mbox{Tr}\; e^{-t P_g^-} = \sum_{k=1}^\infty e^{-t \tau_k^{-}}$ (respectively, $\mathcal{Z}^+ (t) := \mbox{Tr}\; e^{-t P_g^+} = \sum_{k=1}^\infty e^{-t \tau_k^{+}}$)
defined for $t>0$ and monotone decreasing in $t$.
In \cite{Liu-21}, by using the method of the heat trace and ``method of image''  we obtained the following result:

\vskip 0.22 true cm

\noindent{\bf Theorem 1.1.} \ {\it
Let $(\Omega,g)$ be a smooth compact Riemannian manifold of dimension $n$ with smooth boundary $\partial \Omega$, and let $0< \tau_1^-< \tau_2^- \le \tau^-_3\le \cdots \le \tau_k^- \le \cdots$ (respectively, $0\le \tau_1^+ < \tau_2^+ \le \tau_3^+ \le \cdots \le \tau_k^+ \le \cdots $) be the eigenvalues of the Navier--Lam\'{e} operator $P_g^-$ (respectively, $P_g^+$) with respect to the zero Dirichlet (respectively, zero Neumann) boundary condition. Then
\begin{eqnarray} \label{1-7} && \mathcal{Z}^{\mp}(t) =  \sum_{k=1}^\infty e^{-t \tau_k^{\mp}} =\bigg[ \frac{n-1}{(4\pi \mu t)^{n/2}}
 + \frac{1}{(4\pi (2\mu+\lambda) t)^{n/2}}\bigg] {\mbox{Vol}_n}(\Omega) \\
&& \qquad \;\, \quad \, \mp \frac{1}{4} \bigg[  \frac{n-1}{(4\pi \mu t)^{(n-1)/2}}
 +  \frac{1}{(4\pi (2\mu+\lambda) t)^{(n-1)/2}}\bigg]{\mbox{Vol}_{n-1}}(\partial\Omega)+O(t^{{1-n}/2})\quad\;\; \mbox{as}\;\; t\to 0^+.\nonumber\end{eqnarray}
 Here ${\mbox{Vol}}_{n}(\Omega)$ denotes the $n$-dimensional volume of $\Omega$,  ${\mbox{Vol}}_{n-1} (\partial\Omega)$ denotes the $(n-1)$-dimensional volume of $\partial \Omega$. }

 \vskip 0.16 true cm
For any $\Lambda\in (-\infty, +\infty)$, denote by $\mathcal{N}^{-} (\Lambda) := \max\limits_{\substack{ k}} \{ k\big| \tau^{-}_k < \Lambda\}$ (respectively, $\mathcal{N}^{+} (\Lambda) :=\max\limits_{\substack{ k}} \{ k\big| \tau^{+}_k < \Lambda\}$) the {\it eigenvalue counting function} for the elastic Lam\'{e} operator with zero Dirichlet (respectively, zero Neumann) boundary condition. In \cite{CaFrLeVa-23}, Matteo Capoferri, Leonid Friedlander, Michael Levitin and Dmitri Vassiliev proved the following:

 \vskip 0.16 true cm
\noindent{\bf Theorem 1.2.} \ {\it
Let $(\Omega,g)$ be a smooth compact connected $n$-dimensional Riemannian manifold with smooth boundary $\partial \Omega$. Suppose that $(\Omega, g)$ is such that the corresponding billiards is neither dead-end nor absolutely periodic.
 Then
\begin{eqnarray} \label{2023.5.5-1}  \mathcal{N}^\mp (\Lambda) = a \mbox{Vol}_n (\Omega) \Lambda^{n/2} + b^{\mp} \mbox{Vol}_{n-1} (\partial \Omega) \Lambda^{(n-1)/2}+ o(\Lambda^{(n-1)/2}) \;\; \mbox{as}\;\; \Lambda\to +\infty,\end{eqnarray}
where
\begin{eqnarray}\label{2023.5.5-2}\!\!\!\!\!\!\!\! \!\!\!\!&& a= \frac{1}{(4\pi)^{n/2}\, \Gamma (1+\frac{n}{2})} \Big(\frac{n-1}{\mu^{n/2}}+ \frac{1}{(\lambda+2\mu)^{n/2}}\Big),\\
\!\!\!\!\!\!\!\!\!\!\!\!&&\label{202.5.5-10} b^{-} =- \frac{ \mu^{\frac{1-n}{2}}}{ 2^{n+1} \pi^{\frac{n-1}{2}} \Gamma(\frac{n+1}{2})} \bigg( \frac{4(n-1)}{\pi} \int_{\sqrt{\alpha}}^1 \tau^{n-2} \arctan \Big(\sqrt{(1-\alpha \tau^{-2})(\tau^{-2}-1)} \Big) d\tau \\
 &&\quad \;\;\quad + \alpha^{\frac{n-1}{2}} +n-1\bigg),\nonumber\\
\!\!\!\!\!\!\!\!\!\!\!\!&& \label{202.5.5-3}  b^{+} = \frac{ \mu^{\frac{1-n}{2}}}{ 2^{n+1} \pi^{\frac{n-1}{2}} \Gamma(\frac{n+1}{2})} \!\bigg(\! \frac{4(n-1)}{\pi} \!\int_{\sqrt{\alpha}}^1 \!\tau^{n-2} \arctan \!\Big(\frac{(\tau^{-2} \!-2)^2}{4\sqrt{(1-\!\alpha \tau^{-2})(\tau^{-2}\!-1)}} \Big) d\tau \\
&& \quad\;\;\quad + \alpha^{\frac{n-1}{2}} +n-5+4 \gamma_R^{1-n}\!\bigg),\nonumber\end{eqnarray}
 where $\alpha:= \frac{\mu}{\lambda +2\mu}$, $\gamma_R:= \sqrt{\omega_1}$, and $w_1$ is the distinguished real root of the cube equation $R_\alpha (w):=w^3-8w^2 +8(3-2\alpha) w +16 (\alpha-1)=0$ in the interval $(0,1)$.}

\vskip 0.16 true cm
Clearly, (\ref{1-7}) can be rewritten as
\begin{eqnarray}\label{2023.5.5-8}  \mathcal{Z}^{\mp}(t) = \tilde{a} \mbox{Vol}_n( \Omega) t^{-n/2} + \tilde{b}^{\mp} \mbox{Vol}_{n-1} (\partial \Omega) t^{-(n-1)/2} + O(t^{1-n/2})\;\; \mbox{as} \;\; t\to 0^+,\end{eqnarray}
where \begin{eqnarray*}&& \tilde{a}= \frac{1}{(4\pi)^{n/2}} \Big( \frac{n-1}{\mu^{n/2}} +\frac{1}{ (\lambda+ 2\mu)^{n/2}}\Big), \\
&&  \tilde{b}^{\mp} = \mp \frac{ 1}{ 4(4\pi)^{(n-1)/2}}\bigg( \frac{n-1}{\mu^{(n-1)/2} } + \frac{1}{ (\lambda+2\mu)^{(n-1)/2}}\bigg)\\
&&\;\;\;\;\;= \mp \frac{\mu^{\frac{1-n}{2}}}{2^{n+1} \pi^{\frac{n-1}{2}}} \big( \alpha^{(n-1)/2} +n-1\big). \end{eqnarray*}
Since the {\it partition functions}  $\mathcal{Z}^{\mp}(t):=\sum_{k=1}^\infty e^{-t\tau_k^{\mp}}$  are  just
 the Riemann-Stieltjes integrals of $e^{-t\Lambda}$ with respect to the counting functions $\mathcal{N}^{\mp}(\Lambda)$,
 \begin{eqnarray} \label{2023.5.5-11} \mathcal{Z}^{\mp}(t) = \int_{-\infty}^{+\infty} e^{-t\Lambda } d\mathcal{N}^{\mp} (\Lambda),\end{eqnarray}
 one expects that the following relations hold \begin{eqnarray} \label{2023.5.5-6} \tilde{ a} = \Gamma \Big( 1+\frac{n}{2}\Big)\, a, \, \; \,\;  \tilde{b}^{\mp} = \Gamma\Big(1+\frac{n-1}{2}\Big) b^{\mp}.\end{eqnarray}

 However, by a simple calculation it can be seen that \begin{eqnarray} \label{2023.5.5-9} b^{\mp}=\frac{\tilde{b}^\mp}{\Gamma\big(1+\frac{n-1}{2}\big)} = \mp \frac{ \mu^{\frac{1-n}{2}}}{2^{n+1} \pi^{\frac{n-1}{2}}\Gamma\big(\frac{n+1}{2}\big)} \big(\alpha^{(n-1)/2} +n-1 \big), \end{eqnarray}
 which differ from the results (\ref{202.5.5-10}) and (\ref{202.5.5-3}) of \cite{CaFrLeVa-23} by additional integral terms. In other words, by applying (\ref{2023.5.5-11}) and by using the asymptotic expansion (\ref{2023.5.5-1}) it is easy to see that the second term of the heat trace asymptotic expansion in \cite{Liu-21} and that obtained from \cite{CaFrLeVa-23} have different coefficients.
 Clearly, at most one of two results in papers \cite{Liu-21} and \cite{CaFrLeVa-23} can be correct.
The authors of \cite{CaFrLeVa-23} ``predicted'' that
  \begin{eqnarray} \label{2023.5.5-12}  \lim_{\Lambda\to +\infty } \frac{\mathcal{N}^\mp(\Lambda) -a \mbox{Vol}_n (\Omega) \Lambda^{n/2} }{\mbox{Vol}_{n-1} (\partial \Omega) \Lambda^{(n-1)/2}} \end{eqnarray}
  should be the value  $b^\mp$.
 In order to ``support'' their ``prediction'', the authors of  \cite{CaFrLeVa-23} ``numerically calculated'' the values  $$\frac{\mathcal{N}^\mp(\Lambda) -a \mbox{Vol}_n (\Omega) \Lambda^{n/2} }{\mbox{Vol}_{n-1} (\partial \Omega) \Lambda^{(n-1)/2}}$$  in the interval $\Lambda\in [0,3000]$ for the two-dimensional  unit disk and flat cylinders in the interval $\Lambda\in [0,1000]$ as well as the unit square for $\Lambda$ in intervals $[0,1600]$, $[0,2400]$ and $[0, 2800]$). By comparing these  ``numerical results''  in some finite intervals, the authors of \cite{CaFrLeVa-23} ``thought (or guess)'' that their result is ``correct''.

{{\bf In this note, we shall point out that all ``numerically calculations'' and figures  in \cite{CaFrLeVa-23} are wrong because these calculations are based on some incorrect formulas.  Furthermore, by pointing out several serious errors in \cite{CaFrLeVa-23} and especially by Proposition 7.1, Remarks 7.2--7.3, and Section 8 (a result of A. Pierzchalski and B. {\O}rsted),  we show that the conclusions published by Matteo Capoferri, Leonid Friedlander, Michael Levitin and Dmitri Vassiliev \cite{CaFrLeVa-23} are completely wrong. This also implies that the most key ``algorithm'' theory for two-term spectral asymptotics in \cite{SaVa-97} is wrong.}}
 Finally, we explain the correctness of proof of Theorem 1.1 in \cite{Liu-21} by giving some remarks and putting the whole proof in Appendix. Unlike a very special elastic Lam\'{e} operator  was only dealt with on the upper-semi Euclidean space in \cite{CaFrLeVa-23}, our proof in \cite{Liu-21} had precisely investigated the corresponding elastic Lam\'{e} operator on the whole Riemannian manifold
by applying (global) geometric analysis techniques so that no information had been lost in our result.

\vskip 1.109 true cm

\section{The first serious mistake in \cite{CaFrLeVa-23} }

 \vskip 0.48 true cm

First, we point out a very obvious mistake in \cite{CaFrLeVa-23}.

\vskip 0.26 true cm

 In p.35 of \cite{CaFrLeVa-23}, the authors wrote:

\vskip 0.16 true cm

\textcolor{blue}{{\it ``More precisely, we take
$$ \qquad \quad\qquad  \quad \qquad u(r,\phi) = \mbox{grad}\;\psi_1(r,\phi) +\mbox{curl}\,  (\mathbf{z}\psi_2(r,\phi)) \qquad \;\; \;\qquad \quad \; \; \; \;\qquad  (B.1)$$
where $\mathbf{z}$ is the third coordinate vector. Then it is easily seen that the scalar potentials $\psi_j(r,\phi)$, $ j = 1, 2$,
should satisfy the Helmholtz equations
$$ \qquad \quad \quad \quad \qquad \qquad \qquad - \Delta \psi_j = \omega_{j, \Lambda}  \psi_j, \qquad \quad \qquad\qquad \qquad \qquad \;\quad \qquad \quad (B.2)$$
where
$$\qquad \qquad \quad \quad \qquad \omega_{1, \Lambda} :=\frac{\Lambda}{\lambda+2\mu},\;\; \;  \omega_{2, \Lambda}= \frac{\Lambda}{\mu}. \quad \quad \qquad \quad \,\qquad \qquad \qquad\quad \quad \;\; (B.3)$$
The general solution of (B.2) regular at the origin is well-known,
\begin{eqnarray*}\qquad  \qquad\psi_j(r, \phi)= c_{j,0}J_0(\sqrt{\omega_{j,\Lambda}}\,r)+\sum_{k=1}^\infty J_k(\sqrt{\omega_{j,\Lambda} }\, r)\big(c_{j,k,+} e^{i k\phi} + c_{j,k,-} e^{-ik\phi}\big), \;\;\quad  (B.4)\end{eqnarray*}
where the $J_k$ are Bessel functions, and the c's are constants.''
}}

\vskip 0.25 true cm

The above equations (B.2) are wrong, so that (B.4) and all results (in Appendix B: A Two-Dimensional Example: The Disk) of \cite{CaFrLeVa-23} are all wrong. The correct statement should be (cf. Theorem 2.1 below)
$$ -  \Delta (\mbox{grad}\,\psi_1) = \omega_{1, \Lambda} (\mbox{grad}\, \psi_1),\;\;\; \;\;  - \Delta \big(\mbox{curl}\, (\mathbf{z} \psi_2)\big)= \omega_{2, \Lambda} \big(\mbox{curl}\, (\mathbf{z} \psi_2)\big)$$
instead of (B.2).
That is, \begin{eqnarray*} \mbox{grad}\,\big( \Delta \psi_1 + \omega_{1, \Lambda} \psi_1\big)=0,\;\;\; \;\;  \mbox{curl}\,\big( \Delta  (\mathbf{z} \psi_2)+ \omega_{2, \Lambda} \, (\mathbf{z} \psi_2)\big)=0.\end{eqnarray*}
Or equivalently,
\begin{eqnarray} \label{2023.5.14-1} &\Delta \psi_1 + \omega_{1, \Lambda}  \psi_1=C \;\;\mbox{for any constant}\;\;C\in \mathbb{R}^1,\\
  \label{2023.5.14-2} &\Delta ( \mathbf{z} \psi_2)+ \omega_{2, \Lambda} \, (\mathbf{z} \psi_2)=\mathbf{f}, \;\,\mbox{for any} \;\; \mathbf{f}\;\; \mbox{with}\;\; \mbox{curl}\,\mathbf{f}=0.\end{eqnarray}
 In \cite{CaFrLeVa-23}, only special  $C\equiv 0$ and $\mathbf{f}\equiv 0$ are chosen. Thus, a large number of solutions have not been considered in  \cite{CaFrLeVa-23}, and a large number of elastic eigenvalues have been lost.
Obviously, for any $C\ne 0$ or $\mathbf{f}\ne 0$ with $\mbox{curl}\ \mathbf{f}=\mathbf{0}$,  the equations (\ref{2023.5.14-1}) and (\ref{2023.5.14-2})  are non-homogenous equations, their solutions  have different forms except for (B.4).
{{\bf Because  (B.2) in \cite{CaFrLeVa-23} \ is incorrect,  all calculations for elastic eigenvalues in the unit disk (in particular, (B.5) and (B.6)) in  \cite{CaFrLeVa-23} are wrong. Of course, Fig.$\,$5 and Fig.$\,$6 (on p.$\,$ 36 in \cite{CaFrLeVa-23})  are also wrong.}}

\vskip  0.25 true cm
{{\bf This mistake stems from the supplementary of an earlier paper \cite{LeMoSe-21}, in which the same mistake occurred.}}

\vskip 0.45 true cm
To help with a good understanding to the above discussions, here we copy 2.5.Theorem on p.$\,$123--124 of \cite{KGBB} and its proof:

\vskip  0.2 true cm
    \noindent{\bf Theorem 2.1.} \ {\it  Let $D$ be a domain in $\mathbb{R}^3$. The solution $\mathbf{u} = (u^1,u^2,u^3)\in [C^2(D)\cap C^1(\bar D)]^3$ of equation \begin{eqnarray} \label{2023.5.6-3}  \mu \Delta\mathbf{u}+(\lambda+\mu)\, \mbox{grad} \; \mbox{div}\, \mathbf{u}+\Lambda \mathbf{u}=0 \;\, &\mbox{in}\;\, D,\end{eqnarray} is
represented as the sum
\begin{eqnarray} \label{2023.5.6-6} \mathbf{u}=\mathbf{u}^{(p)} +\mathbf{u}^{(s)},\end{eqnarray}
where $\mathbf{u}^{(p)}$  and $\mathbf{u}^{(s)}$  are the regular vectors, satisfying the conditions
\begin{eqnarray}\label{2023.5.6-7} (\Delta +\omega_{1,\Lambda})\mathbf{u}^{(p)} =0, \;\; \mbox{curl}\; \mathbf{u}^{(p)}=0,\\
 (\Delta +\omega_{2,\Lambda})\mathbf{u}^{(s)} =0, \;\; \mbox{div}\; \mathbf{u}^{(s)}=0, \;\;\end{eqnarray}
where
\begin{eqnarray} \omega_{1,\Lambda} := \frac{\Lambda}{ \lambda+2\mu}, \;\;\, \omega_{2,\Lambda} := \frac{\Lambda}{ \mu} \end{eqnarray}
}

     \noindent {\it Proof.}  \ Let $\mathbf{u}$ be a twice differentiable function in $D \subset \mathbb{R}^3$  and
\begin{eqnarray}\label{2023.5.6-10} \left\{ \begin{array}{ll}  \mathbf{u}^{(p)}= \frac{1}{\omega_{2,\Lambda}-\omega_{1,\Lambda} }(\Delta+\omega_{2,\Lambda})\mathbf{u},\\
 \mathbf{u}^{(s)}= \frac{1}{\omega_{1,\Lambda}-\omega_{2,\Lambda} }(\Delta+\omega_{1,\Lambda})\mathbf{u}.\end{array} \right.\end{eqnarray}
Then $\mathbf{u}^{(p)} + \mathbf{u}^{(s)} =\mathbf{u}$. Let $\mathbf{u}\in C^2(D)\cap C^1(\bar D)$  be a  solution of the oscillation equation (\ref{2023.5.6-3}). The theorem will be proved if we show that
\begin{eqnarray} &&\label{2023.5.6-14} (\Delta +\omega_{1,\Lambda}) (\Delta +\omega_{2,\Lambda})\mathbf{u}=0,\\
&&\label{2023.5.6-15}(\Delta +\omega_{2,\Lambda})\,\mbox{curl}\, \mathbf{u} =0,\\
&&\label{2023.5.6-16}(\Delta +\omega_{1,\Lambda})\,\mbox{div}\, \mathbf{u} =0.
 \end{eqnarray}
Applying the operation div.$\;$to equation
(\ref{2023.5.6-3}) and taking into account the identity $\mbox{div}\,\mbox{grad} \equiv \Delta$, we obtain $(\Delta +
\omega_{1,\Lambda})\, \mbox{div}\, \mathbf{u} =0$. Similarly, the operation curl and the identity $\mbox{curl}\;\mbox{grad}\equiv  0$
lead to the equation $(\Delta + \omega_{2,\Lambda})\, \mbox{curl}\, \mathbf{u} = 0$.
Finally, applying to both sides of (\ref{2023.5.6-3}) the operator $(\Delta + \omega_{1,\Lambda})$ and taking
into account (\ref{2023.5.6-16}) which has already been proved, we obtain (\ref{2023.5.6-14}).
Thus, the regular in $D$ solution of equation (\ref{2023.5.6-3}) is represented as the sum of
the irrotational (potential) and solenoidal vectors which satisfy Helmholtz equation in $D$,
\begin{eqnarray} \label{2023.5.6-21} (\Delta + \omega)\mathbf{v} = 0,\end{eqnarray}
for $\omega = \omega_{1,\Lambda}$ and $\omega =\omega_{2,\Lambda}$, respectively. \qed

\vskip 1.19 true cm

\section{The second serious mistake in \cite{CaFrLeVa-23} }

\vskip 0.49 true cm

The similar mistake appears in Appendix C on p.$\,$36 in \cite{CaFrLeVa-23}. More precisely, (C.3) on p.$\,$36 in \cite{CaFrLeVa-23} is also wrong.  In other words, the statement \textcolor{blue}{``Once again, it is easy to see that each potential $\psi_j$ satisfies (B.2), (B.3), with....''} (see,  p.$\,$36 of \cite{CaFrLeVa-23}) is wrong. The correct expressions should be
\begin{eqnarray*}&& \Delta (\mbox{grad}\, \psi_1) +\omega_{1,\Lambda} (\mbox{grad}\, \psi_1)=0, \\
&&  \Delta (\mbox{curl}\, (\mathbf{z}\psi_2)) +\omega_{2,\Lambda} (\mbox{curl}\, (\mathbf{z}\psi_2))=0, \\
&& \Delta (\mbox{curl}\,\mbox{curl}\, (\mathbf{z}\psi_3)) +\omega_{2,\Lambda} (\mbox{curl}\,\mbox{curl}\, (\mathbf{z}\psi_3))=0.\end{eqnarray*}
 In (C.3)  of \cite{CaFrLeVa-23}, the $\psi_1$, $\psi_2$ and $\psi_3$ should be replaced by $\mbox{grad}\, \psi_1$, $\mbox{curl}\,(\mathbf{z} \psi_2)$ and $\mbox{cur}\, \mbox{curl}\, (\mathbf{z}\psi_3)$, respectively. {{\bf Clearly, all calculations on p.$\,$37 are wrong. Fig.$\,$7 and Fig.$\,$8 on p.$\,$37 in  \cite{CaFrLeVa-23} are also wrong.}}

\vskip 1.09 true cm

\section{The third serious mistake in \cite{CaFrLeVa-23} }

\vskip 0.42 true cm

 On p.$\,$32 of \cite{CaFrLeVa-23}, for $\Omega:=\{(x,y)\in \mathbb{R}^2\big|y>0\}$  and the usual Laplacian $\Delta= \partial_{xx}^2 +\partial_{yy}^2$ on $\mathbb{R}^2$, the authors of  \cite{CaFrLeVa-23} showed
\begin{eqnarray} \label{23.5.21-1} (J\circ \Delta)=(\Delta \circ J),\end{eqnarray}
 here the operations on the two sides are evaluated at the same point $(x, -y)$.
Furthermore, for the elasticity operator $\mathcal{L}$ on $\mathbb{R}^2$ the authors of  \cite{CaFrLeVa-23} showed (see p.$\,$33 of \cite{CaFrLeVa-23})
 \begin{eqnarray*}  (J\circ \mathcal{L})\ne (\mathcal{L} \circ J),\end{eqnarray*}
 here the operations on the two sides are still evaluated at the same point $(x, -y)$. The authors of  \cite{CaFrLeVa-23} then claimed  \textcolor{blue}{``The above argument shows that the principal symbol of the Laplacian (or the
Laplace-Beltrami operator when working in curved space) is invariant under reflection, whereas the principal symbol of the operator of linear elasticity is not. This is what makes the method of images work for the Laplacian, but not for the operator of
linear elasticity.''}

Actually, the authors of \cite{CaFrLeVa-23} have not understood the main ideas and the ``images method'' of  H. McKean and I. M. Singer for the Laplace-Beltrami operator on a  Riemannian manifold $(\Omega, g)$; therefore the authors of \cite{CaFrLeVa-23}  have given a series of erroneous and useless remarks as mentioned above.

 Let $\Delta_g$ (respectively, $\mathcal{L}$) be the Laplace-Beltrami operator (respectively, the elastic Lam\'{e} operator) defined on Riemannian $n$-manifold $(\Omega,g)$. The double of $\Omega$ is the manifold $\Omega \cup_{\mbox{Id}} \Omega$, where $\mbox{Id}: \partial \Omega\to \partial \Omega$ is the identity map of $\partial \Omega$; it is obtained from $\Omega \sqcup \Omega$ by identifying each boundary point in one  copy of $\Omega$
 with same boundary point in the other.
 Let  $\tau: (x',x_n)\mapsto (x', -x_n)$ be the reflection with respect to the boundary $\partial \Omega$ in $\mathcal {M}=\Omega \cup (\partial \Omega) \cup \Omega^*$ (here we always assume $x_n\ge 0$ when $(x',x_n)\in \Omega$),  where $x'=(x_1,\cdots, x_{n-1})$. Then  we can get the $\Omega^*$ from the given $\Omega$ and $\tau$.
  Put \begin{eqnarray} \label{18-2} \left\{\!\!\begin{array}{ll} \Delta_g  \;\; \;\; \mbox{in}\;\; \Omega,\\
\Delta_g \;\;\;\; \mbox{in}\;\;  \Omega^*\end{array} \right. \end{eqnarray}
  (respectively, \begin{eqnarray} \label{18-1}\left.  \left\{\!\!\begin{array}{ll} P_g   \;\; \;\;\mbox{in}\;\; \Omega,\\
P_g \;\; \;\;\mbox{in}\;\, \Omega^*\end{array} \right. \right)\end{eqnarray}
  It is easy to verify that $\Delta_g \circ J\ne J\circ \Delta_g$
 (respectively, $P_g \circ J\ne J\circ P_g$) when the values of the two sides  are all evaluated at the same point $(x', -x_n)$.
This follows from the fact that the first-order term of the Laplace-Beltrami operator $\Delta_g$ will not disappear on a (curved) Riemannian manifold, so that the full symbol $A(x, \xi_1, \cdots, x_{n-1}, \xi_n)$ of $\Delta_g$ is not an even function in $\xi_n$.
The same case still occurs when the Laplace-Beltrami operator $\Delta_g$  is replaced by the elastic Lam\'{e} operator $P_g$ on a (curved) Riemannian manifold (The reason is that each entry of the full symbols $P_g (x, \xi_1,\cdots, \xi_{n-1},\xi_n)$ of the elastic Lam\'{e} operator $P_g$ is not an even function in $\xi_n$). It is not an essential place (i.e., it is  useless at all) whether an elliptic differential operator is commutable with $J$ which are all evaluated at the same point $(x',-x_n)\in \Omega^*$ for discussing the heat trace expansion. (\ref{23.5.21-1}) holds, by chance, because the full symbol $-\xi_1^2-\xi_2^2$ of the Laplace operator $\Delta$ on $\mathbb{R}^2$  is an even function of $\xi_2$.

 The main purpose of McKean and Singer in \cite{MS-67} was to find an elliptic differential operator $A^\star$ which is defined on $\Omega^*$ such that  for all $x=(x',x_n)\in \Omega$,
 \begin{eqnarray} \label{2023.5.17-30}  \big(A^\star u(\overset{*}{x})\big)\big|_{\mbox{evaluated at the point $\overset{*}{x}$}}=  \big(\Delta_g u(\tau x) \big)\big|_{\mbox{evaluated at the point $x$,}}\end{eqnarray}
  where $\overset{*}{x}:= \tau x= (x',-x_n)$.
      In (\ref{2023.5.17-30}), the values on two sides of (\ref{2023.5.17-30}) are evaluated at two different points  $\overset{*}{x}$ and $x$ (in addition, $u(\tau x)$ is regards as a function in variable $x$), which play a key role  for considering the action of $\Delta_g$ on the second term of Green's function at the point $x$ in the last part of \cite{MS-67}).

    Let us recall McKean and Singer how to get such an elliptic operator $A^\star$ defined on $(\Omega^*,g)$, so that we can consider whether such a wonderful method may be applied to the corresponding elastic Lam\'{e} problem on the Riemannian manifold $(\Omega, g)$.
 The Laplace-Beltrami operator $\Delta_g$ defined on $(\Omega,g)$ has the following local expression:
\begin{eqnarray} && \Delta_g = \frac{1}{\sqrt{|g|}}
\sum_{j,k=1}^n \frac{\partial}{\partial x_j}\big( \sqrt{|g|}\,g^{jk} \frac{\partial}{\partial x_k}\big)\nonumber\\
  &&\label{2023.5.16}\qquad =\sum\limits_{k,l=1}^n \bigg( g^{kl} \frac{\partial^2 }{\partial x_k\partial x_l} -\sum\limits_{s=1}^n g^{kl}\Gamma_{kl}^s \frac{\partial}{\partial x_s}\bigg).\end{eqnarray}
  One  can rewrite the local expression (\ref{2023.5.16}) of  $\Delta_g$ as \begin{eqnarray} \label{2023.5.16-9}
 &&\qquad \Delta_g:= A \triangleq A  \Big( \{g^{jk} (x)\}_{1\le j,k\le n}, \{\Gamma^s_{kl}(x)\}_{1\le s,k,l\le n},  \frac{\partial }{\partial x_1}, \cdots, \frac{\partial }{\partial x_{n-1}}, \frac{\partial }{\partial x_n}\Big),\end{eqnarray}
 where $(g^{jk})=g^{-1}$ and $\Gamma^{s}_{kl}$ are the Christoffel
symbols for the Levi-Civita connection on $(\Omega,g)$.

On the $\Omega$, the Riemannian metric is still $g$. By the reflection operator $\tau$,  McKean-Singer gave the metric on $\Omega^*$ to be
\begin{eqnarray} \label{2023.5.17-18} \left\{ \begin{array}{ll}   g_{jk} (\overset{*}{x})=- g_{jk} (x) \quad \, \mbox{for}\;\;
  j<k=n \;\;\mbox{or}\;\; k<j=n,\\  g(\overset{*}{x}) = g_{jk} (x)\;\; \;\;\mbox{for}\;\; j,k<n \;\;\mbox{or}\;\; j=k=n,\\
    g_{jk}(x)= 0 \;\; \;\mbox{for}\;\; j<k=n \;\;\mbox{or}\;\; k<j=n \;\;\mbox{on}\;\; \partial \Omega.\end{array}\right. \end{eqnarray}  We can denote the latter as $(\Omega^*, g^*)$.
In this way (pasting $(\Omega,g)$ and $(\Omega^*,g^*)$ together by the identity map of $\partial \Omega$), McKean-Singer obtained a metric tensor on the whole $\mathcal{M}=\Omega
\cup (\partial \Omega)\cup \Omega^*$. Clearly, $\mathcal{M}$ is a closed Riemannian manifold with continuous metric tensor. It is easy to verify that
\begin{eqnarray}\label{2023.5.16-5} \;\;\,\; \left. \begin{array}{ll} \Gamma^{j}_{kl}(\overset{*}{x})=\Gamma^{j}_{kl}(x) \,\,\; \,\mbox{for}\; \; 1\le j,k,l<n;\\
 \Gamma^{j}_{kl}(\overset{*}{x})=-\Gamma^{j}_{kl}(x) \;\;\mbox{for}\; \; 1\!\le \! j,k\!<\!n, \;\, l\!=\!n; \;\; \mbox{or} \;\, 1\!\le \!j,l\!<\!n, \ k\!=\!n, \;\; \mbox{or}\;\, 1\!\le \!k,l\!<\!n,\;  j\!=\!n;\\
  \Gamma^{j}_{kl}(\overset{*}{x})=\Gamma^{j}_{kl}(x) \;\;\mbox{for}\, \; 1\!\le\!  j\!<\!n, \;\, k\!=\!l\!=\!n; \; \mbox{or} \,\, 1\!\le \! k<n, \, j\!=\!l\!=\!n, \,\; \mbox{or}\;\, 1\!\le\! l\!<\!n,\;  j\!=\!k\!=\!n;\\
   \Gamma^{n}_{nn}(\overset{*}{x})=-\Gamma^{n}_{nn}(x). \end{array}\right.
 \end{eqnarray}
More precisely,  McKean and Singer in \cite{MS-67} wanted to construct a new elliptic differential operator $A^\star$ defined on $\Omega^*$ such that the following key relation holds:

\vskip 0.25 true cm 

 \begin{eqnarray*} && \mbox{$w(x):=u(\tau x)$ (regarded as a function of variable $x$) being acted by the} \\
   && \mbox{ differential operator $A$ at the point $x\in \Omega$ just is equal to   $u(\overset{*}{x}):=u(\tau x)$ (regarded}\\
&&\mbox{ as a function of variable $\overset{*}{x}$) being acted by $A^\star$ at the point $\overset{*}{x}=\tau x$}.\end{eqnarray*}

\vskip 0.30 true cm 

 Equivalently, the above  $A^\star$ can also be got by another way as follow: Let $u$ be a function defined on $\Omega^*$. Pull back the function $u$ by the reflection operator  $\tau$ to get the function $w$
(defined on $\Omega$), i.e., 
    \begin{eqnarray*} w(x):=u(\tau x)=u(\overset{*}{x})\;\;\mbox{for}\;\; x\in \Omega,\end{eqnarray*} then define $A^\star$ by
  \begin{eqnarray*} && A^\star u(\overset{*}{x} )\big|_{\mbox{evaluated at the point $\overset{*}{x}$}}: =A(w(x))\big|_{\mbox{evaluated at the point $x$}}\\
  && \quad \quad \quad \;\;\quad =\Delta_g w(x)\big|_{\mbox{evaluated at the point $x$}}\,\quad \;\; \; \;\; \mbox{for all} \;\, \overset{*}{x} \in \Omega^*.\end{eqnarray*}

By some calculations, it can be seen that
\begin{eqnarray}\label{2023.5.14-17}  A^\star :=\Big( \{g^{jk} (x)\}_{1\le j,k\le n}, \{\Gamma^s_{kl}(x)\}_{1\le s,k,l\le n},  \frac{\partial }{\partial x_1}, \cdots, \frac{\partial }{\partial x_{n-1}}, -\frac{\partial }{\partial x_n}\Big).\end{eqnarray}
But we must rewrite the local expression (\ref{2023.5.14-17}) in the language of metric $g^*:=g(\overset{*}{x})$ and Christoffel symbols  $(\Gamma^{j}_{kl})^*=(\Gamma^{j}_{kl}(\overset{*}{x}))$ on Riemannian manifold $(\Omega^*,g^*)$. Substituting (\ref{2023.5.17-18}) and (\ref{2023.5.16-5}) into
(\ref{2023.5.14-17}), one can obtain the exact expression of $A^\star$ on $(\Omega, g^*)$:
\begin{eqnarray}\label{2023.5.14-19} && A^\star :=\Big( \{g^{jk} (\overset{*}{x})\}_{1\le j,k< n}, \{-g^{jn} (\overset{*}{x})\}_{1\le j< n},  \{-g^{nk} (\overset{*}{x})\}_{1\le k< n}, \{\Gamma^s_{kl}(\overset{*}{x})\}_{1\le s,k,l< n},  \\
&& \qquad\quad \{-\Gamma^n_{kl}(\overset{*}{x})\}_{1\le k,l< n},
 \{-\Gamma^s_{nl}(\overset{*}{x})\}_{1\le s,l< n},  \{-\Gamma^s_{kn}(\overset{*}{x})\}_{1\le s,k< n},  \{\Gamma^s_{nn}(\overset{*}{x})\}_{1\le s< n},
 \nonumber \\
&& \qquad \quad   \{\Gamma^n_{kn}(\overset{*}{x})\}_{1\le k< n},
 \{\Gamma^n_{nl}(\overset{*}{x})\}_{1\le l< n},\,-\Gamma^n_{nn}(\overset{*}{x}),
 \frac{\partial }{\partial x_1}, \cdots, \frac{\partial }{\partial x_{n-1}}, -\frac{\partial }{\partial x_n}\Big).\nonumber \end{eqnarray}
McKean and Singer further defined  \begin{eqnarray} \label{2022.10.18-22} \mathcal{A}=\left\{\begin{array}{ll} \! A \;\;\; \;\;\;\,\mbox{on} \;\, (\Omega, g)\\
 \!A^\star \;\;\;\; \;\mbox{on} \;\, (\Omega^*,g^*). \end{array} \right.\end{eqnarray}
By the definitions of $A^\star$ and $\mathcal{A}$, it immediately follows that for any $x\in \Omega$,
 \begin{eqnarray} \label{2023.5.17-20}  \\ \Delta_g \big(w(x)\big)\big|_{\mbox{evaluated at the point $x$}}\!=\! A^{\star}
\big( u(\overset{*}{x})\big)\big|_{\mbox{evaluated at the point $\overset{*}{x}$,}} \;\,\;\;\,\,\;\,\; w(x)= u(\tau(x)), \nonumber  \end{eqnarray}
 which is the most key requirement when calculating  the second term $\big(\Delta_g K(t, \tilde{x}, \overset{*}{x})\big)\big|_{\mbox{evaluated at the point $x$}}$ for the heat kernel Green function of  \begin{eqnarray} \label{2023.5.17-21}&& \big(\Delta_g K^{\mp}K(t,\tilde{x},x)\big)\big|_{\mbox{evaluated at the point $x$}} \\
 && \qquad \quad =\big(\Delta_g \big( (K(t, \tilde{x}, x)-K(t, \tilde{x}, \overset{*}{x})\big)\big)\big)\big|_{\mbox{evaluated at the point $x$},}\nonumber\end{eqnarray}
 here $\tilde{x}\in \Omega$ and the Laplace-Beltrami operator $\Delta_g$ is acted to variable $x$.
 Obviously,  $\mathcal{A}$ is still a linear elliptic differential operator on whole $\mathcal{M}$. It can be seen that the coefficients of the differential operator $\mathcal{A}$ are not smooth on whole $\mathcal{M}$. But the linear elliptic differential operator $\mathcal{A}$ is not ``too bad'' on $\mathcal{M}$ because the top-order coefficients of $\mathcal{A}$ are continuous  on the whole $\mathcal{M}$, the lower-order coefficients are bounded measurable on the whole $\mathcal{M}$ (they are discontinuous when $x$ crosses $\partial \Omega$, see the statement on  p.$\,$53
 of \cite{MS-67}), all coefficients of $\mathcal{A}$ are smooth in $\mathcal{M}\setminus (\partial \Omega)$. These properties of the operator $\mathcal{A}$ are  enough for studying the $W^{2,p}$-estimate and  regularity  of  solutions of the heat equation.
 
Obviously, the authors of \cite{CaFrLeVa-23} did not even know what is  McKean-Singer's ``double operator''. In \cite{CaFrLeVa-23} and  \cite{CaFrLeVa-22}, they had erroneously regarded the ``double operator'' of Laplacian $\Delta_g$ as the following operator 
\begin{eqnarray*} \left\{ \begin{array}{ll} \Delta_g \;\;\;\mbox{in}\;\, \Omega, \\
 \Delta_g \;\;\;\mbox{in}\;\; \Omega^*,\end{array}\right.\end{eqnarray*} 
so that they gave many irrelevant (and useless) remarks. For example, from line 3 to line 5 on p.$\,$34 in \cite{CaFrLeVa-23})  they wrote:  ``\textcolor{blue}{However, McKean and Singer
applied the method of images to the Laplacian, for which the double operator is
self-adjoint.}''\  {\bf In fact, as being discussed above,  the ``double'' operator $\mathcal{A}$  on double manifold $\mathcal{M}$, which was considered by McKean and Singer in  \cite{MS-67}, is given by   (\ref{2023.5.16-9}), (\ref{2023.5.14-19}) and  (\ref{2022.10.18-22}).} 

\vskip 0.12 true cm 

On p.$\,$10177 in \cite{Liu-21},  we have expressed the Navier-Lam\'{e} operator $P_g$ in $\Omega$  as the form of components relative to local coordinates:
 \begin{eqnarray} \label{2023.5.22-12} \\
 P_g\mathbf{u}=\!\!\!\!\!\!&&\!\!\!\left\{ -\mu\Big( \sum_{m,l=1}^n g^{ml} \frac{\partial^2 }{\partial x_m\partial x_l}\Big){\mathbf{I}}_n -(\mu+\lambda)
\begin{bmatrix} \sum\limits_{m=1}^n g^{1m} \frac{\partial^2}{ \partial x_m\partial x_1} & \cdots & \sum\limits_{m=1}^n g^{1m} \frac{\partial^2}{ \partial x_m\partial x_n}
  \\    \vdots & {} & \vdots  \\
 \sum\limits_{m=1}^n g^{nm} \frac{\partial^2}{ \partial x_m\partial x_1} & \cdots & \sum\limits_{m=1}^n g^{nm} \frac{\partial^2}{ \partial x_m\partial x_n}   \end{bmatrix} \right.\nonumber\\
&& +\mu \Big( \sum\limits_{m,l,s=1}^n g^{ml} \Gamma_{ml}^s \frac{\partial}{\partial x_s} \Big){\mathbf{I}}_n - \mu  \begin{bmatrix}  \sum\limits_{m,l=1}^n 2g^{ml} \Gamma_{1m}^1 \frac{\partial }{\partial x_l} & \cdots &  \sum\limits_{m,l=1}^n 2g^{ml} \Gamma_{nm}^1 \frac{\partial }{\partial x_l}\\
\vdots & {}& \vdots \\
\sum\limits_{m,l=1}^n 2g^{ml} \Gamma_{1m}^n \frac{\partial }{\partial x_l} & \cdots & \sum\limits_{m,l=1}^n 2g^{ml} \Gamma_{nm}^n \frac{\partial }{\partial x_l}
\end{bmatrix}\nonumber \\
&& -(\mu+\lambda) \begin{bmatrix}\sum\limits_{m,l=1}^n g^{1m} \Gamma_{1l}^l  \frac{\partial }{\partial x_m} & \cdots & \sum\limits_{m,l=1}^n g^{1m} \Gamma_{nl}^l  \frac{\partial }{\partial x_m} \\       \vdots & {} & \vdots  \\
  \sum\limits_{m,l=1}^n g^{nm} \Gamma_{1l}^l  \frac{\partial }{\partial x_m} & \cdots & \sum\limits_{m,l=1}^n g^{nm} \Gamma_{nl}^l  \frac{\partial }{\partial x_m} \end{bmatrix}\nonumber
 \\
 &&  -\mu \begin{bmatrix} \sum\limits_{l,m=1}^n g^{ml} \big(  \frac{\partial \Gamma^1_{1l}}{\partial x_m} +  \Gamma_{hl}^1 \Gamma_{1m}^h - \Gamma_{1h}^1 \Gamma_{ml}^h \big) & \cdots &  \sum\limits_{l,m=1}^n g^{ml} \big(  \frac{\partial \Gamma^1_{nl}}{\partial x_m} +  \Gamma_{hl}^1 \Gamma_{nm}^h - \Gamma_{nh}^1 \Gamma_{ml}^h \big)
\\ \vdots & {} & \vdots \\
\sum\limits_{l,m=1}^n g^{ml} \big(  \frac{\partial \Gamma^n_{1l}}{\partial x_m} +  \Gamma_{hl}^n \Gamma_{1m}^h - \Gamma_{1h}^n \Gamma_{ml}^h \big) & \cdots &  \sum\limits_{l,m=1}^n g^{ml} \big(  \frac{\partial \Gamma^n_{nl}}{\partial x_m} +  \Gamma_{hl}^n \Gamma_{nm}^h - \Gamma_{nh}^n \Gamma_{ml}^h \big)
\end{bmatrix} \nonumber\\
&& \left.- (\mu+\lambda) \begin{bmatrix} \sum_{l,m=1}^n g^{1m} \frac{\partial \Gamma^l_{1l}}{\partial x_m}& \cdots & \sum_{l,m=1}^n g^{1m} \frac{\partial \Gamma^l_{nl}}{\partial x_m}\\
\vdots & {} & \vdots \\
\sum_{l,m=1}^n g^{nm} \frac{\partial \Gamma^l_{1l}}{\partial x_m}& \cdots & \sum_{l,m=1}^ng^{nm} \frac{\partial \Gamma^l_{nl}}{\partial x_m}\end{bmatrix} - \mu \begin{bmatrix}
R^1_1 & \cdots & R^1_n\\
\vdots & {} & \vdots \\
R^n_1 & \cdots & R^n_n\end{bmatrix} \right\}\begin{bmatrix} u^1\\
\vdots\\
u^n\end{bmatrix},  \nonumber \end{eqnarray}
where ${\mathbf{I}}_n$ is the $n\times n$ identity matrix.
The above expression $P_g \mathbf{u}$ has played a key role in our paper \cite{Liu-21}.
From this, we can also find an elliptic differential operator $P^\star$ defined on $\Omega^*$ such that the following relation holds \begin{eqnarray}
 \label{23.5.22-19} \\
  P_g \big(w(x)\big)\big|_{\mbox{evaluated at the point $x$}} =
 P^{\star}
\big( u(\overset{*}{x})\big)\big|_{\mbox{evaluated at the point $\overset{*}{x}$\,,}} \quad \quad\;\;  w(x):=u(\tau (x)).  \nonumber\end{eqnarray}
This implies that all works can still be done for the elastic Lam\'{e} operator $P_g$ on a Riemannian manifold by overcoming some other difficulties, and the two coefficients of the asymptotic expansion can be obtained (see \cite{Liu-21}). { {\bf  Obviously, the authors of \cite{CaFrLeVa-23} have not understood the key step (i.e., (\ref{2023.5.17-20}) and  (\ref{23.5.22-19}))  for seeking a new (important elliptic differential) operator both for the Laplace-Beltrami operator $\Delta_g$ and for the elastic Lam\'{e} operator $P_g$  on the reflection manifold $(\Omega^*,g^*)$.}}

\vskip 1.12 true cm

\section{The forth serious mistake in \cite{CaFrLeVa-23} }

\vskip 0.36 true cm

  In order to ``overthrow'' the conclusion of \cite{Liu-21}. In Appendix A of \cite{CaFrLeVa-23}, authors wrote (see, line -19 to line -1 from bottom on  p.33):
\textcolor{blue}{``The author defines $\mathcal{M}:=\Omega\cup \Omega^*$ to be the ``double'' of
$\Omega$, and $\mathcal{T}$ to be
 the ``double'' of the operator of linear elasticity in $\mathcal{M}$.
In the simplified setting of this appendix $\mathcal{M}:=\mathbb{R}^2$ and
\begin{eqnarray*}
\mathcal{T}:=\left\{\begin{array}{ll} \mathcal{L}, \;\; &\mbox{in}\;\; \Omega= \{y>0\},\\
\mathcal{L}+2(\lambda+\mu) \begin{pmatrix}0& \partial_x\\
\partial_x&0  \end{pmatrix}\frac{\partial }{\partial y} \;\; & \mbox{in}\;\; \Omega^*=\{y<0\}.\end{array}\right.
\end{eqnarray*}
Given $u, v\in C^\infty_c (\mathbb{R}^2)$ (the space of infinitely smooth functions with compact support),
by a straightforward integration by parts, one obtains
\begin{eqnarray*}\quad \quad \quad  ( \mathcal{T}\mathbf{u}, \mathbf{v})-(\mathbf{u}, \mathcal{T}\mathbf{v})=4(\lambda+\mu)\int_{\{y=0\} } (\partial_x u_1\, \bar{v}_2 -\mathbf{u}_2 \overline{\partial_x \mathbf{v}_1})dx. \;\;\;\;\quad \;  (A.1) \end{eqnarray*}
Here $(\cdot, \cdot)$ denotes the natural $L^2$ inner product and overline denotes complex conjugation. But (A.1) implies that $\mathcal{T}$ is not symmetric; therefore, it does not give rise to
a heat operator.
As a result, the statement ``$\mathcal{T}$ generates a strongly continuous semigroup $(e^{t\mathcal{T}})_{t\ge 0} $ on $L^2(\mathcal{M})$ with integral kernel $\mathbf{K}(t, x, y)$'' in [13, p. 10169, third line after (1.14)] is
wrong, and all the analysis based on it breaks down (including [13, formula (4.3)]).''}

\vskip 0.25 true cm

{\bf The fundamental mistake in \cite{CaFrLeVa-23} is the above statement} \textcolor{blue}{``But (A.1) implies that $\mathcal{T}$ is not symmetric; therefore, it does not give rise to a heat operator.''}
{{\bf Such a statement (for an elliptic operator) is not true!}} For example, let  the (non-divergence form) differential operator $L$ be defined on an open set $\Omega$ in $\mathbb{R}^n$:
\begin{eqnarray*} L=\sum_{|\alpha|\le 2m} a_\alpha(x) D_\alpha, \;\;\;\; D_\alpha=\frac{\partial^{\alpha_1}}{\partial x_1^{\alpha_1}}\,\frac{\partial^{\alpha_2}}{\partial x_1^{\alpha_2}}\,\cdots \,\frac{\partial^{\alpha_n}}{\partial x_1^{\alpha_n}},\; \;\; x\in \Omega \subset \mathbb{R}^n,\end{eqnarray*}
where $L$ is uniformly strongly elliptic, $L$ satisfies the  {\it root condition}, $a_\alpha$  are measurable and uniformly bounded on $\bar \Omega$ by a positive constant $B$, the $a_\alpha$ are uniformly continuous on $\Omega$ for $|\alpha|=2m$. For clarity, we also assume that $a_\alpha(x)$ ($|\alpha|<2m$) are discontinuous on $\Omega$ (for example, $a_\alpha(x)$ jump when $x$ crosses an $(n-1)$-dimensional  hyper-surface $\mathcal{S}$ in $\Omega$). Obviously, for any $u, v\in C^\infty_c (\Omega)$,  \begin{eqnarray} \label{2023.5.15-1}(Lu, v)-(u, Lv)\ne 0,\end{eqnarray} which implies that $L$ is not symmetric, where  $(\cdot, \cdot)$ still denotes the natural $L^2$ inner product. However, it is well-known that $L$ can give rise to
a heat operator; furthermore, $L$ generates strongly continuous semigroups (actually, $L$ can generate analytic semigroups) $(e^{tL})_{t\ge 0}$ with integral kernel $K(t, x,y)$ on $L^2(\Omega)$ space (by Browder \cite{Brow}), $L^p(\Omega)$ space (by Friedman \cite{Fri}), $C_{00}(\bar \Omega)$ space (by Stewart \cite{Ste} if all $a_\alpha$ are continuous, $\Omega$ may be an unbounded domain), respectively. These celebrated results were based on $L^{2,p}$-estimates for such kinds of elliptic differential operators  which have continuous top-order coefficients $a_{\alpha}$ for all $|\alpha|=2m$ (Of course, $a_\alpha$, $|\alpha|<2m$, may be bounded measurable).
These classical $W^{2,p}$-estimates and regularity theory for elliptic operators and parabolic operators had been established in the last century (cf. \cite{ADN}, \cite{GiTr} or \cite{Fri}).

Now, let us go back to the double differential operator $\mathcal{A}$ on the double manifold $\mathcal{M}$ (see the previous section), which is deduced by the Laplace-Beltrami operator $\Delta_g$ on Riemannian manifold $(\Omega, g)$ in \cite{MS-67}. Clearly, $\mathcal{A}$ is an elliptic differential operator on $\mathcal{M}$ whose coefficients are smooth in $\mathcal{M}\setminus (\partial \Omega)$. Note that the endowed metric tensor on $\mathcal{M}$ are continuous. By the local expression of $\mathcal{A}$, we see that the top-order coefficients of $\mathcal{A}$ are continuous on the closed Riemannian manifold $\mathcal{M}$, but the lower-order coefficients of $\mathcal{A}$ are discontinuous on $\partial \Omega$ (although they are bounded measurable on whole  $\mathcal{M}$). Therefore, by applying  the $W^{2,p}$-estimates and regularity theory that are mentioned above, McKean and Singer in \cite{MS-67} discussed the  parabolic equation $u_t- \mathcal{A}u=0$ on $\mathcal{M}$ and considered the strongly continuous semigroups $(e^{t \mathcal{A}})_{t\ge 0}$ as well as the integral kernel on $\mathcal{M}$ so that they obtained the heat asymptotic expansion for Laplace-Beltrami operator $\Delta_g$.

\vskip 0.25 true cm

Next, we discuss  the elastic Lam\'{e} operator as follows. Fortunately, $C^{1,1}_{loc}$-estimates and $W^{2,p}_{loc}$-estimates of viscosity solutions have been given  by J. Xiong \cite{Xi-11} for non-divergence form elliptic differential equations, and $W^{1,2}_2$-estimates and $C^{1,1}_{loc}$-estimates of strong solutions have been  established by H. Dong \cite{Do-12}  for non-divergence parabolic differential equations which have two pieces smooth coefficients $a_\alpha$ (i.e., $a_\alpha$ jump only on a  smooth hypersurfaces). In \cite{Do-12}, $W^{1,2}_p:= \{ \mathbf{u}\big| \mathbf{u},\frac{\partial \mathbf{u}}{\partial t}, D\mathbf{u}, D^2 \mathbf{u} \in L^p\} $.
This type of system arises from the problems of linearly elastic laminates and composite materials (see, for example, \cite{CKVC}, \cite{LiVo}, \cite{LiNi}, \cite{Do-12} and \cite{Xi-11}).
 They are  important and breakthrough works on  regularity theory for non-divergence  elliptic equations and non-divergence form  parabolic differential equations with discontinuous coefficients only occurred in  both sides of a smooth hypersurface in this century.
In fact, for scalar non-divergence form parabolic equations
$Mu :=-u_t + \sum_{j,k=1}^n a_{jk}D_{jk} u + \sum_{j=1}^n b_j D_j u + cu = f$,  H. Dong  proved (p.121 of \cite{Do-12}) that  if the coefficients and $f$ are Dini continuous in $x'=(x_1, \cdots, x_{n-1})$ (irregular only in one spatial direction $x_n$), then any solution
$u$ of the above equation is $C^1$ at $t$, $C^{1,1}$ in $x$, and $u_t$ and $D_{xx'}u$ are continuous.
Obviously, Dong's results (see, Theorem 4 on p.$\,$141 of \cite{Do-12}) still hold for such kinds of (special) non-divergence form  (elastic) parabolic systems  whose coefficients are  bounded smooth on both sides of a smooth hypersurface (but jump as crossing this smooth hypersurface).

  In order to further explain the solutions for non-divergence form elliptic equations, we need to recall the concept of viscosity solutions which was motivated by L. Caffarelli for studying the regularity of non-linear elliptic equation(s) (see, for example, \cite{CaCa}):

\vskip 0.25 true cm

\noindent{\bf Definition 5.1.} \   For $f\in L^p_{loc} (D), p>n/2$, a function $u\in C(\Omega)$ is an $L^p$-viscosity subsolution (supersolution) of
$Lu = f$ in $D$ if for any $\phi\in W^{2,p}_{loc}(\Omega)$ touching
$u$ from above (resp. below) at point $\tilde{x}$ locally one has
\begin{eqnarray*}&& \mbox{ess}\, \lim\limits_{x\to \tilde{x}}\sup  \,\Big(\sum_{|\alpha|\le 2} a_\alpha (x) D_\alpha \phi_\alpha (x)- f(x)\Big)\ge 0,\\
&&\big( \mbox{ess}\, \lim \limits_{x\to \tilde{x}}\inf \, \Big(\sum_{|\alpha|\le 2} a_\alpha (x) D_\alpha \phi_\alpha (x)- f(x)\Big)\le 0\big).
\end{eqnarray*}
We say $u \in C(D)$ is an $L^p$-viscosity of $Lu = f$ if $u$ is both an $L^p$-viscosity subsolution
and supersolution.

\vskip 0.25 true cm

\noindent{\bf Theorem 5.2 (see p.375 of \cite{Xi-11}):}
{\it Let $D$ is a bounded domain in $\subset \mathbb{R}^n$, and let $\mathcal{S}\subset \mathbb{R}^n$ be an $(n-1)$-dimensional
embedded (but not necessarily connected or compact) $C^{1,\alpha}$ hypersurface for some $\alpha\in (0,1)$. Suppose that
$\mathcal{S}\cap D\ne \varnothing$  and for any point $x\in \Omega$ there exists a positive constant $r$, depending on $x$,
such that $a_\alpha$ are uniformly H\"{o}lder continuous on every connected component
of $B_r\setminus \mathcal{S}$ but might be discontinuous cross $\mathcal{S}$.  Then for every boundary value
$\phi \in C(\partial D)$, there exists a unique $L^p$-viscosity solution $u \in C^{1,1}_{loc} (D) \cap C(D)$  of the
Dirichlet problem
$$Lu=f \;\, \mbox{in}\;\;  D, \;\;\, u=\phi\;\;\mbox{on}\;\;\partial D.$$}

\vskip 0.25 true cm

A more important result for parabolic equation is the following:

\vskip 0.12 true cm

\noindent{\bf Theorem 5.3 (see Theorem 4 on  p.141 of \cite{Do-12}):} \  
{\it   Let $\delta\in (0,1)$, $a,b,c\in C_{z'}^{\delta/2, \delta}$ and $f\in C_{z'}^{\delta/2, \delta} (Q_1)$. Assume that $u\in W^{1,2}_2(Q_1)$ is a strong solution in $Q_1$ to 
\begin{eqnarray} \label{2023.6.14-9} Pu:= -u_t + a^{\alpha\beta} D_{\alpha\beta} u+b^\alpha D_\alpha u+cu=f.\end{eqnarray} 
Then we have 
\begin{eqnarray*} |u|_{1,2;Q_{1/2}} +[u_t]_{\delta/2, \delta; Q_{1/2}} +[D_{xx'}u]_{\delta/2,\delta;Q_{1/2}} \le N\big(|f|_{z',\delta/2,\delta;Q_{1}} +\|u\|_{L^2(Q_1)}\big) \end{eqnarray*}
where $Q_r=Q_r(0,0)$ and $Q_r(t,x)=(t-r^2,t)\times B_r (x)$, $\,N=N(d,\delta, \nu, K,[a]_{z',\delta/2, \delta}, [b]_{z',\delta/2,\delta}),  [c]_{z',\delta/2,\delta})$. }
\vskip 0.10 true cm

  By applying Dong's $W^{1,2}_2$-estimates and $C^{1,1}$-estimates (see, from line -12 to line -4 from the bottom on p.$\,$121 of \cite{Do-12}) 
to our case of elastic Lam\'{e} operator, we can immediately obtain the corresponding parabolic equation  system and strongly continuous semigroups in $C^0$ and $L^2$ spaces, respectively.

\vskip 0.25 true cm

 \  Let us come back to previous discussion for the elastic Lam\'{e} operator defined on $\mathbb{R}^n_+$ once again.  The  new constructed (double Lam\'{e}) operator $\mathcal{P}$ is not symmetric  since $\mathcal{P}$ is a non-divergence form elliptic operator with discontinuous coefficients in $\mathcal{M}$. {\bf However, the heat operator corresponding to $\mathcal{P}$ can be given as follows}.

The authors of \cite{CaFrLeVa-23} had not understood our method in \cite{Liu-21} and the strong solutions (see, the definition in \cite{Do-12}) for the (non-divergence form) parabolic equation. Let us introduce the correct approach in \cite{Liu-21} (also see \cite{Liu-22c}). We still put the discussion on the semi-space $\mathbb{R}^2_+$ as done in \cite{CaFrLeVa-23}, although this is not a suitable domain because $\mathbb{\bar R}^2_+$ is a flat Euclidean unbounded domain (It should be discussed in a subdomain in $\mathbb{\bar R}^2_+$ which intersects the hyperplane $\{(x,y)\in \mathbb{R}^2\big|y=0\}$).
 Let $\mathcal{L}$ be the elasticity operator on $\mathbb{\bar R}^2_+:=\{y\ge 0\}$, which acts on vector-valued functions
$\mathbf{u}(x,y)=\begin{pmatrix} u_1(x,y)\\
u_2(x,y)\end{pmatrix}$ for $(x,y)\in \mathbb{\bar R}^2_+$  as
\begin{eqnarray*} \mathcal{L}\mathbf{u}:= -\mu \Delta \mathbf{u}-(\lambda+\mu) \mbox{grad}\, \mbox{div}\,\mathbf{u}=\begin{pmatrix} -\mu \Delta -(\lambda+\mu) \partial_{xx} & -(\lambda+\mu) \partial_{xy}\\
-(\lambda+\mu) \partial_{xy}&  -\mu \Delta -(\lambda+\mu) \partial_{yy}    \end{pmatrix} \begin{pmatrix} u_1 \\u_2\end{pmatrix} . \end{eqnarray*}
   From the above $\mathcal{L}$ defined only on $\mathbb{R}^2_+$, we further introduce the operator defined on $\mathbb{R}^2_-:=\{(x,y)\in \mathbb{R}^2\big|y<0\}$ by
   \begin{eqnarray*} P^\star:=-\mu \Delta \begin{pmatrix} u_1\\
u_2\end{pmatrix} -(\lambda+\mu) \begin{pmatrix} \partial_{xx} & -\partial_{xy}\\
-\partial_{xy} & \partial_{yy} \end{pmatrix}\begin{pmatrix} u_1\\
u_2\end{pmatrix}.\end{eqnarray*}
   Put    $\mathcal{T}$ in $\mathbb{R}^2:=\mathbb{\bar R}^2_{+}\cup
 \mathbb{R}^2_{-}$ by
\begin{eqnarray*} &&\mathcal{T}:=\left\{\begin{array}{ll}\mathcal{L}
\;\; &\mbox{on}\;\; \mathbb{\bar R}^2_+,\\
 P^\star \;\;& \mbox{on}\;\;  \mathbb{R}^2_-.\end{array}\right.
\end{eqnarray*}
Let $\tau : (x, y)\mapsto  (x, -y)$ for $(x, y)\in \mathbb{R}^2_+$ be a reflection with respect to the $x$-axis from $\mathbb{R}^2_+$ to $\mathbb{R}^2_-$.
 Let us point out that $\tau$ has definition only on $\mathbb{R}^2_+$.
 For a vector-valued function
$\mathbf{u}$  consider the involution $J\mathbf{u} := \mathbf{u} \circ \tau$ , so that $(J\mathbf{u})(x, y) =
\mathbf{u}(x, -y)$ for any $(x, y)\in \mathbb{R}^2_+$.
 It can be easily verified that

\begin{eqnarray*} && \mbox{by putting $\mathbf{w}(x,y):=\mathbf{u}(x,-y)$, the vector-valued function $\mathbf{w}(x,y)$}\\
 &&\mbox{being acted by $\mathcal{L}$ at the point $(x,y)$  equals to the vector-valued }\\
 && \mbox{function $\mathbf{u}(x, -y)$ being acted by $P^\star$ at the point $(x,-y)$}.\end{eqnarray*}

\vskip 0.25 true cm 

 {\bf Without doubt, the authors of \cite{CaFrLeVa-23} had not understood the key steps and important idea in \cite{Liu-21} that were mentioned above. They erroneously regarded ``double operator'' $\mathcal{T}$ as
\begin{eqnarray*} \mathcal{\tilde{T}}= \left\{ \begin{array}{ll} \mathcal{L} \;\; \mbox{in}\;\; \mathbb{R}^2_+,\\
\mathcal{L} \;\; \mbox{in}\;\; \mathbb{R}^2_-,\end{array} \right.\end{eqnarray*}
and discussed the commutativity for $\mathcal {\tilde{T}}$ and $J$ (that is useless at all).}

\vskip 0.19 true cm 

Clearly, $\mathcal{T}$ is a (non-divergence form) elliptic operator with bounded measurable coefficients in $\mathbb{R}^2$.
The coefficients of $\mathcal{T}$ are  discontinuous in $\mathbb{R}^2$ (jump on $x$-axis $\partial \mathbb{R}^2_+$), but the coefficients of $\mathcal{T}|_{\mathbb{\bar R}^2_+}$ and $\mathcal{T}|_{\mathbb{R}^2_-}$ are smooth on $\mathbb{\bar R}^2_+$ and (up to) $\mathbb{\bar R}^2_-$, respectively.
For such kind of elliptic operators (their coefficients are discontinuous only on an $(n-1)$-dimensional smooth hyper-surface), J. Xiong in \cite{Xi-11} established the existence and regularity of solutions (in the sense of viscosity solution).
For more general (non-divergence form) elliptic operators $A$ whose coefficient are discontinuous only at one direction,  H. Dong \cite{Do-12} established the regularity theory of  solutions for the parabolic equations $(\frac{\partial}{\partial t}-A)u=f$ (in sense of strong solution). H. Dong in \cite{Do-12} also proved that the solutions are
$C^{1,1}_{loc}$, and gave $W^{1,2}_2$-estimates (In fact, H. Dong's  results are enough for our discussing in regularity of solutions).  
Thus, ``\textcolor{blue}{$\mathcal{P}$ ...., it does not give rise to
a heat operator. ''  in line -4 on p.$\,$33 of \cite{CaFrLeVa-22}} is a wrong statement.

\vskip 0.2 true cm
{{\bf Because the authors of \cite{CaFrLeVa-23} have not known these great progresses in regularity of non-divergence form parabolic (or elliptic) equations with discontinuous points of coefficients only placed on a smooth hypersurface, authors of \cite{CaFrLeVa-23}  gave a very ridiculous and  erroneous statement}}: \textcolor{blue}{``But (A.1) implies that $\mathcal{T}$ is not symmetric; therefore, it does not give rise to a heat operator.''} (see, from line -5 to line -3 on p.33 in \cite{CaFrLeVa-23}).

\vskip 1.19 true cm

\section{The fifth serious mistake in \cite{CaFrLeVa-23} }

\vskip 0.59 true cm

\noindent{\bf Remark 6.1.} \   When authors of \cite{CaFrLeVa-23} can not clearly determine which of two conclusions in \cite{Liu-21} and  \cite{CaFrLeVa-23} are correct, they used a non-profesional (ridiculous)  method, so-called the ``numerical verification''. That is,   ``predict'' a limit value $\lim_{\Lambda\to +\infty} \frac{1}{\mbox{Vol}_{n-1}(\partial \Omega)\,\Lambda^{(n-1)/2}}\big[ \mathcal{N} (\Lambda)- a \mbox{Vol}_n (\Omega) \Lambda^{n/2}\big]$  by only calculating the value in a finite interval $[0,3000]$ for $\Lambda$. 
 In fact, such a numerical verification can not describe any true asymptotic behaviour of $\frac{1}{\mbox{Vol}_{n-1}(\partial \Omega)\,\Lambda^{(n-1)/2}}\big[ \mathcal{N} (\Lambda)- a \mbox{Vol}_n (\Omega) \Lambda^{n/2}\big]$  as $\Lambda\to +\infty$,  where $a$ is the one-term coefficient of $\mathcal{N}(\Lambda)$. Numerical experimental verification should consider the case for sufficiently large $\Lambda$.  More precisely,  in order to verify $$\lim_{\Lambda\to +\infty} \frac{1}{\mbox{Vol}\,(\partial \Omega)\,\Lambda^{(n-1)/2} }\left[ \mathcal{N} (\Lambda)- a \mbox{Vol}_n (\Omega) \Lambda^{n/2} \right] = b'$$ for some constant $b'$, by the definition of limit, one should prove that for every number $\epsilon>0$ there is a number $M>0$ which depends on $\epsilon$ such that if $\Lambda>M$ then $$\bigg|\frac{1}{\Lambda^{(n-1)/2}}\left[ \mathcal{N} (\Lambda)- a \mbox{Vol}_n (\Omega) \Lambda^{n/2} \right]  - b'\bigg|<\epsilon.$$
       The authors of {\rm\cite{CaFrLeVa-22}} ``claimed'' that their result by using (non-professional)   ``numerically'' verifying for $0\le \Lambda\le 3000$ or $0\le \Lambda\le 2400$ or $0\le \Lambda\le 2600$ for the unit disk, flat cylinders and unit square. {\bf Such a so-called ``numerical verification'' method is not believable at all. The more serious problem is that these results are all  wrong because their eigenvalues calculations were based on some erroneous formulas for a unit disk and flat cylinders (see, Section 2 and Section 3)}. {{\bf This implies that all figures of \cite{CaFrLeVa-23} in numerical verification are wrong.}} 

\vskip 0.55 true cm

\noindent{\bf Remark 6.2.} \ \  {\bf Unlike we get our result in \cite{Liu-21} by applying (global) geometric analysis technique on the whole Riemannian manifold for the corresponding elastic Lam\'{e} operator (more importantly,  we  have enough regularity results for the elastic parabolic system (for example, Dong's theorem 3 and Theorem 4 in \cite{Do-12})), the authors of {\rm\cite{CaFrLeVa-23}} only discussed the elastic operator in a very special domain (i.e., the Euclidean upper half-space with flat metric $g_{jk}=\delta_{jk}$)}. Their method is \textcolor{blue}{``stretch $\Omega$ by a linear factor $\kappa>0$, note that the
eigenvalues then rescale as $\kappa^{-2}$, and check the rescaling of the geometric invariants and of (1.19)''} (see, from line -3 to line -1 on p.$\,$6 of {\rm\cite{CaFrLeVa-23}}). {\bf However, for studying the two-term asymptotics, such a method has lost a large amount of useful information. The reason is that the elastic waves have completely different propagation paths in a flat Euclidean space and a curved Riemannian manifold (elastic wave will be more complicated when the wave arrives at the boundary)}.
In a flat Euclidean space, $P$-wave and $S$-wave have many simple properties, but on a (curved) Riemannian manifold such waves have not similar properties. 
 In fact, in a Riemannian manifold $(\Omega, g)$, it follows from Lemma 2.1.1 of {\rm\cite{Liu-19}} that the elastic operator has the following local representation: \begin{eqnarray*} &&\mu\sum\limits_{j=1}^n\! \Big\{\! \Delta_g u^j \! +2 \!\sum\limits_{k,s,l=1}^n\! g^{kl} \Gamma_{sk}^j \frac{\partial u^s}{\partial x_l} \!+\!
   \sum_{k,s,l=1}^n \!\Big(g^{kl} \frac{\partial \Gamma^j_{sl}}{\partial x_k}\! +\!\sum\limits_{h=1}^n g^{kl} \Gamma_{hl}^j \Gamma_{sk}^h\! -\!\sum\limits_{h=1}^n g^{kl} \Gamma_{sh}^j \Gamma_{kl}^h \Big)u^s\Big\}\frac{\partial}{\partial x_j}\\
   &&\;\, +(\lambda+\mu) \operatorname{grad}\,\operatorname{div}\, \mathbf{u} \!+\! \mu\,\operatorname{Ric} (\mathbf{u}).\nonumber\end{eqnarray*}
{\bf In {\rm\cite{CaFrLeVa-23}} (and {\rm\cite{CaFrLeVa-22}}), by  regarding the above elastic operator defined on $(\Omega,g)$ as $\mu\sum_{j=1}^n \frac{\partial^2\mathbf{u}}{\partial x_j^2} +(\lambda+\mu) \nabla(\nabla\cdot \mathbf{u})$ on $\mathbb{\bar R}^{n}_{+}$, and further by erroneously taking a neighborhood of the boundary as the Euclidean upper half-space, {\bf the authors of {\rm\cite{CaFrLeVa-23}} had completely changed the original eigenvalues problems into other different spectral problems}. {\bf The reader can not see where the branching Hamiltonian billiards condition (i.e., the corresponding billiard is neither dead-end nor absolutely periodic) on the Riemannian manifold is used in the proof of {\rm\cite{CaFrLeVa-23}}.}} In a Riemannian manifold, in order to get the two-term asymptotics of the counting function, more geometric analysis tools should be applied (See Response 6 in Section 3 in \cite{Liu-22b}).

\vskip 1.19 true cm

\section{Erroneous conclusions in  \cite{CaFrLeVa-23}  and \cite{SaVa-97}}

\vskip 0.62 true cm

 First, the (elastic) parabolic equations   
\begin{eqnarray}\label{2023.6.16-2}  \left\{ \begin{array}{ll} \big(\frac{\partial}{\partial t} +P_g\big) \mathbf{u} (t, x)=0,\;\;\,  \, t>0, \, x\in M,\\ [1mm]
\mathbf{u}(0,x)=\boldsymbol{\phi} (x), \,\;\; x\in M \end{array}\right. \end{eqnarray}
defined on a Riemannian manifold $(M,g)$ ($M$ may be closed or unbounded) has the symbol equations:
\begin{eqnarray}\label{2023.6.17-1} \left\{ \begin{array}{ll}  \big( \frac{\partial}{\partial t} +\mathbf{A}_g (x, \xi)\big)\mathbf{\hat{u}}(t, \xi)=0, \;\; \, t>0,\, x\in M,\, \xi\in \mathbb{R}^n,\\ [1mm]
\hat{\mathbf{u}}(0,\xi)=\boldsymbol{\hat{\phi}}(\xi),\end{array}\right.  \end{eqnarray}  
where $P_g$ and  $\mathbf{A}_g(x, \xi)$ are the (elastic) Lam\'{e} operator and  the full symbol of $P_g$, respectively,  which are explicitly given on p.$\,$10177--10178 in \cite{Liu-21}. 
Clearly, (\ref{2023.6.17-1}) is a vector-valued  ordinary differential equations in $t$, which has a unique solution 
\begin{eqnarray} \label{2023.6.17-2} \mathbf{\hat{u}}(t, \xi) =\boldsymbol{\hat{\phi}}(\xi)\, e^{-t\mathbf{A}_g (x,\xi)}.\end{eqnarray}
According to the Cauchy integral formula we have 
\begin{eqnarray}\label{2023.6.30-1} e^{-t\mathbf{A}_g(x, \xi)} =\frac{1}{2\pi i } \int_{\mathcal{C}} e^{-t\tau} \big(\tau- \mathbf{A}_g (x, \xi)\big)^{-1} d\tau, \end{eqnarray} 
where $\mathcal{C}$ is a  contour around the positive real axis in the complex plane because the principal symbol of $P_g$ is a positive-definite matrix. It follows that 
 \begin{eqnarray*} \label{2023.6.17-3} \mathbf{\hat{u}}(t, \xi) =\Big(\frac{1}{2\pi i } \int_{\mathcal{C}} e^{-t\tau} \big(\tau- \mathbf{A}_g (x, \xi)\big)^{-1} d\tau \Big)\boldsymbol{\hat{\phi}}(\xi),\end{eqnarray*}
so that  
\begin{eqnarray} \mathbf{u}(t, x)=\frac{1}{(2\pi)^n} \int_{\mathbb{R}^n} e^{i x\cdot \xi}\Big(\frac{1}{2\pi i } \int_{\mathcal{C}} e^{-t\tau} \big(\tau- \mathbf{A}_g (x, \xi)\big)^{-1} d\tau \Big)\boldsymbol{\hat{\phi}}(\xi)\,d\xi.  \end{eqnarray} 
In particular, from this and (\ref{2023.6.30-1}) we have $\lim\limits_{t\to 0^+} \mathbf{u}(t,x)=\boldsymbol{\phi}(x)$ for $x\in M$.
Noting that $\big(\tau- \mathbf{A}_g (x, \xi)\big)^{-1}\in S^{-2}_{cl}$ (see \cite{Gr-86}, or p.$\,$10174 of \cite{Liu-21}), we have  \begin{eqnarray*}\big(\tau- \mathbf{A}_g (x, \xi)\big)^{-1}\,\sim \,\sum_{l\ge 0} \mathbf{q}_{-2-l} (x, \xi,\tau), \end{eqnarray*} 
 where $\mathbf{q}_{-2-l}(t,x,\xi)$ is homogeneous of degree $-2-l$ in $\xi$ for $|\xi|\ge 1$ (see (3.10) on p.$\,$10180 in \cite{Liu-21}).
Thus we obtain 
\begin{eqnarray} \mathbf{u}(t, x)= \frac{1}{(2\pi)^{n}}\int_{\mathbb{R}^n} e^{i x\cdot \xi}\Big(\frac{1}{2\pi i } \int_{\mathcal{C}} e^{-t\tau} \sum_{l\ge 0} \mathbf{q}_{-2-l} (x, \xi,\tau)\; d\tau \Big)\boldsymbol{\hat{\phi}}(\xi)\,d\xi, \;\; \, t>0, \; x\in M.  \end{eqnarray} 
This implies that a fundamental solution $\mathbf{F}(t,x,y)$ (an $n\times n$ matrix-valued function) of the  initial value problem (\ref{2023.6.16-2}) (i.e., $\mathbf{F}(t,x,y)$ satisfies 
 \begin{eqnarray}\label{2023.6.7-11}\left\{ \begin{array}{ll}
  \Big( \frac{\partial}{\partial t}+{P}_g\Big) \mathbf{F} (t,x, y)=
 0, \quad t>0, \; x,y\in M,\\
  \mathbf{F}(0, x,y)=\boldsymbol{\delta}(x-y)) \end{array}\right.\end{eqnarray}
can be represented by 
\begin{eqnarray} \label{2023.6.17-6} \mathbf{F}(t, x,y)=  \frac{1}{(2\pi)^{n}}\int_{\mathbb{R}^n} e^{i (x-y)\cdot \xi}\Big(\frac{1}{2\pi i } \int_{\mathcal{C}} e^{-t\tau}  \sum_{l\ge 0}^{+\infty} \mathbf{q}_{-2-l} (y, \xi, \tau)\; d\tau \Big)\,d\xi, \,\;\;\; x, y\in M.  \end{eqnarray}
Since the Lam\'{e} operator are smooth on smooth Riemannian manifold $(M, g)$, we see from (\ref{2023.6.17-6}) and $\sum_{l\ge 0} \mathbf{q}_{-2-l} (x$, $\xi$,$\tau)\in S^{-2}_{cl}$ that  $\mathbf{F}\in [C^\infty ((0,+\infty)\times M \times M)]_{n\times n}$. In addition, for any $0<\epsilon <T<+\infty$, $\,\mathbf{F}(t, x,y)$ is uniformly continuous on $[\epsilon, T]\times M\times M$. Of course, for any $x, y\in M$ and $x\ne y$, we have $\lim_{t\to 0^+} \mathbf{F}(t,x,y)=\boldsymbol{\delta} (x-y)$.   

  When $M= \mathbb{R}^n$, the $\mathbf{F} (t,x,y)$ in  (\ref{2023.6.17-6}) is just a fundamental solution of  the initial value problem for the (elastic) parabolic equations 
\begin{eqnarray} \left\{ \begin{array}{ll} \frac{\partial \mathbf{u}(t,x)}{\partial t}-
\mu \Delta \mathbf{u}(t,x)  -(\lambda+\mu)\, \mbox{grad}\, \mbox{div}\; \mathbf{u}(t,x) =0, \;\;\; t>0,\, x\in \mathbb{R}^n, \\ [1mm]
\mathbf{u}(0,x)=\boldsymbol{\phi}(x), \;\;\, x\in \mathbb{R}^n.\end{array}\right.\end{eqnarray}
In particular, if $M=\mathbb{R}^n$ and $\lambda+\mu=0$, then (\ref{2023.6.17-6}) reduces to  the  fundamental solution of the classical heat equations, 
\begin{eqnarray} \label{2023.6.17-9} \mathbf{F}(t,x,y) = \frac{1}{(4\pi \mu t)^{n/2}} e^{- |x-y|^2 /4t}\, \mathbf{I}_n,\end{eqnarray} 
where $|x-y|=\sqrt{\sum_{k=1}^n (x_k-y_k)^2}\;$ for $x,y\in \mathbb{R}^n$.

\vskip 0.39 true cm
By using symbol analysis, we have the following:

\vskip 0.24 true cm
\noindent{\bf Proposition 7.1.} \ {\it Let $(\Omega,g)$ be a smooth compact connected $n$-dimensional Riemannian manifold with smooth boundary $\partial \Omega$. Suppose that $(\Omega, g)$ is such that the corresponding billiards is neither dead-end nor absolutely periodic. If $ \mathcal{N}^\mp (\Lambda)$ has the following asymptotic expansion: \begin{eqnarray} \label{2023.6.5-1}  \mathcal{N}^\mp (\Lambda) = a \mbox{Vol}_n (\Omega) \Lambda^{n/2} + b_1^{\mp} \mbox{Vol}_{n-1} (\partial \Omega) \Lambda^{(n-1)/2}+ o(\Lambda^{(n-1)/2}) \;\; \mbox{as}\;\; \Lambda\to +\infty,\end{eqnarray}
where $a$ is given by (\ref{2023.5.5-2}), and $b_1^-$ and $b_1^+$ are two constants depending only on the Lam\'{e} parameters $\mu$ and $\lambda$, 
 then \begin{eqnarray} \label{2023.6.8/1} b_1^-+ b_1^+=0.\end{eqnarray}} 

 \vskip 0.39 true cm

  \noindent {\it Proof.}  \  Since  $ \mathcal{N}^\mp (\Lambda)$ have the two-term asymptotic expansions (\ref{2023.6.5-1}), it is easy to verify from this and  (\ref{2023.5.5-11}) that 
\begin{eqnarray}\label{2023.5.5-8}  \mathcal{Z}^{\mp}(t) = c\, \mbox{Vol}_n( \Omega)\, t^{-n/2} + d_1^{\mp}\, \mbox{Vol}_{n-1} (\partial \Omega) \,t^{-(n-1)/2} + o(t^{-(n-1)/2})\;\; \mbox{as} \;\; t\to 0^+,\end{eqnarray}
where  \begin{eqnarray} \label{2023.5.5-6} c = \Gamma \Big( 1+\frac{n}{2}\Big)\, a, \, \; \,\;  d_1^{\mp} = \Gamma\Big(1+\frac{n-1}{2}\Big) b_1^{\mp}.\end{eqnarray}     
Let $\mathbf{G}^{-} (t,x,y)$ (respectively, $\mathbf{G}^+(t,x,y)$) be Green's function of the following (elastic) parabolic  system (defined on $\Omega$) with Dirichlet (respectively, Neumann) boundary condition:
  \begin{eqnarray} \label{2023.6.7-2}  \left\{  \begin{array}{ll}  \frac{\partial \mathbf{G}^-(t,x,y)}{\partial t} + {P}_g\mathbf{G}^-(t,x,y)=0,\;\;\; t>0,\, x,\, y\in\Omega,\\
 \mathbf{G}^- (t, x,y)=0, \;\;\;\;  t>0,\; x\in \Omega, \,\;  y\in \partial \Omega,\\
 \mathbf{G}^-(0, x,y)=\boldsymbol{\delta}(x-y), \;\;\;  x,\,y\in \Omega\end{array} \right. \end{eqnarray}
(respectively, 
 \begin{eqnarray}  \label{2023.6.7-3} \left\{   \begin{array}{ll}   \frac{\partial \mathbf{G}^+(t,x,y)}{\partial t} + {P}_g\mathbf{G}^+(t,x,y)=0, \,\;\;\; t>0, \;x, \,y\in \Omega,\\
\mathcal{F} \big(\mathbf{G}^+ (t, x,y)\big)=0, \;\; \,\,t>0,\;  x\in \Omega, \;\;  y\in \partial \Omega,\\
 \mathbf{G}^+(0, x,y)=\boldsymbol{\delta}(x-y),\;\;\;  x,\,y\in \Omega),\end{array} \right. \end{eqnarray}
 where the Lam\'{e} operator is acted in the variable $y$, and $\mathcal{F} \big( \mathbf{G}^+\big):= \mu \big(\nabla \mathbf{G}^+ +(\nabla \mathbf{G}^+)^T\big)\boldsymbol{\nu} +\lambda (\mbox{div}\, \mathbf{G}^+)\boldsymbol{\nu}$ on $\partial \Omega$. It is clear that $\mathbf{G}^{\mp}\in [C^\infty ((0, +\infty)\times \bar \Omega \times \bar\Omega)]_{n\times n}$ and that $\lim_{t\to 0^+} \mathbf{G}^\mp(t,x,y)=\boldsymbol{\delta} (x-y)$  for any fixed $x, y\in \bar \Omega$ and $x\ne y$. This implies that
for any $\boldsymbol{\phi}\in [C^\infty_0(\Omega)]^n$ (respectively, $\boldsymbol{\phi}\in [C^\infty (\Omega)]^n$ with $\big[ \mu \big(\nabla \boldsymbol{\phi} +(\nabla \boldsymbol{\phi})^T\big)\boldsymbol{\nu} +\lambda (\mbox{div}\, \boldsymbol{\phi})\boldsymbol{\nu}\big)\big]\big|_{\partial \Omega}=0$), the vector-valued functions $\mathbf{u}^\mp(t, y):=\int_{\Omega} \mathbf{G}^\mp (t, x, y)\,\boldsymbol{\phi}(x)\,dx$ satisfy 
  \begin{eqnarray} \label{2023.7.7-2}  \left\{  \begin{array}{ll}  \frac{\partial \mathbf{u}^-(t,y)}{\partial t} + {P}_g\mathbf{u}^-(t,y)=0,\;\;\; t>0,\, y\in\Omega,\\
 \mathbf{u}^- (t,y)=0, \;\;\;\;  t>0,\;\;  y\in \partial \Omega,\\
 \lim\limits_{t\to 0^+} \mathbf{u}^-(t, y)=\boldsymbol{\phi}(y), \;\;\;  y\in \Omega\end{array} \right. \end{eqnarray}
(respectively, 
  \begin{eqnarray} \label{2023.7.7-3}  \left\{  \begin{array}{ll}  \frac{\partial \mathbf{u}^+(t,y)}{\partial t} + {P}_g\mathbf{u}^+(t,y)=0,\;\;\; t>0,\, x,\, y\in\Omega,\\
 \mathbf{u}^+ (t,y)=0, \;\;\;\;  t>0,\;   y\in \partial \Omega,\\
 \lim\limits_{t\to 0^+} \mathbf{u}^+(t, y)=\boldsymbol{\phi}(y), \;\;\;  y\in \Omega\,).\end{array} \right. \end{eqnarray}

Let $\{{\mathbf{u}}_k^\mp\}_{k=1}^\infty$ be the orthonormal eigenvectors of the elastic operators $P_g^\mp$ corresponding to the eigenvalues $\{\tau_k^\mp\}_{k=1}^\infty$, then  Green's functions  ${\mathbf{G}}^\mp(t, x, y)=e^{-t P_g^\mp} \boldsymbol{\delta}(x-y)$ are given by \begin{eqnarray} \label{18/12/18} {\mathbf{G}}^\mp(t,x,y) =\sum_{k=1}^\infty e^{-t \tau_k^\mp} {\mathbf{u}}_k^\mp(x)\otimes {\mathbf{u}}_k^\mp(y),\end{eqnarray}
where \begin{equation*}  {\mathbf{u}}_k^\mp(x)\otimes {\mathbf{u}}_k^\mp(y)= \renewcommand\arraystretch{1.5} \left( \begin{matrix} u_{k1}^\mp(x) u_{k1}^\mp (y)& u_{k2}^\mp (x) u_{k1}^\mp (y)
 &\cdots & u_{kn}^\mp(x) u_{k1}^\mp(y) \\ 
u_{k1}^\mp(x) u_{k2}^\mp(y) & u_{k2}^\mp (x) u_{k2}^\mp(y)  &\cdots & u_{kn}^\mp(x) u_{k2}^\mp(y)\\
\cdots & \cdots & \cdots & \cdots\\
u_{k1}^\mp(x) u_{kn}^\mp(y) & u_{k2}^\mp (x) u_{kn}^\mp(y)  &\cdots &u_{kn}^\mp(x) u_{kn}^\mp(y)\end{matrix}\right)
 \end{equation*}
and  $\mathbf{u}_k^\mp(x)= (u_{k1}^\mp(x), u_{k2}^\mp(x), \cdots, u_{kn}^\mp (x))^T$. 
Note that 
\begin{eqnarray} \label{1-0a-21}\int_{\Omega} \mbox{Tr}\,({\mathbf{G}}^\mp(t,x,x)) dV=\sum_{k=1}^\infty e^{-t \tau_k^\mp}=\mathcal{Z}^\mp (t) .\end{eqnarray}
From this and (\ref{2023.5.5-8}), we get   \begin{eqnarray} \label{2023.6.7-5}&& \frac{1}{2}\left[ \sum_{k=1}^\infty e^{-t \tau_k^-} +\sum_{k=1}^\infty e^{-t \tau_k^+}\right] = c\, \mbox{Vol}_n( \Omega) \,t^{-n/2} + \frac{1}{2}(d_1^-+ d_1^+)\, \mbox{Vol}_{n-1} (\partial \Omega)\, t^{-(n-1)/2}\\
 && \qquad\qquad \qquad \qquad \quad \;\quad \;\;\, \qquad + o(t^{-(n-1)/2})\;\; \; \mbox{as} \;\; t\to 0^+.\nonumber\end{eqnarray}

On the other hand, we can also structure the Green's functions $\mathbf{G}^\mp(t,x,y)$ by the following method: 
Let $\mathcal{M}$ is the double manifold of $\Omega$, and $\mathcal{P}$ (defined on $\mathcal{M}$) is the double operator of the Lam\'{e} $P_g$ on $\Omega$ (see, the previous discussion or (\ref{2021.2.6-3})--
(\ref{2022.10.18-2}) in Appendix below). Since  $\mathcal {M}$ is a closed manifold, there is a fundamental solution of $\mathbf{K} (t,x,y)$ defined on $(0,+\infty)\times \mathcal{M}\times \mathcal{M}$ such that \begin{eqnarray}\label{2023.6.29-8}\left\{ \begin{array}{ll}
   \Big( \frac{\partial}{\partial t}+\mathcal{P}\Big) \mathbf{K} (t,x, y)=0,\;\;\, x, y\in \mathcal{M},\\
  \mathbf{K}(0, x,y)=\boldsymbol{\delta}(x-y) \end{array}\right.\end{eqnarray}
Note  that (see Appendix below) $\mathbf{K}\in [C^\infty ((0,+\infty)\times (\mathcal{M}\setminus \partial \Omega) \times (\mathcal{M}\setminus\partial \Omega))]_{n\times n} \cap  [C^1 ((0,+\infty)\times \mathcal{M}\times \mathcal{M})]_{n\times n}$. It follows that 
$\mathbf{G}^\mp(t,x,y):=\mathbf{K}(t, x,y) \mp\mathbf{K} (t,x, \overset{*}{y}(y)) -\mathbf{H}(t,x,y)$ are just the Green's functions with the Dirichlet and Neumann (i.e. free) boundary conditions, where $\overset{*}{y}(y):=\overset{*}{y}= (x',-x_n)$ for any $y=(x',x_n)\in M$.  
Thus \begin{eqnarray} \label{23.12.24-2} \;\;\;\;\;\mathbf{G}(t,x,y):=\frac{1}{2}\big(\mathbf{G}^-(t,x,y)+\mathbf{G}^+(t,x,y)\big)= \mathbf{K}(t,x,y)-\mathbf{H}(t,x,y) \;\;\;\mbox{for all}\;\; x, y\in {M}. \end{eqnarray}  Since $|\mathbf{H}(t,x,y)|\le C$ for all $0\le t\le 1$ and  $x,y\in \bar M$ (see Appendix below), we get 
$$\bigg|\int_{\Omega} \mbox{Tr}\big(\mathbf{H}(t,x,x)\big)\,dx\bigg| \le nC \big(\mbox{Vol}(\Omega) \big)\;\;\, \mbox{for}  \;\; 0\le t\le 1,$$  so that for $n\ge 2$,
\begin{eqnarray}\label{23.12.24-1}\int_\Omega\mbox{Tr}\big( \mathbf{H}(t,x,x)\big)\,dx= o(t^{-n/2})\;\;\;\mbox{as}\;\; t\to 0^+.\end{eqnarray} 

  Now, by virtue of (\ref{2023.6.17-6}) we have  
\begin{eqnarray} \label{2023.6.17-16} \mathbf{K}(t, x,x)=  \frac{1}{(2\pi)^{n}}\int_{\mathbb{R}^n} \Big(\frac{1}{2\pi i } \int_{\mathcal{C}} e^{-t\tau}  \sum_{l\ge 0}^{+\infty} \mathbf{q}_{-2-l} (x, \xi, \tau)\; d\tau \Big)\,d\xi.  \end{eqnarray}
From the result on p.$\,$10182 of \cite{Liu-21}, we know that
   \begin{eqnarray}\label{2023.5.7-1}
   && \;\;\;\;{\mathbf{q}}_{-2} (x,\xi,\tau) =\frac{1}{\tau- \mu \sum\limits_{l,m=1}^n g^{lm}\xi_l \xi_m }\,{\mathbf{I}}_n\\
 && \qquad \;\;\;\; \; \; +\frac{\mu+\lambda}{  \big(\tau- \mu \sum\limits_{l,m=1}^n g^{lm}\xi_l \xi_m \big)\big(\tau- (2\mu+\lambda)  \sum\limits_{l,m=1}^n g^{lm}\xi_l \xi_m \big)}\begin{bmatrix} \sum\limits_{r=1}^n g^{1r} \xi_r \xi_1 &\cdots &  \sum\limits_{r=1}^n g^{1r} \xi_r \xi_n\\
\vdots& {} &\vdots \\
 \sum\limits_{r=1}^n g^{nr} \xi_r \xi_1 &\cdots &  \sum\limits_{r=1}^n g^{nr} \xi_r \xi_n\end{bmatrix}\nonumber\end{eqnarray}
 and   \begin{eqnarray}\label{2023.5.7-2}  && \mbox{Tr}\, \big({\mathbf{q}}_{-2} (x,\xi,\tau)\big)=
   \frac{n}{\big(\tau \!-\!\mu \!\sum_{l,m=1}^n \!g^{lm} \xi_l\xi_m\big)} \\
 && \qquad \qquad \qquad \qquad \;\;  + \frac{(\mu+\lambda)\sum_{l,m=1}^n g^{lm} \xi_l\xi_m}{
  \big(\!\tau\! -\!\mu\! \sum_{l,m=1}^n \!g^{lm} \xi_l\xi_m\!\big)
 \big(\!\tau \!-\!(2\mu\!+\!\lambda)  \!\sum_{l,m=1}^n \!g^{lm} \xi_l\xi_m\!\big)}
.\nonumber \end{eqnarray}
Also, it is easy to find from (3.15) on p.$\,$10182 in \cite{Liu-21} that $\mbox{Tr}\,(\mathbf{q}_{-3} (x,\xi,\tau))$ is an odd function in $\xi=(\xi_1, \cdots,\xi_n)\in \mathbb{R}^n$.    
   For each $x\in \Omega$, we use a geodesic normal coordinate system centered at this $x$. It follows from \S11 of Chap.1 in \cite{Ta-1}
        that in such a coordinate system, $g_{jk}(x)=\delta_{jk}$ and $\Gamma_{jk}^l(x)=0$. Then (\ref{2023.5.7-2}) reduces to
        \begin{eqnarray}\label{2020.7.6-3} \quad \quad \,\;\mbox{Tr} \,\big({\mathbf{q}}_{-2} (x, \xi,\tau)\big)=
     \frac{n}{(\tau -\mu |\xi|^2)}+
      \frac{(\mu+\lambda)|\xi|^2}{(\tau -\mu |\xi|^2) (\tau-(2\mu+\lambda)|\xi|^2)},\end{eqnarray}
  where $|\xi|=\sqrt{\sum_{k=1}^n \xi^2_k}$ for any $\xi\in {\mathbb{R}}^n$.
   By applying the residue theorem (see, for example, Chap.$\,$4, \S5 in \cite{Ahl}) we get
 \begin{eqnarray} \label{3.10} &&\; \frac{1}{2\pi i} \int_{\mathcal{C}} e^{-t\tau} \bigg(\frac{n}{(\tau -\mu |\xi|^2)}\!+\!
      \frac{(\mu+\lambda)|\xi|^2}{(\tau -\mu |\xi|^2) (\tau-(2\mu+\lambda)|\xi|^2)} \bigg) d\tau= (n\!-\!1) e^{-t\mu|\xi|^2}\! +\!  e^{-t(2\mu+\lambda)|\xi|^2}.\end{eqnarray}
       It follows that
      \begin{eqnarray}\label{2020.7.10-1} \frac{1}{(2\pi)^n}\!\!\!\!&\!\!\!&\!\!\!\!\! \!\!\! \!\! \int_{{\mathbb{R}}^n}\Big( \frac{1}{2\pi i} \int_{\mathcal{C}} e^{-t\tau}\, \mbox{Tr}\,({\mathbf{q}}_{-2} (x, \xi,\tau) ) d\tau \Big) d\xi \\
      \!\!\! &=\!\!\!&
    \frac{1}{(2\pi)^n} \int_{{\Bbb R}^n} \bigg((n-1) e^{-t\mu|\xi|^2} +  e^{-t(2\mu+\lambda)|\xi|^2}\bigg)
         d\xi \nonumber\\
               \!\!\! &=\!\!\!&\frac{n-1}{(4\pi \mu t)^{n/2}} +  \frac{1}{(4\pi (2\mu+\lambda) t)^{n/2}} \;\;\;\mbox{uniformly for all} \;\; x\in \Omega,
                  \nonumber\end{eqnarray}
                  and hence
        \begin{eqnarray}\label{2023.7.12-3}  && \int_{\Omega}\! \left\{\! \frac{1}{(2\pi)^n}\! \int_{{\mathbb{R}}^n}\!\Big( \frac{1}{2\pi i} \int_{\mathcal{C}} e^{-t\tau}\, \mbox{Tr}\,({\mathbf{q}}_{-2} (x, \xi,\tau) ) d\tau\! \Big) d\xi\!\right\}\! dV\\
        && \qquad \;\;=\Big(\frac{n-1}{(4\pi \mu t)^{n/2}} \! +\!  \frac{1}{(4\pi (2\mu\!+\!\lambda) t)^{n/2}}\!\Big){\mbox{Vol}}(\Omega).\nonumber\end{eqnarray}
                Then,  for $l\ge 1$, it can be verified that
                  $\mbox{Tr}\, ({\mathbf{q}}_{-2-l} (x, \xi, \tau))$ is a sum of finitely many terms, each of which  has the following form:
                  $$\frac{r_k(x, \xi)}{(\tau-\mu \sum_{l,m=1}^n g^{lm} \xi_l \xi_m )^s (\tau -(2\mu +\lambda)\sum_{l,m=1}^n g^{lm} \xi_l\xi_m )^j },$$
where $k-2s-2j=-2-l$, and $r_k(x,\xi)$ is the symbol independent of $\tau$ and homogeneous of degree $k$.
 Again we take the geodesic normal coordinate systems center at $x$ (i.e., $g_{jk}(x)=\delta_{jk}$ and $\Gamma_{jk}^l(x)=0$), by applying residue theorem we see  that, for $l\ge 1$ , \begin{eqnarray*}\label{2023.7.12-1} && \frac{1}{(2\pi)^n} \int_{{\mathbb{R}}^n}\Big( \frac{1}{2\pi i} \int_{\mathcal{C}} e^{-t\tau}\, \mbox{Tr}\,({\mathbf{q}}_{-2-l} (x, \xi,\tau) ) d\tau \Big) d\xi=O(t^{l-\frac{n}{2}}) \;\;\mbox{as}\;\, t\to 0^+ \;\;\, \mbox{uniformly for} \;\, x\in \Omega.
 \end{eqnarray*}
Therefore
\begin{eqnarray}\label{2023.7.12-10} && \int_{\Omega} \bigg\{\frac{1}{(2\pi)^n} \int_{{\mathbb{R}}^n}\Big( \frac{1}{2\pi i} \int_{\mathcal{C}} e^{-t\tau}\, \sum_{l\ge 1}\mbox{Tr}\,({\mathbf{q}}_{-2-l} (x, \xi,\tau) ) d\tau \Big) d\xi\bigg\} dV =O(t^{1-\frac{n}{2}}) \;\;\mbox{as}\;\, t\to 0^+.
 \end{eqnarray} 
Combining (\ref{2023.6.17-16}), (\ref{2023.7.12-3}) and (\ref{2023.7.12-10}), we obtain
       \begin{eqnarray} \label{43.11}&& \int_{\Omega}\mbox{Tr}\,({\mathbf{K}}(t, x, x)) \, dV =    \bigg[  \frac{n-1}{(4\pi \mu t)^{n/2}} +  \frac{1}{(4\pi (2\mu+\lambda) t)^{n/2}}\bigg]{\mbox{Vol}}(\Omega)+O(t^{1-\frac{n}{2}})\;\;\mbox{as}\;\; t\to 0^+, \end{eqnarray}
and hence  by (\ref{23.12.24-2})
\begin{eqnarray} \label{43.11}&& \int_{\Omega}\mbox{Tr}\,({\mathbf{G}}(t, x, x)) \, dV =    \bigg[  \frac{n-1}{(4\pi \mu t)^{n/2}} +  \frac{1}{(4\pi (2\mu+\lambda) t)^{n/2}}\bigg]{\mbox{Vol}}(\Omega)+O(t^{1-\frac{n}{2}})\;\;\mbox{as}\;\; t\to 0^+. \end{eqnarray}
From  (\ref{1-0a-21}), (\ref{2023.6.7-5}) and (\ref{43.11}),  we have   \begin{eqnarray} \label{44.11}&& \frac{1}{2}\left[ \sum_{k=1}^\infty e^{-t \tau_k^-} +\sum_{k=1}^\infty e^{-t \tau_k^+}\right]   =   c\, \mbox{Vol}_n( \Omega) t^{-n/2}+O(t^{1-\frac{n}{2}})\;\;\mbox{as}\;\; t\to 0^+, \end{eqnarray}
Comparing the growth order of  time variable $t$ as $t\to 0^+$ from (\ref{2023.6.7-5}) and (\ref{44.11}), we immediately find that $d_1^-+d_1^+=0$, and hence $b_1^-+b^+_1=0$ by (\ref{2023.5.5-6}).  \qed

\vskip 0.52 true cm  

\noindent{\bf Remark 7.2.} \ {\bf From (\ref{202.5.5-10}) and (\ref{202.5.5-3})  of Theorem 1.2, it is easy to verify that $b^-+b^+\ne 0$ (see, for example,  Table 1 on p.$\,$9 and Table 2 on p.$\,$10 in \cite{CaFrLeVa-23}). However,  $b^-+b^+$ must be zero according to Proposition 7.1.  This is a contradiction,}  {{\bf which implies that the conclusions in \cite{CaFrLeVa-23} (and  \S6.3 of \cite{SaVa-97}) are wrong}}.

\vskip 0.39 true cm 

\noindent{\bf Remark 7.3.} \  Let us point out that the range  of the Lam\'{e} coefficients $\mu$ and $\lambda$ can be taken as $\mu>0$ and $\lambda+2\mu> 0$ because in this range, the Lam\'{e} operator $P_g$ is strongly elliptic (see  \cite{Liu-39} or p.$\,$208 of \cite{LiQin-13}, \cite{McL}); therefore the corresponding elastic eigenvalue problems can well be done.  In the special case when $\mu>0$ and  $\lambda+\mu\to 0^+$, the elastic Lam\'{e} operator $P_g$ reduces to the classical  Laplace-type operator $\mu \nabla^*\nabla -\mu\, \mbox{Ric}$ on smooth Riemannian manifold $(\Omega,g)$, where $\mu \nabla^*\nabla$ is the Bochner Laplacian on $\Omega$. 
In fact, for $\mu>0$ and $\lambda+\mu= 0$, the elstic Dirichlet (respectively, Neumann) eigenvalue problem  is 
\begin{eqnarray} \label {2023.7.11-1}\left\{ \begin{array}{ll} \mu \nabla^* \nabla \mathbf{u}  -\mu \, \mbox{Ric} (\mathbf{u})= \tau \mathbf{u}, \;  \;\; &\mbox{in} \;\; \Omega,
\\
 \mathbf{u}=0\;\;\; &\mbox{on}\;\; \partial \Omega\end{array} \right.\end{eqnarray} (respectively,  
\begin{eqnarray} \label {2023.7.10-1}\left\{ \begin{array}{ll} \mu \nabla^* \nabla \mathbf{u}  -\mu \, \mbox{Ric} (\mathbf{u})= \tau \mathbf{u}, \;  \;\; &\mbox{in} \;\; \Omega,
\\
2\mu \, (\mbox{Def}\, \mathbf{u})^\# \boldsymbol{\nu} -\mu\, (\mbox{div}\, \mathbf{u} )\boldsymbol{\nu}=0\;\;\; &\mbox{on}\;\; \partial \Omega),\end{array} \right.\end{eqnarray} 
 where $\tau$ is the elastic Dirichlet (respectively,  Neumann) eigenvalue.   
 From p.$\,$10  of \cite{Liu-19} and (2.19)  of  \cite{Liu-19}, we have 
 \begin{eqnarray*}  -\mu \nabla^*\nabla \mathbf{u} +(\lambda+ \mu) \,\mbox{grad}\, \mbox{div}\, \mathbf{u}  +\mu \, \mbox{Ric}(\mathbf{u}) = -2\mu \, \mbox{Def}^* \mbox{Def} \, \mathbf{u} +\lambda \,\mbox{grad}\, \mbox{div}\, \mathbf{u}.\end{eqnarray*}
Note that  in \cite{Liu-19} we have actually given an alterative proof for  the Weitzenbock formula
\begin{eqnarray} \label{2023.7.10-9} 2\, \mbox{div}\; \mbox{Def}\; \mathbf{u} =- \nabla^* \nabla \mathbf{u} +\mbox{grad}\; \mbox{div}\; \mathbf{u} +\mbox{Ric}\, (\mathbf{u}).\end{eqnarray}
 Thus,  for $\mu>0$ and $\mu+\lambda=0$, the eigenvalue problem (\ref{2023.7.11-1}) and (\ref{2023.7.10-1}) can also be re-written as 
\begin{eqnarray} \label {2023.7.11-18}\left\{ \begin{array}{ll}2 \mu \,\mbox{Def}^*\, \mbox{Def} \, \mathbf{u} +\mu\, \mbox{grad}\; \mbox{div}\; \mathbf{u}= \tau \mathbf{u}, \;  \;\; &\mbox{in} \;\; \Omega,
\\
 \mathbf{u}=0\;\;\; &\mbox{on}\;\; \partial \Omega,\end{array} \right.\end{eqnarray}
and 
\begin{eqnarray} \label {2023.7.10-8}\left\{ \begin{array}{ll}2 \mu \,\mbox{Def}^*\, \mbox{Def} \, \mathbf{u} +\mu\, \mbox{grad}\; \mbox{div}\; \mathbf{u}= \tau \mathbf{u}, \;  \;\; &\mbox{in} \;\; \Omega,
\\
2\mu  (\mbox{Def}\, \mathbf{u})^\# \boldsymbol{\nu} -\mu\, (\mbox{div}\, \mathbf{u} )\boldsymbol{\nu}=0\;\;\; &\mbox{on}\;\; \partial \Omega,\end{array} \right.\end{eqnarray}
 respectively.  Here we have used the fact that  $\mbox{Def}^* \mathbf{v} =-\,\mbox{div}\, \mathbf{v}$ 
  for any tensor field $\mathbf{v}$ of type $(0,2)$,  where $(\mbox{div}\, \mathbf{v})^j =v^{jk}_{\;\;\;\; ;k}$.
It follows from  (2.18) of \cite{Liu-19} that  
\begin{eqnarray} \label{23/7.18-1} && - \int_{\Omega} \big(\lambda(\mbox{div}\,\mathbf{u})(\mbox{div}\, \mathbf{v}) \!+\!2\,\mu\,\langle \mbox{Def}\, \mathbf{u}, \mbox{Def}\, \mathbf{v}\rangle \big) dV\\
&& \qquad \quad =  \int_{\Omega} \langle (\lambda\, \mbox{grad}\, \mbox{div} \!-\!2\mu \, \mbox{Def}^* \mbox{Def})\mathbf{u}, \mathbf{v}\rangle dV\nonumber 
- \int_{\partial \Omega} \langle 2\mu (\mbox{Def}\, \mathbf{u})^\# \boldsymbol{\nu} \!+\!\lambda (\mbox{div}\, \mathbf{u})\boldsymbol{\nu},\mathbf{v}\rangle dS.\nonumber\end{eqnarray}
 Therefore,  for both boundary conditions, the corresponding Rayleigh quotient is   
\begin{eqnarray} \label{2023.7.10-2} \frac{\int_{\Omega} \big[ 2\mu \langle  \mbox{Def}\, \mathbf{u}, \mbox{Def}\, \mathbf{u}\rangle -\mu (\mbox{div}\, \mathbf{u})^2 \big] dV}{ \int_\Omega |\mathbf{u}|^2dV}. \end{eqnarray} 

Since the operator $\mu \nabla^* \nabla   -\mu \, \mbox{Ric}\, $ (i.e., $2\mu \,\mbox{Def}^*\, \mbox{Def}+ \mu \,\mbox{grad}\,\mbox{div}$) with zero Dirichlet boundary condition $\mathbf{u}\big|_{\partial \Omega}=0$ (respectively,  zero Neumann  boundary condition $\big[2\mu  (\mbox{Def}\, \mathbf{u})^\# \boldsymbol{\nu} - \mu\, (\mbox{div}\, \mathbf{u} )\boldsymbol{\nu}\big]\big|_{\partial \Omega}=0$)  is an unbounded, self-adjoint and  non-negative semidefinite operator in $[H_0^1(\Omega)]^n$ (respectivley, $[H^1(\Omega)]^n$) with discrete spectrum $0<  \tau_1^- < \tau_2^- \le \cdots \le \tau_k^- \le \cdots \to +\infty$ (respectively, $0\le   \tau_1^+ < \tau_2^+ \le \cdots \le \tau_k^+ \le \cdots \to +\infty$ ), one has
\begin{eqnarray} \label{2023.7.10-4} \left\{ \begin{array} {ll} \mu \nabla^* \nabla \mathbf{u}_k^\mp  -\mu \, \mbox{Ric} (\mathbf{u}^\mp_k)= \tau \mathbf{u}^\mp_k, \;\;&\mbox{in}\;\; \Omega,\\
2\mu \, (\mbox{Def}\, \mathbf{u}^\mp_k)^\# \boldsymbol{\nu} -\mu\, (\mbox{div}\, \mathbf{u}^\mp_k )\boldsymbol{\nu}=0, \;\;&\mbox{on}\;\; \partial \Omega.\end{array} \right.\end{eqnarray}
 where ${\mathbf{u}}_k^-\in  [H^1_0(\Omega)]^n$) (respectively,  ${\mathbf{u}}_k^+\in  [H^1(\Omega)]^n$) is the eigenvector corresponding to the elastic Dirichlet eigenvalue  $\tau_k^{-}$  (respectively, Neumann eigenvalue $\tau_k^+$).

Note that 
 \begin{eqnarray} \!\!\!\!\!\!\!\!\!\!\! &&\label{2023.7.10-20}\qquad \qquad \\ 
\!\!\!\!\!\!\!\!\!\!\!\!\! && \mu\nabla^*\nabla \mathbf{u}:=-\mu\sum\limits_{j=1}^n\! \Big\{\! \Delta_g u^j \! +2 \!\sum\limits_{k,s,l=1}^n\! g^{kl} \Gamma_{sk}^j \frac{\partial u^s}{\partial x_l} \!+\!
   \sum_{k,s,l=1}^n \!\Big(g^{kl} \frac{\partial \Gamma^j_{sl}}{\partial x_k}\! +\!\sum\limits_{h=1}^n g^{kl} \Gamma_{hl}^j \Gamma_{sk}^h\! -\!\sum\limits_{h=1}^n g^{kl} \Gamma_{sh}^j \Gamma_{kl}^h \Big)u^s\Big\}\frac{\partial}{\partial x_j},\nonumber\end{eqnarray}
which does not contain other second-order derivative terms except for $\Delta_g \mathbf{I}_n$. Therefore $\mu \nabla^* \nabla \mathbf{u}  -\mu \, \mbox{Ric} (\mathbf{u})$ is the (matrix-valued) Laplace-type operator on $\Omega$.  It can easily be verified that (\ref{2023.7.11-1})  just is the Dirichlet eigenvalue problem of the  (matrix-valued) Laplace-type operator. When taking $\mathbf{u}=(u, \cdots, u)$ for any $u\in H^2(\Omega)$,  then  (\ref{2023.7.10-1}) becomes  eigenvalue problem for the (matrix-valued) Laplace-type  operator with the free boundary condition because the free boundary condition becomes $((\nabla  u)\cdot \boldsymbol{\nu}, \cdots, (\nabla  u)\cdot \boldsymbol{\nu}  )=(0, 0, \cdots, 0)$. 
We  can easily find  that $\mathbf{q}_{-2} (x,\xi, \tau)= (\tau -\mu \,|\xi|^2)^{-1} \mathbf{I}_n$, and $\mathbf{q}_{-3}(x,\xi,\tau)$ is odd (vector-valued function) in $\xi \in \mathbb{R}^n$,  where $\iota \big((\tau {I}\! -\!\mathcal{P})^{-1})\,\sim \,\sum_{j\le -2} \mathbf{q}_j (x, \xi, \tau)  $, and $\mathcal{P}$ is the double operator of $\mu \nabla^* \nabla \mathbf{u}  -\mu \, \mbox{Ric} (\mathbf{u})$ on the double closed manifold $\mathcal{M}:=\Omega\cup (\partial \Omega)\cup \Omega^*$.
 It is completely similar to McKean-Singer's method (see, \cite{MS-67}) to show that  
\begin{eqnarray} \label{23.6.8-1}\sum_{k=1}^\infty e^{-t \tau_k^\mp} \sim  \frac{n}{(4\pi\mu)^{n/2}\, } t^{-n/2} \mp \frac{1}{4} \frac{n}{(4\pi \mu)^{(n-1)/2}} t^{-(n-1)/2} +O(t^{1-n/2}) \;\;\;\mbox{as}\;\; t\to 0^+.  \end{eqnarray} This is a classical result (cf. \cite{MS-67} or \cite{BrGi}). 
On the other hand, in this case, one has  $\alpha:=\frac{\mu}{\lambda+2\mu}\to 1^-$ (In fact, one can also take $\alpha=1$ in \cite{CaFrLeVa-23} because the Lam\'{e} operator is strongly elliptic for $\alpha=1$ and  uniformly strongly elliptic for all $\alpha \in [\alpha_0, 1]$, where $\;0<\alpha_0<1$). Obviously, according to (1.23)--(1.24) of \cite{CaFrLeVa-23}, the cube equation $R_\alpha (w)=0$ has six real roots $\gamma_R$ as the following:  $\,0,\,0, \,\sqrt{4+ 2\sqrt{2}}, \,-\sqrt{4+ 2\sqrt{2}}, \,\sqrt{4-2\sqrt{2}}, \,-\sqrt{4-2\sqrt{2}}$. It follows from (1.27) and (1.28) of Theorem 1.8 on p.$\,$8 in \cite{CaFrLeVa-23} that 
as $\alpha\to 1^-$,  
\begin{eqnarray} \label{23.6.8-5} && b^{-}= - \frac{ \mu^{\frac{1-n}{2}}}{ 2^{n+1} \pi^{\frac{n-1}{2}} \Gamma(\frac{n+1}{2})}\,n ,\\
&& \label{2024.5.5-3}  b^{+} =  \frac{ \mu^{\frac{1-n}{2}}}{ 2^{n+1} \pi^{\frac{n-1}{2}} \Gamma(\frac{n+1}{2})} 
\big(n-4+4\gamma_R^{1-n}\big).\end{eqnarray}

Now, suppose by contradiction that Theorem 1.8 of \cite{CaFrLeVa-23} is correct for Neumann (i.e., free) boundary condition. Then, by the corresponding coefficient (\ref{2024.5.5-3}) of the two-term asymptotic expansion   in \cite{CaFrLeVa-23}, one immediately obtains the following heat trace asymptotic expansion from  (\ref{2023.5.5-11}) and  (\ref{2023.5.5-6}):
  \begin{eqnarray}\label{23.5.5-81} && \sum_{k=1}^\infty e^{-t \tau_k^+} \sim  \frac{n}{(4\pi\mu)^{n/2}\, } t^{-n/2} \!+\! \Big(\frac{1}{4} \frac{n}{(4\pi \mu)^{(n-1)/2}}- 
 \frac{(4-4\gamma_R^{1-d} )\mu^{\frac{1-n}{2}}}{ 2^{n+1} \pi^{\frac{n-1}{2}} }   \! \Big)  t^{-(n-1)/2}\\
 && \qquad \qquad \quad \;\; +o(t^{-(n-1)/2}) \;\;\mbox{as}\;\; t\to 0^+, \nonumber \end{eqnarray}
which is different from the McKean-Singer's result (\ref{23.6.8-1}) for the Neumann (i.e., free) boundary condition (here an additional constant $- 
 \frac{(4-4\gamma_R^{1-d}) \mu^{\frac{1-n}{2}}}{ 2^{n+1} \pi^{\frac{n-1}{2}} } $ appears). {{\bf  This also implies that the conclusions in \cite{CaFrLeVa-23}  and \S6.3 of \cite{ SaVa-97} are wrong}.}

\vskip 1.29 true cm

\section{A counterexample to \cite{SaVa-97} and \cite{CaFrLeVa-23} in the case of the Dirichlet boundary condition}

\vskip 0.58 true cm

We first show that the $1$-form representation of the Lam\'{e} operator $P_g$ (defined on a Riemannian manifold $(\Omega, g)$) is just the generalized Ahlfors Laplacian. The Ahlfors Laplacian originated from the conformal geometry and was introduced by Ahlfors in 1974 and 1976 (see \cite{Ahl-74} and \cite{Ahl-76}).     

 Let $T\Omega$ and $T^*\Omega$ be the tangent and cotangent bundles of the Riemannian $n$-manifold $(\Omega,g)$, respectively. The space of all $C^\infty$ vector fields will be denoted by $\mathscr{X}$. 
Recall that in the terms of vector fields, the Lam\'{e} operator $ P_g$ defined on $(\Omega, g)$  can be written as (cf. (\ref{1-1})): \begin{eqnarray} \label{23.9.22-2} P_g \mathbf{u}:\!\!\!\!\!\!&& \!\!\!\!= \mu \big(\nabla^* \nabla \mathbf{u}\big)  -(\mu +\lambda) \,\mbox{grad}\; \mbox{div}\, \mathbf{u} -\mu \, \mbox{Ric}(\mathbf{u})\\
\!\!\!\!\!\!&&\!\!\!\!= \sum_{j=1}^n \Big[ \sum_{k=1}^n\Big(-\mu\nabla_k\nabla^k u^j -(\lambda+\mu) \nabla^j\nabla_k u^k - \mu\,\mbox{Ric}^{j}_{\;k}\, u^k\Big)\Big]\frac{\partial }{\partial x_j}\nonumber \\
\!\!\!\!\!\!&&\!\!\!\!=\sum_{j=1}^n \Big[ \sum_{k=1}^n\Big( - \mu \nabla_k \nabla^k  u^j -\lambda \nabla^j \nabla_ku^k -\mu \nabla_k \nabla^j u^k\Big)\Big]\frac{\partial}{\partial x_j}, \;\;\;\; \mathbf{u}=\sum_{j=1}^n u^j \frac{\partial}{\partial x_j} \in  \mathscr{X},\nonumber\end{eqnarray}
 where $\nabla^*\nabla \mathbf{u}$ is given by (\ref{2023.7.10-20}),  $\nabla_k u^j = u^{j}_{\,\;; k}$  are the components of the covariant derivative  of the vector field $\mathbf{u}$,  and $\nabla^k=\sum_{l=1}^n g^{kl}\nabla_l$.  
 On the other hand,  the  Lam\'{e} operator as well as the associated boundary value problems can equivalently be discussed in the language of $1$-form.  
  If   $\alpha$ is  the $1$-form dual to the vector field $\mathbf{u}=\sum_{j=1}^n u^j \frac{\partial}{\partial x_j}$ in the sense that 
 \begin{eqnarray*} \alpha(\mathbf{X})=g(\mathbf{u}, \mathbf{X}), \;\; \mathbf{X}\in \mathscr{X}, \end{eqnarray*}
 then  \begin{eqnarray} \label{23.9.22-5} \alpha= \sum_{j=1}^n u_j dx_j \;\;\, \mbox{and}\;\; u_j=\sum_{l=1}^n g_{jl}u^l. \end{eqnarray}
   Noting  that  $\nabla g=0$,  one has,  in  index notation, 
   \begin{eqnarray*}&&\sum_{k=1}^n \nabla_k \nabla^k u^j= \sum_{k,l=1}^n \nabla_k\nabla^k g^{jl} u_l 
  =
  \sum_{k,l,m=1}^n  g^{jl} g^{km} \nabla_k \nabla_m u_l
   =\sum_{l,m=1}^n g^{jl}\nabla^m \nabla_m u_l ,\\
 && \sum_{k=1}^n \nabla^j \nabla_ku^k=\sum_{k,l=1}^n  \nabla^j \nabla_k  g^{kl} u_l=\sum_{l=1}^n \nabla^j \nabla^l u_l =\sum_{l,m=1}^n  g^{jm}\nabla_m \nabla^l u_l,\\
 &&\sum_{k=1}^n   \nabla_k \nabla^j u^k   = \sum_{k,l=1}^n \nabla_k \nabla^j g^{kl} u_l =   \sum_{k,l,m=1}^ng^{kl} \nabla_k g^{jm} \nabla_m u_l=\sum_{l,m=1}^n g^{jm} \nabla^l \nabla_m u_l.\end{eqnarray*} 
Combining the above facts and the last line of (\ref{23.9.22-2}) we find   that for $\mathbf{u}\!=\!\sum_{j=1}^n \!u^j \frac{\partial}{\partial x_j} \in  \mathscr{X}$, 
\begin{eqnarray} \label{23.9.21-1} && \;P_g \mathbf{u} = \sum_{j=1}^n \Big[ \sum_{l,m=1}^n \!\Big(\!- \mu  g^{jl}\nabla^m \nabla_m u_l  -\lambda g^{jm}\nabla_m \nabla^l u_l-\mu g^{jm} \nabla^l \nabla_m u_l\Big)\Big] \frac{\partial }{\partial x_j}\\
&& \quad \quad \; = \sum_{j=1}^n \Big[ \sum_{l,m=1}^n \!g^{jm}\Big(\!- \mu  \nabla^l \nabla_l u_m  -\lambda \nabla_m \nabla^l u_l-\mu \nabla^l \nabla_m u_l\Big)\Big] \frac{\partial }{\partial x_j}\nonumber\\
&& \quad \quad \;:=\sum_{j=1}^n  \phi^j \frac{\partial}{\partial x_j}, \nonumber\end{eqnarray}
where $$\phi^j:= \sum_{l,m=1}^n \!g^{jm}\Big(\!- \mu  \nabla^l \nabla_l u_m  -\lambda \nabla_m \nabla^l u_l-\mu \nabla^l \nabla_m u_l\Big).$$
 Therefore,  we have  that 
  \begin{eqnarray} \label{23.9.22-3}  && P_g^{\flat} (\alpha) : =\big( P_g\mathbf{u}\big)^\flat= \sum_{k=1}^n \psi_k \,dx_k, \end{eqnarray}
  where 
  \begin{eqnarray} &&\psi_k =\sum_{j=1}^n g_{kj}\phi^j =  \sum_{j=1}^n \Big[ \!\sum_{l,m=1}^n   g_{kj} g^{jm}\Big(\!- \mu  \nabla^l \nabla_l u_m  -\lambda \nabla_m \nabla^l u_l-\mu \nabla^l \nabla_m u_l\Big)\Big] \\
 && \quad\;\,  =  \sum_{l=1}^n \!\Big(\!- \mu  \nabla^l \nabla_l u_k  -\lambda \nabla_k \nabla^l u_l-\mu \nabla^l \nabla_k u_l\Big) \nonumber\\
  &&\quad\;\,  = \sum_{l=1}^n\Big( -(\lambda +2\mu) \nabla_k \nabla^l u_l +\mu \nabla^l \nabla_k u_l -\mu \nabla^l \nabla_l u_k -2\mu \,\mbox{Ric}^l_{\, k}u_l\Big),\nonumber
   \end{eqnarray}
    $\alpha$ is given by (\ref{23.9.22-5}),  and  $\flat$  is the flat operator (for a vector field) by lowering an index. 
 Hence \begin{eqnarray} \label{23.9.23-10} P^\flat_g (\alpha) = \sum_{k=1}^n \Big[  \sum_{l=1}^n\Big( -(\lambda +2\mu) \nabla_k \nabla^l u_l +\mu \nabla^l \nabla_k u_l -\mu \nabla^l \nabla_l u_k -2\mu \,\mbox{Ric}^l_{\, k}u_l\Big)\Big] dx_k.\end{eqnarray} 
   Let $d: \Lambda  T^*\Omega\to \Lambda T^*\Omega$ be the exterior differential operator, where $\Lambda T^*\Omega=\bigoplus_{p=0}^n \Lambda^p T^*\Omega$.    The adjoint operator $\delta$ of  $d$ acting on a $p$-form $\alpha$ is defined in terms of $d$ and the Hodge star operator by formula 
 \begin{eqnarray*} \delta \alpha =(-1)^{np+n+1} * d* \alpha,\end{eqnarray*}  and the Hodge star operator $*: \Lambda^p T^*\Omega \to \Lambda^{n-p}T^* \Omega,\;  p=0, \cdots, n$, is defined by \begin{eqnarray*}  \langle \gamma, \eta\rangle =\gamma \wedge  *\eta \end{eqnarray*} 
 for any $\gamma,\eta\in \Lambda^pT^*\Omega$.
  It is well known (see, for example, p.$\,$16--17 of \cite{CLN-06})  that $\sum_{l=1}^n \nabla^l u_l= \mbox{div}\, \alpha= -\delta \alpha$ 
   for the $1$-form $\alpha=\sum_{j=1}^n u_j dx_j$.
   Furthermore,  from Exercise 5 on p.$\,$561 in \cite{Ta3} we know that $\delta d \alpha = \sum_{k,l=1}^n \big(\nabla^l \nabla_k u_l-  \nabla^l\nabla_l u_k\big) dx_k$ for the above $1$-form $\alpha$.  Combining these facts and   (\ref{23.9.23-10}) we  obtain that
 \begin{eqnarray} \label{23.9.21-2} P_g^{\flat} (\alpha) = (\lambda+2\mu) d\delta \alpha + \mu \delta d \alpha - 2\mu \, \mbox{Ric}\, (\alpha).\end{eqnarray}
 Here $\mbox{Ric\,}(\alpha)$ denotes the Ricci action on $1$-forms $\alpha$:
\begin{eqnarray*} \mbox{Ric} \,(\alpha) =\mbox{Ric}(\mathbf{u}, \cdot)=\sum_{k=1}^n \big(\sum_{l=1}^n \mbox{Ric}^l_{\,k} u_l \big)\, dx_k,\end{eqnarray*} 
where $\mathbf{u}$ is the vector field dual to $\alpha$.
(\ref{23.9.21-2}) is the $1$-form representation of the Lam\'{e} operator. 

  The {\it Ahlfors operator} $S$ and the {\it Ahlfors Laplacian}  
 $ L=S^*S:= \left(1-\frac{1}{n}\right) d\delta+ \frac{1}{2} \delta d  -\mbox{Ric}\,$ on $\,\Lambda^1T^*\Omega$ play a fundamental role in the study of quasiconformal geometry (see, for example, \cite{Ahl-74}, \cite{Ahl-76}  or \cite{PiOr-96}). 
 The Ahlfors operator $S$ acting on $1$-form $\alpha$ is defined by  (see \cite{BGOP}) 
 \begin{eqnarray*} S\alpha =\nabla^s \alpha +\frac{1}{n} \delta \alpha \cdot g, \end{eqnarray*}
 where $\nabla^s$ denotes the symmetrized $\nabla$:
 \begin{eqnarray*} (\nabla_\alpha^s) (\mathbf{X}, \mathbf{Y})=\frac{1}{2} \Big( (\nabla_{\mathbf{X}}\alpha\big) (\mathbf{Y})+ (\nabla_{\mathbf{Y}} \alpha )(\mathbf{X})\Big), \;\;\, \;\mathbf{X}, \mathbf{Y}\in \mathscr{X}.\end{eqnarray*} 
  The {\it generalized Ahlfors Laplacian}  (see,  \cite{BGOP})  is defined as the operator $P=ad\delta +b\delta d -\epsilon \rho$ on $\Lambda^1 T^*\Omega$, where $a$ and $b$ are positive constants and where $\epsilon \rho$ is an arbitrary constant multiple of the Ricci tensor.  
 Clearly, from (\ref{23.9.21-2}) we see that  the $1$-form representation $P_g^\flat$ of the elastic Lam\'{e} operator $P_g$  is just the {\it generalized Ahlfors Laplacian} on $\Lambda^1T^* \Omega$ when $\lambda+2\mu=a>0$ and $\mu=b>0$.

 In \cite{PiOr-96},  A.  Pierzchalski and B.  {\O}rsted considered  the two-term heat trace asymptotic expansions for  the  Ahlfors Laplacian $L=S^*S$ with the following three boundary conditions: 
  \begin{eqnarray} \label{23.9.19-1} \mathscr{X}_D =\{\mathbf{Z}\in \mathscr{X} \big| \mathbf{Z}_p=0\;\; \mbox{for all}\;\, p\in \partial \Omega\}\end{eqnarray}
for the boundary problem of the Dirichlet type $(D)$. 
Also, define \begin{eqnarray} \label{23.9.19-2}  \mathscr{X}_N=\{\mathbf{Z}\in \mathscr{X} \big|{Z}_p^N=0 \;\, \mbox{and}\,\;  (\nabla_N\mathbf{Z})^T_p =0\;\;\mbox{for all}\;\; p\in \partial \Omega\}\end{eqnarray} 
for the boundary problem of the Neumann type $(N)$. 
In coordinates, the condition ($N$) reads: 
\begin{eqnarray*} \label{2.9.19-3} {Z}^n=0\;\; \mbox{and}\;\; \frac{\partial}{\partial r} Z^k =(A_N\mathbf{Z})^k\;\;\, (K=1, \cdots, n-1)\;\;\, \mbox{on}\;\; \,\partial \Omega.\end{eqnarray*}
Finally, define 
\begin{eqnarray}\mathscr{ X}_E=\{\mathbf{Z}\in \mathscr{X}\big| \mathbf{Z}_p^T =0\;\;\mbox{and}\;\; \mbox{div}\; \mathbf{Z}_p=0\;\;\mbox{for all}\;\, p\in \partial \Omega\}\end{eqnarray}
for the boundary problem of the theory of elasticity ($E$). In coordinates, we  obtain 
\begin{eqnarray*} \label{ 23.9.19-4}  Z^1=\cdots = Z^{n-1} =0\;\;\, \mbox{and}\;\; \, \frac{\partial}{\partial r} Z^n=- Z^n \, \mbox{tr}\; h,\end{eqnarray*} 
$h$ is the second fundamental form of $\partial \Omega$.

Each of the boundary conditions of type $D$, $N$, and $E$ is self-adjoint in the sense: if $\mathbf{Z}_1,\mathbf{Z}_2\in {\mathscr{X}}_B:=\{\mbox{either}\;\, \mathscr{X}_D, \mbox{or}\,\, \mathscr{X}_N, \mbox{or}\,\, \mathscr{X}_E\}$, then $(L\mathbf{Z}_1, \mathbf{Z}_2)=(\mathbf{Z}_1,L\mathbf{Z}_2)$. 
Obviously, each of the boundary conditions $N$ and  $E$ is different from our free boundary condition (see the definition in Section 1).

\vskip 0.2 true cm 

 By using a technique of invariant theory, which was entirely different from our method in \cite{Liu-21},  A. Pierzchalski and B.  {\O}rsted  in 1996 showed the following result:

\vskip 0.22 true cm
 
\noindent {\bf Theorem  8.1 (see Theorem 5.1 of \cite{PiOr-96}). \    {\it Let $L$ be the self-adjoint elliptic extension of the Ahlfors Laplacian with boundary conditions $D$, $N$, or $E$ on the manifold $M$. Then the small-time asymptotic of the heat kernel have the first two terms as follows.

 \vskip 0.1 true cm 
  \mbox{Case D}: 
    \begin{eqnarray*}  \mbox{tr}\; e^{-tL} \!\!\!\!&& \!\!\!\!\! \sim (4\pi t)^{-n/2} \cdot \mbox{vol}(\Omega) \cdot \big[ a^{-n/2} +(n-1)b^{-n/2}\big] \\
\!\!\!\! &&\;\; - \frac{1}{4} (4\pi t)^{-(n-1)/2} \cdot \mbox{vol}(\partial \Omega)\cdot \big[ a^{-(n-1)/2 } + (n-1)b^{-(n-1)/2}\big]. \end{eqnarray*}
 
   \mbox{Case N}: 
    \begin{eqnarray*}  \mbox{tr}\; e^{-tL} \!\!\!\!&& \!\!\!\!\! \sim (4\pi t)^{-n/2} \cdot \mbox{vol}(\Omega) \cdot \big[ a^{-n/2} +(n-1)b^{-n/2}\big] \\
 \!\!\!\! &&\;\; + \frac{1}{4} (4\pi t)^{-(n-1)/2} \cdot \mbox{vol}(\partial \Omega)\cdot \big[ a^{-(n-1)/2 } + (n-3)b^{-(n-1)/2}\big]. \end{eqnarray*}

    \mbox{Case E}: 
    \begin{eqnarray*}  \mbox{tr}\; e^{-tL}  \!\!\!\!&& \!\!\!\!\! \sim (4\pi t)^{-n/2} \cdot \mbox{vol}(\Omega) \cdot \big[ a^{-n/2} +(n-1)b^{-n/2}\big] \\
 \!\!\!\! &&\;\; - \frac{1}{4} (4\pi t)^{-(n-1)/2} \cdot \mbox{vol}(\partial \Omega)\cdot \big[ a^{-(n-1)/2 } + (n-3)b^{-(n-1)/2}\big]. \end{eqnarray*}
 Here $a, b$ are the values in the Ahlfors Laplacian: $a= (n-1)/n$ and  $b=\frac{1}{2}$. (Actually, in Case $E$ they could be arbitrary positive; $L$ would still be self-adjoint and elliptic.) }
 
\vskip 0.39 true cm

\noindent {\bf Remark 8.2}. \    In the Lam\'{e} operator, if  we take $\mu=b=\frac{1}{2}$ and $\mu+2\lambda=a= 1-\frac{1}{n}$, then,  by applying Theorem 8.1 (i.e.,  a result  of A. Pierzchalski and  B. {\O}rsted in \cite{PiOr-96})   we see that the two-term asymptotic expansion for the elastic heat trace with Dirichlet boundary condition is the same as that in (\ref{1-7}) of Theorem 1.1 (i.e., the conclusion of \cite{Liu-21}).  But the conclusion for the Dirichlet boundary condition  in \cite{CaFrLeVa-23} (or p.$\,$237 of \S6.3 in \cite{SaVa-97})  is different from the second term in Case D of Theorem 8.1.  A  superfluous integral term appears in the second term  of \cite{CaFrLeVa-23} (or p.$\,$237 of \S6.3 in \cite{SaVa-97}).  A. Pierzchalski and B.  {\O}rsted's conclusion  for Dirichlet boundary condition (i.e., Case D of Theorem 8.1) is certainly correct because it is based on the invariant theory. This also implies that the conclusion in \cite{CaFrLeVa-23} and  \S6.3 of \cite{SaVa-97}  is completely incorrect.

 \vskip 0.39 true cm

\noindent {\bf Remark 8.3}. \  (i) \  This serious error in \cite{CaFrLeVa-23} or \S6.3 of \cite{SaVa-97} for the elastic spectral asymptotics formula stem from the so-called  ``algorithm'' theory  in \cite{SaVa-97}, which, throughout the whole book \cite{SaVa-97},
is essentially wrong. It is impossible to correct these kinds of fundamental errors at all.

(ii) \ Clearly, most conclusions in this book \cite{SaVa-97} are wrong because they are based
on an erroneous ``algorithm'' theory. The book \cite{SaVa-97}, which is theoretically wrong,  has misled a large number of readers for twenty-six years.} 

 \vskip 0.39 true cm 
 
\noindent {\bf Remark 8.4}. \  In \cite{BGOP},  T. P. Branson, P. B. Gilkey, B. {\O}rsted and A. Pierzchalski gave the first three terms
  of asymptotic expansion of the heat trace for the generalized Ahlfors Laplacian with absolute or relative boundary conditions. These two boundary conditions originated from geometry,  each of which are different from (D),  (N), (E) and free boundary condition.    

\vskip 1.19 true cm

\section{Two remarks and two conjectures}

\vskip 0.58 true cm

We will first give two remarks for erroneous statements in \cite{CaFrLeVa-23} (see also  \cite{CaFrLeVa-22}):

\vskip 0.35 true cm

\noindent {\bf Remark 9.1}. \    {\it In the footnote 3 on p.\,3 of {\cite{CaFrLeVa-23}}, the authors of {\cite{CaFrLeVa-23}} wrote: \textcolor{blue}{ ``We will not call (1.6) the Neumann condition in order to avoid confusion  with erroneous ``Neumann'' condition in \cite{Liu-21}.'}' And for the remark 1.12, from  line -8 to line -7, on  p.\,10 of {\cite{CaFrLeVa-23}}, \textcolor{blue}{``Therefore, it is hard to assign a meaning to
Liu's result in this case [13, Theorem 1.1, the lower sign version of formula (1.10)].''}}

\vskip 0.20 true cm

Sometimes, a system of physical equations can be written as the differential $1$-form equation, it can be equivalently written as a vector field equation. Similarly, one can equivalently write the vector field (i.e., a tensor of type $(1,0)$) as a differential $1$-form (i.e., a tensor of type $(0,1)$). For the zero Neumann (i.e., free) boundary, they are equivalent when this boundary condition is  written as a tensor of type $(1,0)$, or a tensor of type $(0,1)$, or a tensor of type $(1,1)$. On one hand, the authors of \cite{CaFrLeVa-23} plagiarized (did not cite) the tensor expression of type  $(1,1)$ (or $(1,0)$) for the Lam\'{e} operator as well as  Neumann (i.e., free) boundary condition of \cite{Liu-19} on a Riemannian manifold and then revised Neumann (i.e., free) boundary into the vector field expression. On the other hand, the authors of \cite{CaFrLeVa-23}  misinterpreted our definition (see \cite{Liu-19} or \cite{Liu-21}) of the elastic Neumann boundary condition. In fact, in \cite{Liu-21} the Neumann boundary condition is $2\mu (\mbox{Def}\,\mathbf{u})^\#\boldsymbol{\nu} +\lambda(\mbox{div}\, \mathbf{u})\boldsymbol{\nu}$ on $\partial \Omega$ (we simply denoted it as $\frac{\partial \mathbf{u}}{\partial \boldsymbol{\nu}}$), as clearly pointed out from line 18 to line 21 on p.\,10166 of \cite{Liu-21}, and the author of \cite{Liu-21} wrote ``For the derivation of the Navier--Lam\'{e} elastic wave equations, its mechanical meaning, and the explanation of the Dirichlet and Neumann boundary conditions, we refer the reader to \cite{Liu-19} for the case of Riemannian manifold and ....''. In \cite{Liu-19}, the Neumann boundary condition is $2\mu (\mbox{Def}\, \mathbf{u})^\#\boldsymbol{\nu} + \lambda(\mbox{div}\, \mathbf{u})\boldsymbol{\nu}$ on $\partial \Omega$. 
In line 9-13, p.$\,$10165 of \cite{Liu-21}, we also wrote: ``We denote by $P_g^-$ and $P_g^+$ the Navier-Lam\'{e} operators with the Dirichlet and Neumann boundary conditions, respectively.
    Since $P_g^-$ (respectively, $P_g^+$) is an unbounded, self-adjoint and positive (respectively,  nonnegative) operator in $[H^1_0(\Omega)]^n$ (respectively, $[H^1(\Omega)]^n$) with discrete spectrum $0< \tau_1^- < \tau_2^- \le \cdots \le \tau_k^- \le \cdots \to +\infty$ (respectively, $0\le \tau_1^+ < \tau_2^+ \le \cdots \le \tau_k^+ \le \cdots \to +\infty$),....''
 \cite{Liu-19} is an earlier paper  which was posed on arXiv on Aug.\,14, 2019 by the author (Note that the paper \cite{Liu-21} was submitted to {\it The Journal Geometric Analysis} on Jul.\,21, 2020).
 It is very clear from Lemma 2.1.1 of \cite{Liu-19} that $\frac{\partial \mathbf{u}}{\partial \boldsymbol{\nu}}:=2\mu (\mbox{Def}\, \mathbf{u})^\#\boldsymbol{\nu} +\lambda(\mbox{div}\, \mathbf{u})\boldsymbol{\nu}$, which equals to
$ \sum_{j,k=1}^n \big( \mu (u^j_{\;\,;k}+ u_k^{\;\,;j}) \nu^k +\lambda u^k_{\;\, ;k} \nu^j\big) \frac{\partial}{\partial x_j}$
(or equivalently, $\sum_{k=1}^n\big(\lambda  n^j \nabla_k u^k +\mu ( n^k \nabla_k u^j+ n_k \nabla_j u^k)\big)$), where $\{\frac{\partial}{\partial x_j}\}_{j=1}^n$ is the nature coordinate basis.

{{\bf Obviously, Remark 1.12, on p.$\,$10 in \cite{CaFrLeVa-23} is completely wrong. It seems that the authors of \cite{CaFrLeVa-23} attempted to mislead the reader (see also Response 2 of \cite{Liu-22b}).}}

\vskip 0.49 true cm

\noindent {\bf Remark 9.2}. \  {\it From line 6 to line 21 on p.$\,$34 of \cite{CaFrLeVa-23}, the authors of \cite{CaFrLeVa-23} wrote: \textcolor{blue}{``Let us conclude this appendix with a brief historical account. We note that the
expression for $\tilde{b}_{\text{Dir}}$  was already found${}^{19}$ in the 1960 paper by M. Dupuis, R. Mazo, and L.Onsager${}^{20}$ [6]. Remarkably, this paper includes the critique of the 1950 paper by
E. W. Montroll who presented {\it exactly} Liu's expression (1.32) for the second asymptotic coefficient, modulo some scaling, see [16, formulae (3)--(5)]. Dupuis, Mazo, and
Onsager wrote, we quote: ``Montroll ...pointed out in 1950 a defect in the usual counting process of the normal modes of vibration and derived a corresponding correction
term for the Debye frequency spectrum, ... proportional to the area of the solid.....longitudinal modes.''}}

\vskip 0.16 true cm

Obviously, the authors of \cite{CaFrLeVa-23} have misunderstood all results in these papers mentioned above, and have  misled the reader by completely wrong statements. In \cite{Mo-50},
in order to discuss the effect of the volume and surface area to the heat capacities in low temperature,
E. W. Montroll gave the expression for the counting function of eigenvalues for the elastic normal modes of a three-dimensional rectangular solid $\{(x,y,z)\in \mathbb{R}^3\,\big|\, 0\le x\le L_x, 0\le y\le L_y, 0\le z\le L_z\}$ with the boundary condition $\boldsymbol{\nu} \cdot \nabla \mathbf{u}$ (according to the explanation of the boundary condition in \cite{CaFrLeVa-22}). In \cite{Mo-50}, the considered domain is a very special rectangular domain, the boundary conditions neither Dirichlet boundary condition nor the free boundary condition. In \cite{DuMaOn-60}, M. Dupuis, R. Mazo and L. Onsager investigated an isotropic solid at low temperatures whose model is a rectangular plate of thickness $l_3$ and other dimensions
$l_1$  and $l_2$ with realistic boundary conditions. The faces parallel to the plane of the plate are
supposed to be free of stresses, whereas the periodic
boundary conditions are given on the other faces (see, p.$\,$1453 of \cite{DuMaOn-60}).
  Clearly, these boundary conditions are completely different from the (whole) Dirichlet or Neumann (i.e., free) boundary conditions in \cite{Liu-21}. In \cite{Mo-50} and \cite{DuMaOn-60}, the authors respectively calculated the counting functions with different boundary conditions for a very special domain (i.e., three-dimensional rectangular solid) by elementary calculations. The domains and the boundary conditions in three papers (\cite{Liu-21}, \cite{Mo-50} and \cite{DuMaOn-60}) are quite different.
 In \cite{Liu-21}, we proved the asymptotic formulas for the heat traces of an elastic body   for {\it general compact smooth manifold with smooth boundary} (for the corresponding Dirichlet and the Neumann (i.e., free) boundary conditions, respectively).

 It seems that the authors of  \cite{CaFrLeVa-23} always confuse the huge differences or essentially technical difficulties (of a mathematical problem) among different domains (for example, for a rectangle domain, for a bounded domain in the Euclidean space, and for a bounded domain in a Riemannian manifold, etc). 

\vskip 0.39 true cm

Clearly, the following two conjectures are still open:
 \vskip 0.1 true cm
 {\bf  Conjecture 9.3}:  \ Let $\Omega$ be a compact, connected smooth Riemannian $n$-manifold with smooth boundary $\partial \Omega$.
 The following two-term asymptotics hold:
 \begin{eqnarray} \label{2022.9.7-1} &&\;\;\, \, \mathcal{N}_{\mp} (\Lambda) \sim \frac{1}{ \Gamma(1+\frac{n}{2})} \Big( \frac{n-1}{(4\pi \mu)^{n/2}} +\frac{1}{(4\pi (\lambda+2\mu))^{n/2}}\Big)\mbox{vol}_n(\Omega) \Lambda^{n/2}\\
&& \quad \;  \mp \frac{1}{4 \,\Gamma(1+\frac{n-1}{2})} \Big( \frac{n-1}{(4\pi\mu)^{(n-1)/2}} +\frac{1}{(4\pi(\lambda+2\mu))^{(n-1)/2}}\Big)\mbox{vol}_{n-1}(\partial \Omega) \,\Lambda^{(n-1)/2}\; \quad \;\;  \mbox{as}\;\, \Lambda \to +\infty, \nonumber \end{eqnarray}
 where $\mathcal{N}_{-} (\Lambda):=\#\{ k\big| \tau^{-}_k < \Lambda\}$ (respectively, $\mathcal{N}_{+} (\Lambda):=\#\{ k\big| \tau_k^{+}< \Lambda\}$) is the counting function of elastic eigenvalues for the Dirichlet (respectively, Neumann) boundary condition.

\vskip 0.1 true cm
 When $\lambda+\mu=0$, Conjecture 1 just is the famous Weyl conjecture, which has about 110-year history and remains open.

  \vskip 0.2 true cm

 {\bf Conjecture 9.4}:  Let $\Omega$ be a compact, connected smooth Riemannian $n$-manifold with smooth boundary $\partial \Omega$. Suppose that the corresponding elastic billiards is neither dead-end nor absolutely periodic for the $\Omega$. The above two-term asymptotics (\ref{2022.9.7-1}) hold.

 \vskip 0.2 true cm

 Obviously,  for any  given compact smooth Riemannian manifold $\Omega$ with smooth boundary, under the elastic billiards condition, {\bf an explicit two-term  asymptotic expansion of the counting function $\mathcal{N}_{\mp}(\Lambda)$ as $\Lambda\to +\infty$  had not been given in  \cite{SaVa-97} or \cite{CaFrLeVa-23}. 
{The erroneous algorithm for the second coefficients of the asymptotics expansion in \cite{CaFrLeVa-23} stemmed from the earlier erroneous approaches and wrong methods  in} \textcolor{blue}{[Va-84, \S6]}, \cite{Va-86}, \cite{SaVa-97} and \cite {CaVa-22} (cf. line -8 to line -1 from bottom, p.$\,$12 in \cite{CaFrLeVa-23}).}
In [De-12] and  [DuMaOn-60], the counting functions of elastic eigenvalues were explicitly calculated only for very special domain (a three-dimensional rectangular solid) with special boundary conditions.
 However, in Theorem 1.1 of \cite{Liu-21}, our two-term asymptotic formulas (\ref{1-7}) of the heat traces hold for any given compact connected smooth Riemannian manifold $\Omega$ with smooth boundary (do not need any additional assumption). Obviously, in Theorem 1.8 of \cite{CaFrLeVa-23},   the following key assumption has been lost:  $(\Omega,g)$ is such the corresponding billiards is neither dead-end nor absolutely periodic.

 \vskip 0.15 true cm

     However, not every manifold satisfies the billiards condition. For example, when $\lambda+\mu=0$ (the corresponding elastic operator becomes the Laplacian), the semi-sphere $\mathbb{S}^{n-1}_{+}:=\{(x_1,\cdots, x_{n-1}$, $x_n)\in \mathbb{R}^n \big|x_1^2 +\cdots +x_{n-1}^2+x_n^2=1,\; x_n\ge 0\}$ does not satisfy the billiards condition because every point $(x,\xi)\in S^* (\mathbb{S}^{n-1}_{+})$  is periodic  (Here $S^* (\mathbb{S}^{n-1}_{+})$ is the fiber of cotangent unit sphere over $\mathbb{S}^{n-1}_{+}$).
\vskip 0.1 true cm

\vskip 1.16 true cm
\section{Appendix}

\vskip 0.55 true cm
For the sake of convenience and completeness, in this appendix we give the proof of Theorem 1.1. The proof is the same as in \cite{Liu-21} except for some additional explanations in \cite{Liu-22c}.

\vskip 0.46 true cm

 \noindent  {\it Proof of Theorem 1.1.} \  From the theory of elliptic operators (see \cite{GiTr},  \cite{Mo3}, \cite{Pa},  \cite{Ste}), we see that the Navier-Lam\'{e} operator $-P_g$ can generate  strongly continuous semigroups $(e^{-tP_g^\mp})_{t\ge 0}$ with respect to the zero Dirichlet and zero  traction boundary conditions, respectively, in suitable spaces of vector-valued functions (for example, in $[C_0(\Omega)]^n$ (by Stewart \cite{Ste}) or in $[L^2(\Omega)]^n$ (by Browder \cite{Brow})), or in $[L^p(\Omega)]^n$ (by Friedman \cite{Fri}). Furthermore,
   there exist  matrix-valued functions ${\mathbf{K}}^\mp (t, x, y)$, which are called the integral kernels, such that (see \cite{Brow} or p.$\,$4 of \cite{Fri})
        \begin{eqnarray*}  e^{-tP^\mp_g}{\mathbf{w}}_0(x)=\int_\Omega {\mathbf{K}}^\mp(t, x,y) {\mathbf{w}}_0(y)dy, \quad \,
        {\mathbf{w}}_0\in  [L^2(\Omega)]^n.\end{eqnarray*}

Let $\{{{u}}_k^\mp\}_{k=1}^\infty$ be the orthonormal eigenvectors of the elastic operators $P_g^\mp$ corresponding to the eigenvalues $\{\tau_k^\mp\}_{k=1}^\infty$, then the integral kernels  ${\mathbf{K}}^\mp(t, x, y)=e^{-t P_g^\mp} \delta(x-y)$ are given by \begin{eqnarray} \label{18/12/18} {\mathbf{K}}^\mp(t,x,y) =\sum_{k=1}^\infty e^{-t \tau_k^\mp} {{u}}_k^\mp(x)\otimes {{u}}_k^\mp(y).\end{eqnarray}
This implies that the integrals of the traces of ${\mathbf{K}}^\mp(t,x,y)$ are actually  spectral invariants:
\begin{eqnarray} \label{1-0a-2}\int_{\Omega} \mbox{Tr}\,({\mathbf{K}}^\mp(t,x,x)) dV=\sum_{k=1}^\infty e^{-t \tau_k^\mp}.\end{eqnarray}

 We will combine calculus of symbols (see \cite{Gr-86}) and ``method of images'' to deal with asymptotic expansions for the integrals of traces of integral kernels.
Let $\mathcal{M}=\Omega \cup (\partial \Omega)\cup \Omega^*$ be the (closed) double of $\Omega$, and $\mathcal{P}$ the double to $\mathcal{M}$ of
 the  operator $P_g$ on $\Omega$.

Let us explain the double Riemannian manifold $\mathcal{M}$ and the differential operator $\mathcal{P}$ more precisely, and introduce how to get them from the given Riemannian manifold $\Omega$ and the Navier-Lam\'{e} operator $P_g$.
The double of $\Omega$ is the manifold $\Omega \cup_{\mbox{Id}} \Omega$, where $\mbox{Id}: \partial \Omega\to \partial \Omega$ is the identity map of $\partial \Omega$; it is obtained from $\Omega \sqcup \Omega$ by identifying each boundary point in one  copy of $\Omega$
with same boundary point in the other. It is a smooth manifold without boundary, and contains two regular domains diffeomorphic to $\Omega$  (see, p.$\,$226 of \cite{Lee}). When considering the double differential system $\mathcal{P}$ crossing the boundary, we make use of the coordinates as follows.  Let $x'=(x_1, \cdots, x_{n-1})$ be any local coordinates for $\partial \Omega$.  For each point $(x',0)\in \partial \Omega$, let $x_{n}$ denote the parameter along the unit-speed geodesic starting at $(x',0)$ with initial direction given by the inward boundary normal to $\partial \Omega$ (Clearly, $x_n$ is the geodesic distance from the point $(x',0)$ to the point $(x',x_n)$). In such coordinates $x_{n}>0$ in
 $\Omega$, and $\partial \Omega$ is locally characterized by $x_{n}=0$ (see, \cite{LU} or \cite{Ta-2}).
 Since the Navier-Lam\'{e} operator is a linear differential operator defined on $\Omega$, it can be further denoted as (see \ref{2023.5.22-12}) $P_g:=P(x, \{g^{jk}(x)\}_{1\le j,k\le n}$,
$\{\Gamma^j_{kl}(x)\}_{1\le j,k,l\le n}$, $\{\frac{\partial \Gamma^s_{jk}(x)}{\partial  x_l}\}_{1\le s,j,k,l\le n}$, $\{R^j_k(x)\}_{1\le j,k\le n}$, $\frac{\partial}{\partial x_1}, \cdots, \frac{\partial}{\partial x_{n-1}}$,$\frac{\partial }{\partial x_n})$.
  Let  $\varsigma: (x_1, \cdots, x_{n-1}, x_n)\mapsto (x_1, \cdots, x_{n-1}, -x_n)$ be the reflection with respect to the boundary $\partial \Omega$ in $\mathcal {M}$ (here we always assume $x_n\ge 0$).  Then  we can get the $\Omega^*$ from the given $\Omega$ and $\varsigma$.
    Now, we discuss the change of the metric $g$ from $\Omega$ to $\Omega^*$ by $\varsigma$. Recall that the Riemannian metric $(g_{ij})$ is given in the local coordinates $x_1, \cdots, x_n$, i.e., $g_{ij}(x_1,\cdots, x_n)$. In terms of the new coordinates $z_1,\cdots, z_n$, with  $x_i=x_i(z_1,\cdots, z_n), \,\, i=1,\cdots, n,$  the same metric is given by the functions $\tilde{g}_{ij} =\tilde{g}_{ij} (z_1, \cdots, z_n)$, where
\begin{eqnarray} \label{2022.10-2} \tilde{g}_{ij}= \sum_{k,l=1}^n \frac{\partial x_k}{\partial z_i} g_{kl} \frac{\partial x_l}{\partial z_j}.\end{eqnarray}
  If $\varsigma$ is a coordinate change in a neighborhood intersecting with $\partial \Omega$
  \begin{eqnarray}\label{2022.10.28-1} \left\{ \begin{array}{ll} x_1= z_1, \\
                   \cdots \cdots\\
                   x_{n-1}=z_{n-1},\\
                   x_n=-z_n,\end{array}\right.\end{eqnarray}
                                       then its Jacobian matrix is
\begin{eqnarray} \label{2022.11.8-1} J:=\begin{pmatrix} \frac{\partial x_1}{\partial z_1}& \cdots &   \frac{\partial x_1}{\partial z_{n-1}} & \frac{\partial x_1}{\partial z_{n}}\\
   \vdots & \ddots & \vdots & \vdots\\
  \frac{\partial x_{n-1}}{\partial z_{1}} & \cdots & \frac{\partial x_{n-1}}{\partial z_{n-1}} &  \frac{\partial x_{n-1}}{\partial z_{n}}\\
   \frac{\partial x_n}{\partial z_{1}}& \cdots & \frac{\partial x_n}{\partial z_{n-1}}&\frac{\partial x_n}{\partial z_{n}}\end{pmatrix} =   \begin{pmatrix} 1& \cdots &   0 & 0\\
   \vdots & \ddots & \vdots & \vdots\\
  0 & \cdots & 1 &  0\\
   0& \cdots & 0&-1\end{pmatrix}.\end{eqnarray}
Using this and (\ref{2022.10-2}), we immediately obtain the corresponding metric on the $\Omega^*$:
  (see \cite{MS-67},  or p.\,10169, p.\,10183 and p.\,10187 of \cite{Liu-21})
\begin{eqnarray} \label{2021.2.6-3}  g_{jk} (\overset{*}{x})\!\!\!&\!=\!&\!\!\!- g_{jk} (x) \quad \, \mbox{for}\;\;
  j<k=n \;\;\mbox{or}\;\; k<j=n,\\ g(\overset{*}{x}) \!\!\!&\!=\!&\!\!\! g_{jk} (x)\;\; \;\;\mbox{for}\;\; j,k<n \;\;\mbox{or}\;\; j=k=n,\\
  \label{2021.2.6-4}  g_{jk}(x)\!\!\!&\!=\!&\!\!\! 0 \;\; \;\mbox{for}\;\; j<k=n \;\;\mbox{or}\;\; k<j=n \;\;\mbox{on}\;\; \partial \Omega,\end{eqnarray} where $x_n(\overset{*}{x})= -x_n (x)$.
  We denote such a new (isometric) metric on $\Omega^*$ as $g^*$.
  It is easy to verify that \begin{eqnarray*} \begin{bmatrix} g_{11}(x) & \cdots & g_{1,n-1}(x)& -g_{1n}(x)\\
  \vdots & \ddots & \vdots & \vdots\\
  g_{n-1,1}(x) & \cdots & g_{n-1,n-1}(x) &- g_{n-1,n}(x)\\
 - g_{n1}(x) & \cdots & -g_{n,n-1}(x) & g_{nn}(x)\end{bmatrix}^{-1}=\begin{bmatrix} g^{11}(x) & \cdots & g^{1,n-1}(x)& -g^{1n}(x)\\
  \vdots & \ddots & \vdots & \vdots\\
  g^{n-1,1}(x) & \cdots & g^{n-1,n-1}(x) &- g^{n-1,n}(x)\\
 - g^{n1}(x) & \cdots & -g^{n,n-1}(x) & g^{nn}(x)\end{bmatrix}, \end{eqnarray*}
  where $[g^{jk}(x)]_{n\times n}$ is the inverse of $[g_{jk}(x)]_{n\times n}$.
    In addition, by this reflection $\varsigma$,
the differential operators $\frac{\partial }{\partial x_1}$, $\cdots$, $\frac{\partial}{\partial x_{n-1}}$, $\frac{\partial }{\partial x_n}$ (defined on  $\Omega$) are changed to $\frac{\partial }{\partial x_1}$, $\cdots$, $\frac{\partial}{\partial x_{n-1}}$, $-\frac{\partial }{\partial x_n}$ (defined on  $\Omega^*$), respectively. 
 It is easy to verify that
\begin{eqnarray}  \label{2023.2.22-8}\Gamma^{j}_{kl}(\overset{*}{x})=a_{jkl} \Gamma^{j}_{kl}(x),\end{eqnarray}
where \begin{eqnarray*} a_{jkl} =\left\{ \begin{array}{ll} 1 \;\;\; & \mbox{if there is no}\;\; n \;\; 
\mbox{among} \;\; j,k,l,\\
-1 \;\;\;\; & \mbox{if there is an}\;\; n \;\; 
\mbox{among} \;\; j,k,l,\\
1 \;\;\;\; & \mbox{if there are two}\;\; n \;\; 
\mbox{among} \;\; j,k,l,\\
-1 \;\;\; \; & \mbox{if there are three}\;\; n \;\; 
\mbox{among} \;\; j,k,l, \end{array}\right. \end{eqnarray*}
 \begin{eqnarray}\label{2022.5.22-9}
\frac{\partial \Gamma^s_{jk}}
{\partial {x}_l}(\overset{*}{x}) =b_{sjkl}\frac{\partial \Gamma^s_{jk}}{\partial x_l}(x),\end{eqnarray}
where \begin{eqnarray*} b_{sjkl}=\left\{ \begin{array} {ll} 1 \;\;\; &\mbox{if there is no}\;\; n \;\; \mbox{among}\;\; s,j,k,l,\\
-1\;\;\;& \mbox{if there is an}\;\; n\;\;  \mbox{among}\;\; s,j,k,l,\\
1\;\;  \;\; &\mbox{if there are two}\;\; n\;\; \mbox{among}\;\; s,j,k,l,\\
-1 \;\;\;& \mbox{if there are three}\;\; n \;\; \mbox{among}\;\; s,j,k,l,\\
1 \;\;\;& \mbox{if there are  four}\;\; n \;\;\mbox{among}\;\; s,j,k,l,\end{array}\right. \end{eqnarray*}
and \begin{eqnarray}\label{2023.5.22-10} R^j_k(\overset{*}{x}) =c_{jk} R^j_k(x),\end{eqnarray}
where  \begin{eqnarray*}
c_{jk}= \left\{ \begin{array} {ll}1 \;\;\; &\mbox{if there no} \;\; n \;\; \mbox{among}\;\; j,k, \\
- 1 \;\;\; &\mbox{if there is an} \;\; n\;\; \mbox{among}\;\; j,k, \\
1 \;\;\; &\mbox{if there are two} \;\; n  \;\;\mbox{among}\;\; j,k.\end{array}\right.\end{eqnarray*}
  We define \begin{eqnarray} \label{2022.10.18-2} \mathcal{P}=\left\{\begin{array}{ll} \! P_g \;\;\; \;\;\;\,\mbox{on} \;\, \Omega\\
 \!P^\star \;\;\;\; \;\mbox{on} \;\, \Omega^*, \end{array} \right.\end{eqnarray}
where \begin{eqnarray} \label{2022.10.6-8}&& P^\star:=P\Big({g}^{\alpha\beta}(\overset {*}{x}), - {g}^{\alpha n}(\overset{*}{x}),- {g}^{n\beta}(\overset{*}{x}), {g}^{nn}(\overset{*}{x}),a_{jkl}\Gamma^{j}_{kl}(\overset{*}{x}), b_{sjkl}\frac{\partial \Gamma^s_{jk}}
{\partial{x}_l}(\overset{*}{x}),\\
&& \qquad \quad \;\quad c_{jk}R^j_k(\overset{*}{x}), \frac{\partial }{\partial x_1}, \cdots, \frac{\partial }{\partial x_{n-1}},- \frac{\partial }{\partial x_n}\Big),\nonumber \end{eqnarray}
and $\overset{*}{x}=(x',-x_n)\in \Omega^*$.
 Roughly speaking, $P^\star$ is obtained from the expression of $P_g$ by replacing $\frac{\partial}{\partial x_n}$ by $-\frac{\partial}{\partial x_n}$. But we must rewrite such a $P^\star$ in the language of 
the corresponding metric, Christoffel symbols and Ricci curvatures in $\Omega^*$. 
That is, $P^\star$ is got if we replace $g^{\alpha\beta}(x)$, $g^{\alpha n}(x)$, $g^{n\beta} (x)$, $g^{nn}(x)$,
$\{\Gamma^j_{kl}(x)\}_{1\le j,k,l\le n}$, $\{\frac{\partial \Gamma^s_{jk}}{\partial  x_l}(x)\}_{1\le s,j,k,l\le n}$, $\{R^j_k(x)\}_{1\le j,k\le n}$,
 $\frac{\partial }{\partial x_n}$ by ${g}^{\alpha\beta}(\overset{*}{x})$, $-{g}^{\alpha n}(\overset{*}{x})$, $-{g}^{n\beta} (\overset{*}{x})$, ${g}^{nn}(\overset{*}{x})$,
$\;a_{jkl}\Gamma^{j}_{kl}(\overset{*}{x})$, $\,b_{sjkl}\frac{\partial \Gamma^s_{jk}}
{\partial {x}_l}(\overset{*}{x})$, $c_{jk}R^j_k(\overset{*}{x})$, 
 $-\frac{\partial }{\partial x_n}$ in $P_g=P\Big( g^{\alpha\beta} (x)$, $g^{\alpha n}(x)$, $g^{n\beta}(x)$, $g^{nn}(x)$,
$\{\Gamma^j_{kl}(x)\}_{1\le j,k,l\le n}$, $\{\frac{\partial \Gamma^s_{jk}}{\partial  x_l}(x)\}_{1\le s,j,k,l\le n}$, $\{R^j_k(x)\}_{1\le j,k\le n}$,
 $\frac{\partial }{\partial x_1}$, $\cdots$, $\frac{\partial }{\partial x_{n-1}}$, $\frac{\partial }{\partial x_n}\Big)$, respectively. 
Note that we have used the relations (\ref{2021.2.6-3})--(\ref{2023.5.22-10}).
In view of the metric matrices $g$ and $g^*$ have the same order principal minor determinants, we see that $\mathcal{P}$ is still a linear elliptic differential operator on $\mathcal{M}$.

Let $\mathbf{K}(t,x,y)$ be the fundamental solution of the parabolic  system
 \begin{eqnarray*} \left\{ \begin{array}{ll} \frac{\partial \mathbf{u}}{\partial t} + \mathcal{P}\mathbf{u}=0 \;\; &\mbox{in}\;\, (0,+\infty)\times \mathcal{M},\\
  \mathbf{u}=\boldsymbol{\phi} \;\; &\mbox{on}\;\; \{0\}\times \mathcal{M}.\end{array}\right.\end{eqnarray*}
  That is, for any $t\ge 0$ and $x,y\in \mathcal{M}$,
\begin{eqnarray}\label{2020.10.28-10}\left\{\begin{array}{ll}
    \frac{\partial \mathbf{K}(t,x,y)}{\partial t} + \mathcal{P}\mathbf{K}(t,x,y)=0 \;\;  &\mbox{for}\;\, t>0, \, x, y\in \mathcal{M},\\
      \mathbf{K}(0,x,y)=\boldsymbol{\delta}(x-y) \;\; &\mbox{for}\;\;  x,y \in \mathcal{M}.\end{array}\right.\end{eqnarray}
   Here the operator $\mathcal{P}$ is acted in the third argument $y$ of $\mathbf{K}(t,x,y)$.

 Clearly, the coefficients occurring in $\mathcal{P}$ jump as $x$ crosses the $\partial \Omega$ (since the extended metric $g$ is $C^0$-smooth on whole $\mathcal{M}$ and $C^\infty$-smooth in $\mathcal{M}\setminus \partial \Omega$), but $\frac{\partial \mathbf{u}}{\partial t}+\mathcal{P} \mathbf{u}=0$ with $\mathbf{u}(0,x)=\boldsymbol{\phi}(x)$ still has a nice fundamental solution $\mathbf{K}$  of class $[C^1((0,+\infty)\times \mathcal{M}\times  \mathcal{M})]_{n\times n} \cap [ C^\infty((0,+\infty)\times (\mathcal{M} \setminus \partial \Omega) \times (\mathcal{M}\setminus \partial \Omega))]_{n\times n}$, approximable even on $\partial \Omega$  by Levi's sum (see \cite{Liu-21}, or another proof below). Now, let us restrict $x,y\in \Omega$.
 Note that for an (elastic) vector field $\mathbf{u}$ defined in $\Omega$, the boundary traction operator can also  be equivalently written as 
$\mathcal{F} \mathbf{u} := 2\mu  (\mbox{Def}\; \mathbf{u})^\# \, \boldsymbol{\nu} +\lambda\, (\mbox{div}\; \mathbf{u})\,\boldsymbol{\nu}$ on $\partial \Omega$, 
 where  $\mbox{Def} \;\! \mathbf{u}$ is the deformation tensor of the vector field $\mathbf{u}$ (see, \cite{Liu-19}),  $\#$ is the sharp operator by raising index, and $\boldsymbol{\nu}=(\nu_1,\cdots, \nu_n)$ is the unit inner normal to $\partial \Omega$.   
 More precisely,  \begin{eqnarray} \label{23.12.12-1} \;\;\;\quad \; (\mathcal{F} \mathbf{u})^k :=\mu \sum_{l=1}^n\Big( \nu^{\,l} \nabla_l u^k + \nu_l \nabla^k  u^l\Big)+  \lambda\, \boldsymbol{\nu}^k
\sum_{l=1}^n \nabla_l u^l \; \;\mbox{on} \;\, \partial \Omega, \,\;\; k=1,\cdots, n,\;\;\end{eqnarray}  
 where $\nabla_k u^m=u^m_{\;\;\;;\,k}\!:= \frac{\partial u^m}{\partial y_k} +  \sum_{l=1}^n \Gamma^{m}_{kl} u^l$
    is the covariant derivative of the vector field $\mathbf{u}=(u^1, \cdots, u^n)$, and $\nabla^k {u}^m\! :=u^{m;\,k}$ is the raising of  index.     
  Write  $\mathbf{K} (t,x,y)=\big(K^{jk} (t,x,y)\big)_{n\times n}$,  $ \;\mathbf{K}_j(t, x, y)= (K^{j1}(t,x,y)$, $\cdots$, $K^{jn} (t,x,y)$,  $\;(j=1,\cdots, n)$,   and denote   \begin{eqnarray}\label{23.12.13-2} \mathcal{F} \mathbf{K} (t,x,y) =\big( \mathcal{F} \mathbf{K}_1 (t,x,y), \cdots,  \mathcal{F} \mathbf{K}_n(t,x,y)\big)  \end{eqnarray}  
 for all $t>0, x\in \Omega, y\in \partial \Omega$.  More precisely,  for each $j=1,\cdots, n$, the $k$-th component of $\mathcal{F} \mathbf{K}_j$ is 
 \begin{eqnarray} \label{23.12.13-6}  (\mathcal{F} \mathbf{K}_j)^k = \mu \sum_{l=1}^n ( \nu^l \nabla_l K^{jk} +\nu_l \nabla^k K^{jl})
  +\lambda \nu^k \sum_{l=1}^n\nabla_l K^{jl}.\end{eqnarray}
 It is easy to see that $\frac{\partial \big(\mathbf{K}(t,x,y)+\mathbf{K}(t,x,\overset{*}{y})\big)}{\partial \nu_y}\big|_{\partial \Omega}=0$
for  all $t>0$,  $x\in \Omega$ and $y\in \partial \Omega$.  Changing all  terms $\frac{\partial K^{jk}(t,x,y)}{\partial \nu_y}\big|_{\partial \Omega}$ into $0$ in the expression  of $2\mathcal{F}\mathbf{K}(t,x,y)$ (i.e., replacing  all terms $\frac{\partial K^{jk}(t,x,y)}{\partial \nu_y}\big|_{\partial \Omega}= \sum_{m=1}^n \frac{\partial K^{jk}(t,x,y)}{\partial y_m} \, \nu_m$ by $0$ in the expression of   $2\mathcal{F}\mathbf{K}(t,x,y)$),  we obtain  a matrix-valued function $\boldsymbol{\Upsilon}(t,x,y)$  for $t>0$, $x\in \Omega$ and $y\in \partial \Omega$.  This implies that  $\boldsymbol{\Upsilon}(t,x,y)$ only contains the (boundary) tangent derivatives of $K^{jk}(t, x, y)$ with respect to $y\in \partial\Omega$ (without normal derivative of $K^{jk}(t,x,y)$) in the local expression). That is,  in local boundary normal coordinates (the inner normal $\boldsymbol{\nu}$ of $\partial \Omega$ is in the direction of  $x_n$-axis), \begin{eqnarray*}  \boldsymbol{\Upsilon} (t,x,y)= \Big( \boldsymbol{\Upsilon}_1 (t,x,y), \cdots,   \boldsymbol{\Upsilon}_n (t,x,y)\Big), \end{eqnarray*}
   \begin{eqnarray*}\!\!\!\!\!\!\! \!&\!\!\!&\!\!\! \boldsymbol{\Upsilon}_j (t,x,y) \\
\!\!\!\!\!\! \!\!\! &\!\!\!\!\!&\!\!\! \! =\!
  2\mu \!\begin{small}   \begin{pmatrix} \! 2\frac{\partial K^{j1}}{\partial x_1}+2\!\sum\limits_{m=1}^n \!\Gamma^{1}_{1m} K^{jm}  &\! \!\!\cdots\!  \!  & \!\!\!\frac{\partial K^{j1}}{\partial x_{n\!-\!1}} \!\!+\!\!\frac{\partial K^{\!j,n\!-\!1}}{\partial x_1} \!\!+\!\!\! \sum\limits_{m=1}^n\!\! \big( \Gamma^1_{\!n\!-\!1, m}\! \!+\!\!\Gamma^{n\!-\!1}_{\!1m} \big) K^{jm} 
            &   \frac{\partial K^{jn}}{\partial x_1}\!\!+\!\!\!\sum\limits_{m=1}^n\!\! \big( \Gamma^1_{\!nm} \!\!+\!\!\Gamma^n_{1m}\big) K^{jm}  
      \\
\!\!\! \!\cdots \!& \!\cdots \! \!&\!\!\cdots  & \cdots \\
 \!  \! \frac{\partial K^{\!j,n\!-\!1}}{\partial x_{1}}\!\!+\! \!\frac{\partial K^{j1}} {\partial x_{n\!-\!1}}\! \!+\! \!\!\sum\limits_{m=1}^n \!\! \big( \Gamma^{n\!-\!1}_{\!1m}\!\! +\!\! \Gamma^1_{\!n\!-\!1, m} \big)  K^{jm}
 \! \!&\! \cdots    &\!\!\!2 \frac{\partial K^{j,n\!-\!1}}{\partial x_{n\!-\!1}}\!\!+\!\!2\!\sum\limits_{m=1}^n \!\!\Gamma^{n\!-\!1}_{\!n\!-\!1,m} K^{jm} &  \! \frac{\partial K^{jn}}{\partial x_{n\!-\!1}}\!\!+\!\!\sum\limits_{m=1}^n\!\! \big(\Gamma^{n\!-\!1}_{\!\!nm}\! \!+\!\!\Gamma_{\!\!n\!-\!1,m}^n \!\big) K^{\!jm} \!\!
   \\
  \! \frac{\partial K^{jn}}{\partial x_1} \!+\!\!\sum\limits_{m=1}^n \!\!\big( \Gamma_{\!1m}^n\! +\!\Gamma^1_{\!nm}\big) K^{jm}   \! &\!\!\! \cdots  \! \! & 
\frac{\partial K^{jn}}{\partial x_{n\!-\!1}}\! +\!\big(\Gamma^n_{n-1,m}+\Gamma_{\!nm}^{n\!-\!1} \big)K^{jm}   & 
    \;2 \sum\limits_{m=1}^n 
 \Gamma^n_{nm}  K^{jm}\!\!
 \end{pmatrix}\end{small} \!\!\begin{pmatrix} \nu_1 \\ \vdots \\ \nu_n \end{pmatrix}\\  [3mm]
\!\!\!\!\!\!&\!\!\!&\!\!\! \;\;\; + \, 2\lambda\bigg( \frac{\partial K^{j1}}{\partial x_1} + \sum_{m=1}^n \Gamma^{1}_{1m} K^{jm}  + \cdots + \frac{\partial K^{j,n-1}}{\partial x_{n-1}} + \sum_{m=1}^n \Gamma^{n-1}_{n-1,m} K^{jm} +  \sum_{m=1}^n \Gamma^{n}_{nm} K^{jm} \bigg) \begin{pmatrix} \nu_1 \\ \vdots \\ \nu_n \end{pmatrix},\;\;\;  \; j=1,\cdots, n.\end{eqnarray*}
 It is easy to see that  $\boldsymbol{\Upsilon} (t,x,y)$  is a continuous (matrix-valued) function for all $t> 0$, $x\in \Omega$ and $y\in \partial \Omega$. Further, for any  fixed $x\in \Omega$ and any $y\in \partial \Omega$,  since $x\ne y$ we see that 
 \begin{eqnarray} \label{23.12.18-1}  \lim\limits_{t\to 0^+} \mathbf{K} (t,x,y)=0,\end{eqnarray} 
 which implies \begin{eqnarray} \label{23.12.18-2} \lim\limits_{t \to 0^{+}} \frac{\partial \mathbf{K}(t,x,y)}{\partial T} \big|_{\partial \Omega} =0\;\;\;\mbox{for}\;\;\,   x\in \Omega, \;\;\, y\in \partial \Omega,  \end{eqnarray} 
  where $\frac{\partial}{\partial T}$ denotes the tangent derivative along the boundary $\partial \Omega$ in variable $y$. 
 From (\ref{23.12.18-1})--(\ref{23.12.18-2}) we get 
 \begin{eqnarray*} \label{23.12.18-3}  \lim\limits_{t\to 0^+} \nabla_l \mathbf{K} (t,x,y) \big|_{\partial \Omega} =0 \;\; \; \mbox{for all}\;\; x \in \Omega, \; y\in \partial \Omega, \;\;\, 1\le l <n,\end{eqnarray*}
   so that  \begin{eqnarray} \label{23.12.18-4} \lim\limits_{t\to 0^{+}} \boldsymbol{\Upsilon}(t,x,y)=0  \;\;\,\mbox{for any } x\in \Omega, 
   \; y\in \partial \Omega,\end{eqnarray}
    where $\nabla_l \mathbf{K}= (\nabla_l \mathbf{K}_1, \cdots, \nabla_l \mathbf{K}_n)$ and $\nabla_l {K}^{jk}:= \frac{\partial 
       K^{jk} (t,x,y)}{\partial y_l}  +\sum_{m=1}^n\Gamma^k_{lm} K^{jm} (t,x,y)$,  $\,(1\le j,k,l\le  n)$.
     Let $\mathbf{H}(t,x,y)$ be the solution of  
    \begin{eqnarray*} \left\{ \begin{array}{ll}  \frac{\partial \mathbf{u} (t,x,y)}{\partial t} = P_g \mathbf{u} (t,x,y) \; \;\;\mbox{for}\;\,  t>0, \; x,y\in \Omega,\\
2  \mu  \big(\mbox{Def}\, \mathbf{u}(t,x,y)\big)^\#  \,\boldsymbol{\nu}  +\lambda \big(\mbox{div}\, (\mathbf{u}(t,x,y))\big)\, \boldsymbol{\nu} =\boldsymbol{\Upsilon}(t,x,y)\,  \;\;\mbox{for} \,\;  t>0, \, x\in \Omega, \, y\in \partial \Omega,\\
  \mathbf{u}(0, x,y) = \mathbf{0} \,\;\;\mbox{for} \,\;  x, y \in \Omega.\end{array} \right. \end{eqnarray*}
  From (\ref{23.12.18-4}), we get that the above parabolic system satisfy the compatibility condition.   
 Thus,  the  matrix-valued solution $\mathbf{H}(t,x,y)$, is  smooth in $(0, \infty) \times  \Omega \times  \Omega$ and  continous on  $[0, \infty) \times \bar \Omega \times \bar \Omega$.  Then there exists a constant $C>0$ such that 
  \begin{eqnarray*}  | \mathbf{H} (t,x,y) |\le C \,\, \;\;\mbox{for all} \;\;   0\le t\le 1, \; x,\,y\in \bar \Omega, \;\;\,\end{eqnarray*}
   and hence, for dimensions $n\ge 2$,   \begin{eqnarray}\label{23.12.9-1} \int_{\Omega} \mbox{Tr}\; \mathbf{H} (t,x,x)\, dx =nC\,\mbox{vol}(\Omega) =o(t^{-\frac{n-1}{2}})\,  \;\;\mbox{as}\;\; t\to 0^+.\end{eqnarray}
(Actually, we have  $\,\lim_{t \to 0^{+}} \!\int_{\Omega} \mbox{Tr}\; \mathbf{H} (t,x,x)\, dx=0 $). 
   It  can easily be verified  that 
\begin{eqnarray*} 
 && \mathbf{K}^{-} (t,x,y)= \mathbf{K}(t,x,y)- \mathbf{K}(t,x,\overset{*}{y}), \\
&& \mathbf{K}^{+} (t,x,y)= \mathbf{K}(t,x,y)+ \mathbf{K}(t,x,\overset{*}{y})-  \mathbf{H}(t,x,y) ,
\end{eqnarray*}   
 are the Green functions of \begin{eqnarray*}\left\{ \begin{array}{ll}\frac{\partial \mathbf{u}}{\partial t} +{P}_g\mathbf{u}=0\;\;&\mbox{in}\;\, (0,+\infty) \times \Omega,\\
\mathbf{u}=\boldsymbol{\phi} \;\; &\mbox{on}\;\; \{0\}\times \Omega\end{array}\right.\end{eqnarray*} with zero Dirichlet and zero traction (i.e., free) boundary conditions, respectively,
where $y=(y',y_n)$, $y_n\ge 0$, and  $\overset{*}{y}:=\varsigma(y',y_n)=(y^{\prime},-y_n)$.  In other words,
 \begin{eqnarray*}\left\{    \begin{array}{ll} \frac{\partial \mathbf{K}^-(t,x,y)}{\partial t} + {P}_g\mathbf{K}^-(t,x,y)=0,\;\;\; t>0,\, x,\, y\in\Omega,\\
 \mathbf{K}^- (t, x,y)=0, \;\;\;\;  t>0,\; x\in \Omega, \,\;  y\in \partial \Omega,\\
 \mathbf{K}^-(0, x,y)=\boldsymbol{\delta}(x-y), \;\;\;  x,y\in \Omega\end{array} \right. \end{eqnarray*}
and
 \begin{eqnarray*}\left\{  \begin{array}{ll}   \frac{\partial \mathbf{K}^+(t,x,y)}{\partial t} + {P}_g\mathbf{K}^+(t,x,y)=0, \,\;\;\; t>0, \;x, \,y\in \Omega,\\
 \mathcal {F}(\mathbf{K}^+ (t, x,y))=0, \;\; \,\,t>0,\;  x\in \Omega, \;\;  y\in \partial \Omega,\\
 \mathbf{K}^+(0, x,y)=\boldsymbol{\delta}(x-y),\;\;\;  x,y\in \Omega,\end{array} \right. \end{eqnarray*}
 where $ \mathcal {F}(\mathbf{K}^+ (t, x,y))$ is the traction of  $\mathbf{K}^+ (t, x,y) $ on $\partial \Omega$.
   By combining  the fact  $\frac{\partial (\mathbf{K} (t,x,y)+ \mathbf{K}(t,x,\overset{*}{y}))}{\partial \boldsymbol{\nu}_y}\big|_{\partial \Omega} =0$  and $\mathcal{F} \big(\mathbf{K} (t,x,y)+\mathbf{K}(t, x, \overset{*}{y})\big)=\boldsymbol{\Upsilon} (t,x,y) =\mathcal{F} \mathbf{H} (t,x,y)$ for $t>0$,  $x\in \Omega$ and $y\in \partial \Omega$,  we  get $\mathcal{F} \mathbf{K}^+ (t,x,y) =0$ for all $t>0$, $x\in \Omega$ and $y\in \partial \Omega$. 
  Now, we will show  $\mathbf{K}^{-}(t,x,y)$ and $\mathbf{K}^+(t,x,y)$ satisfy the associated parabolic system.  
  In fact, for any $t>0$, $x,y \in \Omega$, we have
  $P_g\mathbf{K}(t,x, y)=\mathcal{P}\mathbf{K}(t,x, y)$,  so that
  \begin{eqnarray}\label{2022.10.29-8}\left\{ \begin{array}{ll}
  \Big( \frac{\partial}{\partial t}+{P}_g\Big) \mathbf{K} (t,x, y)=
  \Big( \frac{\partial}{\partial t}+\mathcal{P}\Big) \mathbf{K} (t,x, y)=0,\\
  \mathbf{K}(0, x,y)=\boldsymbol{\delta}(x-y) \end{array}\right.\end{eqnarray}
  by (\ref{2020.10.28-10}).
Noting that the  Jacobian matrix of the reflection $\varsigma$ is $J$ (see (\ref{2022.11.8-1})), it follows from chain rule that for any fixed $t>0$ and $x\in \Omega$, and any $y=(y',y_n)\in \Omega$,
\begin{align*} &\left[{P}_g ( \mathbf{K}(t,x,\overset{*}{y}))\right]\bigg|_{\text{\normalsize evaluated at the point $y$}}\\
&=\left[{P}_g ( \mathbf{K}(t,x,\varsigma(y',y_n)))\right]\bigg|_{\text{\normalsize evaluated at the point $(y',y_n)$}}
 \\
&= \left[{P}_g (\mathbf{K}(t, x, (y',-y_n))\big)\right]\bigg|_{\text{\normalsize evaluated at the point $(y',y_n)$}}\\
& \! =\!
\left\{\!\left[P \big(g^{\alpha\beta}(y), g^{\alpha n}(y), g^{n\beta} (x), g^{nn}(x),
\big\{\Gamma^j_{kl}(x)\big\}_{1\le j,k,l\le n}, \Big\{\frac{\partial \Gamma^s_{jk}}{\partial  x_l}(x)\Big\}_{1\le s,j,k,l\le n}, \right.\right.\\
&
\left.\left. \quad\quad  \;\big\{R^j_k(x)\big\}_{1\le j,k\le n}, \,\frac{\partial}{\partial y_1}, \cdots, \frac{\partial}{\partial y_{n-1}}, \frac{\partial}{\partial y_n}\big)\right]\! \! \mathbf{K} (t, x, (y'\!,-y_n))\!\right\}\!\Bigg|_{\text{\normalsize evaluated at $\!(y'\!,y_n)$}}\\
\!& \! =
 \left\{\!\left[P \big(g^{\alpha\beta}(\overset{*}{y}), -g^{\alpha n}(\overset{*}{y}), -g^{n\beta}(\overset{*}{y}), g^{nn}(\overset{*}{y}), 
\big\{a_{jkl}\Gamma^{j}_{kl}(\overset{*}{x})\big\}_{1\le j,k, l\le n},\Big\{b_{sjkl}\frac{\partial \Gamma^s_{jk}}
{\partial {x}_l}(\overset{*}{x})\Big\}_{1\le s,j,k,l\le n},\right. \right.\\
& \left.\left.\quad \quad \; \;  \big\{c_{jk}R^j_k(\overset{*}{x})\big\}_{1\le j,k\le n},\,
 \frac{\partial}{\partial y_1}, \cdots, \frac{\partial}{\partial y_{n-1}}, -\frac{\partial}{\partial y_n}\big)\right]\!  \mathbf{K}(t,x,y)\!\right\}\Bigg|_{\text{\normalsize evaluated  at $\overset{*}{y}=(y',-y_n)$}} \\
&= P^\star ( \mathbf{K}(t, x, \overset{*}{y}))\Big|_{\text{\normalsize evaluated at the point $\overset{*}{y}=(y',-y_n)$}}.  \end{align*}
 That is, the action of $P_g$ to $\mathbf{K}(t,x,\overset{*}{y})$ (regarded as a vector-valued function of $y$) at the point ${y}=(y',y_n)$ is just the action of $P^\star$
 to $\mathbf{K}(t,x,\overset{*}{y})$ (regarded as a vector-valued of $\overset{*}{y}$) at the point $\overset{*}{y}=(y',-y_n)$. Because of  $\varsigma(y',y_n)=(y',-y_n)\in \Omega^*$, we see $$P^\star( \mathbf{K}(t,x,\overset{*}{y}))\big|_{\text{\normalsize evaluated at the point $\overset{*}{y}=(y',-y_n)$}}=\mathcal{P}(\mathbf{K} (t,x,\overset{*}{y}))\big|_{\text{\normalsize evaluated at the point $\overset{*}{y}=(y',-y_n)$}}.$$
For any $t>0$, $x\in \Omega$ and $(y',-y_n)\in \Omega^*$, we have $$(\frac{\partial}{\partial t}+ {\mathcal{P}}) (\mathbf{K} (t, x,  (y',-y_n)))=0.$$
In addition, $\mathbf{K} (t, x,  (y',-y_n))= \mathbf{K} (t, x,  \varsigma(y))$ for any $t>0$, $x,y\in \Omega$. By virtue of $x\ne (y',-y_n)$,
 this leads to $\mathbf{K} (0, x,  (y',-y_n))=0$ and
 \begin{eqnarray*} \label{2020.10.29-1} (\frac{\partial}{\partial t}+ {P}_g) \big(\mathbf{K} (t, x,  (y',-y_n))\big)=0\;\; \mbox{for any}\;\,t>0, x\in \Omega \;\,\mbox{and}\;\,(y',-y_n)\in \Omega^*,\end{eqnarray*}
 i.e.,
  \begin{eqnarray} \label{2020.10.29-2}\left\{\! \begin{array}{ll}  (\frac{\partial}{\partial t}+ {{P}_g}) \mathbf{K} (t, x,  \overset{*}{y})=0\;\; \mbox{for any}\;\,t>0, x\in \Omega \;\,\mbox{and}\;\,\overset{*}{y}\in \Omega^*,\\
   \mathbf{K} (0, x,  \overset{*}{y})=0 \;\;\, \mbox{for any}\,\; x,y\in \Omega.\end{array} \right.\end{eqnarray}
 Combining (\ref{2022.10.29-8}) and (\ref{2020.10.29-1}), we obtain that
 \begin{eqnarray} \label{2020.10.29-3} \left\{ \begin{array}{ll} (\frac{\partial}{\partial t}+ {{P}_g}) \Big(\mathbf{K} (t, x, {y})-\mathbf{K} (t, x,  \overset{*}{y})\Big)=0\;\;\mbox{for any}\;\,t>0, \; x,y\in \Omega,\\
 \mathbf{K} (0, x, {y})-\mathbf{K} (0, x,  \overset{*}{y})=\boldsymbol{\delta} (x-y)\;\;\mbox{for any}\;\, x,y\in \Omega.\end{array}\right.\end{eqnarray}
  $\mathbf{K}(t,x, y)$ is $C^1$-smooth with respect to $y$ in $\mathcal{M}$ for any fixed $t>0$ and $x\in \Omega$, so does it on the hypersurface $\partial \Omega$. Therefore, we get that  $\mathbf{K}^-(t,x,y)$ (respectively $\mathbf{K}^+(t,x,y)$) is the Green function in $\Omega$ with the Dirichlet (respectively, free) boundary condition on $\partial \Omega$.

\vskip 0.12 true cm

To show $C^{1}$-regularity of the fundamental solution $\mathbf{K}(t,x,y)$, it suffices to prove  $C^{1,1+\alpha}_{loc}$-regularity for a $W^{1,2}_2$ strong solution  $\mathbf{u}$ of the  parabilic system $\big(\frac{\partial}{\partial t} +\mathcal{P}\big)\mathbf{u}=0$ in $(0,+\infty)\times \mathcal{M}$, where $W^{1,2}_p:= \{ \mathbf{u}\big| \mathbf{u},\frac{\partial \mathbf{u}}{\partial t}, D\mathbf{u}, D^2 \mathbf{u} \in L^p\} $.
 When the coefficients of an elliptic system are smooth on both sides of an $(n-1)$-dimensional hypersurface (may be discontinuous  crossing this hypersurface), the corresponding $C^{1,1+\alpha}$-regularity for solutions of a parabolic equation system is a special case of Dong's result (see Theorem 4 of p.$\,$141 in \cite{Do-12}).
In fact, Dong in \cite{Do-12} has given regularity results to the strong solutions for parabolic equation with  more general coefficients. 
This type of system arises from the problems of linearly elastic laminates and composite materials (see, for example, \cite{CKVC}, \cite{LiVo}, \cite{LiNi}, \cite{Do-12} and \cite{Xi-11}).  

Let us discuss the parabolic (elastic) system in more detail. In fact, from the local expression  (\ref{2023.5.22-12}) of $P_g$, we see that the top-order coefficients of $\mathcal{P}$ are not ``too bad'' since only the first $(n-1)$ coefficients of the $n$-th column in the second matrix in $P_g$ defined on $\Omega$ are changed their signs in $P^\star$ at the reflection points of $\Omega^*$. In other words, only  $(\sum_{m=1}^n g^{1m} (x)
\frac{\partial^2 u^n}{\partial x_m\partial x_n}$,$\cdots$, $\sum_{m=1}^n g^{n-1,m} (x)
\frac{\partial^2 u^n}{\partial x_m\partial x_n})^T$ in $\Omega$ is changed into $(-\sum_{m=1}^n g^{1m} (\overset{*}{x})
\frac{\partial^2 u^n}{\partial x_m\partial x_n}$,$\cdots$, $-\sum_{m=1}^n g^{n-1,m} (\overset{*}{x})
\frac{\partial^2 u^n}{\partial x_m\partial x_n})^T$ in $\Omega^*$ in the second term in $\mathcal{P}$ (see  (\ref{2023.5.22-12}) of $P_g$).
  For any small coordinate chart $V\subset \mathcal{M}$, if  $V\subset \mathcal{M}\setminus (\partial \Omega)$, then the solution $\mathbf{u}$ belongs to $[C^\infty((0, +\infty)\times V)]^n \cap [W^{1,2}_2 ((0, +\infty)\times V)]^n$ since the coefficients of parabolic system 
\begin{eqnarray} \label{2023.6.13-3} (\frac{\partial }{\partial t}+\mathcal{P})\mathbf{u}=0\end{eqnarray} are smooth in $V$.
If the coordinate chart $V\subset \mathcal{M}$ and $V\cap \partial \Omega \ne \varnothing$, then we can find a (local) diffeomorphism $\Psi$ such that $\Psi(V)= U\subset \mathbb{R}^n$ and $\partial \Omega$ is mapped onto $U\cap \{x\in \mathbb{R}^n \big| x_n=0\}$ (i.e., $V\cap \partial \Omega$ is flatten into hyperplane  $\{x\in \mathbb{R}^n\big|x_n=0\}$ by $\Psi$).  By this coordinate transformaton $\Psi$, the parabolic system (\ref{2023.6.13-3}) is changed into another parabolic system $(\frac{\partial}{\partial t}- L)\mathbf{v}=0$, whose coefficients are smooth on both sides of $(n-1)$-dimensional hyperplane  $\{x\in \mathbb{R}^n\big|x_n=0\}$ (may be discontinuous crossing this hyperplane, i.e., the coefficients of $L$ have jump only on this hyperplane). It follows from Dong's regularity result (Theorem 4 of \cite{Do-12}) that $\mathbf{v}\in [C^{1,1+\alpha}((0,+\infty)\times U)]^n$, so that $\mathbf{u}\in  [C^{1,1+\alpha}((0,+\infty)\times V)]^n$. 
 Note that $\mathcal{M}$ is a compact closed Riemannian manifold. Thus we find from the above discussion  and (global) geometric analysis technique that for any initial value $\boldsymbol{\phi}\in [C^\infty (\mathcal{M})]^n$, there exists a (global) strong solution in $[W^{1,2}_2((0,+\infty)\times \mathcal{M})]^n$ which is in  $[C^{1,1+\alpha}((0,+\infty)\times \mathcal{M})]^n$.  Of course, this result  also holds to our case for $C^{1,1+\alpha}$-regularity of the fundamental solution $\mathbf{K} (t,x,y)$ on $(0,+\infty)\times \mathcal{M}\times \mathcal{M}$ (see \cite{Liu-21}).

\vskip 0.20 true cm

Therefore, the integral kernels $\mathbf{K}^{\mp}(t,x,y)$
 of $\frac{\partial \mathbf{u}}{\partial t}+P_g^{\mp}\mathbf{u}=0$ can be expressed on $(0,\infty)\times \Omega\times \Omega$ as
 \begin{eqnarray} \label{c4-23} & {\mathbf{K}}^{-} (t,x,y) =\mathbf{K}(t,x,y)- \mathbf{K}(t,x,\overset{\ast} {y}),\\
 & {\mathbf{K}}^{+} (t,x,y) =\mathbf{K}(t,x,y)+ \mathbf{K}(t,x,\overset{\ast} {y}) -\mathbf{H} (t,x,y)
 \end{eqnarray}
 $\overset{*} {y}$ being the double of $y\in \Omega$ (see, p.$\,$53 of \cite{MS-67}).
 Since the strongly continuous semigroup $(e^{-t\mathcal{P}})_{t\ge 0}$ can also be represented as  \begin{eqnarray*} e^{-t\mathcal{P}} =\frac{1}{2\pi i} \int_{\mathcal{C}} e^{-t\tau} (\tau I- \mathcal{P})^{-1} d\tau,\end{eqnarray*}
where $\mathcal{C}$ is a suitable curve in the complex plane in the positive direction around the spectrum of $\mathcal{P}$ (i.e., a contour around the positive real axis). It follows that
 \begin{eqnarray*}\; \, {\mathbf{K}} (t,x,y) \!=\! e^{-t\mathcal{P}}\delta(x\!-\!y)\! =\! \frac{1}{(2\pi)^n}\!\! \int_{{\Bbb R}^n}\!\! e^{i(x\!-\!y)\cdot\xi} \bigg(\!\frac{1}{2\pi i} \! \int_{\mathcal{C}}\!\! e^{-t\tau}\; \iota \big((\tau {I}\! -\!\mathcal{P})^{\!-1}\big) d\tau\!\bigg) d\xi, \;\; \forall\, t>0,\,  x,y\in \mathcal{M}.\end{eqnarray*}

 We claim that \begin{eqnarray}\; \label{2022.11.11-1}\;\;\, \frac{1}{2\pi i} \int_{\mathcal{C}} (\tau I\!-\! \mathcal{P} )^{-1} e^{-t\tau}\,\delta(x\!-\!y)\,d\tau  = \frac{1}{2\pi i} \int_{\mathcal{C}} \Big( \int_{\mathbb{R}^n} e^{i(x\!-\!y)\cdot \xi} \sum_{j\le -2} \mathbf{q}_j ( x,\xi, \tau) \,d\xi \Big) e^{-t \tau} d\tau.\end{eqnarray}
 In fact, for any smooth vector-valued function $\boldsymbol{\phi}$ with compact support we have \begin{eqnarray*} \big(e^{-t\mathcal{P}} \boldsymbol{\phi} \big)(x) \!\!&\!\!=\!\!&\! \! \Big( \frac{1}{2\pi i} \int_{\mathcal{C}} e^{-t\tau} (\tau I- \mathcal{P} )^{-1}  d\tau \Big) \boldsymbol{\phi}(x)\\
 \!&=\!& \frac{1}{2\pi i} \int_{\mathcal{C}} e^{-t\tau} \Big( \int_{\mathbb{R}^n} e^{i x\cdot \xi} \sum_{j\le -2} \mathbf{q}_j (x, \xi, \tau) \hat{\boldsymbol{\phi}} (\xi) \,d\xi \Big) d\tau.\end{eqnarray*}
On the one hand, from the left-hand side of (\ref{2022.11.11-1}), we get 
\begin{eqnarray} \label{2022.11.11-2}
&& \int\Big[ \Big( \frac{1}{2\pi i} \int_{\mathcal{C}} (\tau I- \mathcal{P} )^{-1}  e^{-t\tau} d\tau\Big) (\delta (x-y))\Big]\boldsymbol{\phi}(y)dy \\
&&   \;\quad\quad \quad =  \Big( \frac{1}{2\pi i} \int_{\mathcal{C}} (\tau- \mathcal{P})^{-1} e^{-t\tau} d\tau\Big) \boldsymbol{\phi}(x)
= e^{-t\mathcal{P} } \boldsymbol{\phi} (x).\nonumber\end{eqnarray}
On the other hand, from the right-hand side of (\ref{2022.11.11-1}) we obtain
\begin{eqnarray} \label{2022.11.11-3} && \int \Big[ \frac{1}{2\pi i} \int_C \Big(\int_{\mathbb{R}^n} e^{i(x-y)\cdot \xi} \sum_{j\le -2} \mathbf{q}_j ( x, \xi,\tau) d\xi\Big) e^{-t\tau} d\tau \Big] \boldsymbol{\phi}(y)dy \\
&& \quad \quad \quad =\frac{1}{2\pi i} \int_C \Big(\int_{\mathbb{R}^n} e^{ix\cdot \xi} \sum_{j\le -2} \mathbf{q}_j ( x, \xi,\tau) d\xi\Big) e^{-t\tau} d\tau \int e^{-y\cdot \xi} \boldsymbol{\phi}(y) dy \nonumber\\
  && \quad \quad \quad = \frac{1}{2\pi i} \int_C \Big(\int_{\mathbb{R}^n} e^{ix\cdot \xi} \sum_{j\le -2} \mathbf{q}_j (x, \xi,\tau) \hat {\boldsymbol{\phi}} (\xi) d\xi\Big) e^{-t\tau} d\tau
= e^{-t\mathcal{P}} \boldsymbol{\phi} (x).\nonumber\end{eqnarray}
 Thus, the desired identity (\ref{2022.11.11-1}) is asserted by (\ref{2022.11.11-2}) and (\ref{2022.11.11-3}).

  In particular, for every $t>0$ and $x\in \Omega$, \begin{eqnarray}\label{2020.7.5-1}&&  {\mathbf{K}} (t,x,x) = e^{-t\mathcal{P}}\delta(x-x) = \frac{1}{(2\pi)^n} \int_{{\Bbb R}^n}\bigg(\frac{1}{2\pi i} \int_{\mathcal{C}} e^{-t\tau}\; \iota \big((\tau {I} -\mathcal{P})^{-1}\big) d\tau\bigg) d\xi\\
    && \qquad\,\; \qquad = \frac{1}{(2\pi)^n} \int_{{\Bbb R}^n}\bigg(\frac{1}{2\pi i} \int_{\mathcal{C}} e^{-t\tau}\; \iota \big((\tau {I} -{P_g})^{-1}\big) d\tau\bigg) d\xi\nonumber\\
    && \qquad\,\;\qquad = \frac{1}{(2\pi)^n} \int_{{\mathbb{R}}^n} \Big( \frac{1}{2\pi i} \int_{\mathcal{C}} e^{-t\tau} \sum_{l\ge 0} q_{-2-l} (x, \xi, \tau) \,d\tau \Big) d\xi,\nonumber\\
\label{2020.7.5-2} &&  {\mathbf{K}} (t,x,\overset{*}{x}) = e^{-t\mathcal{P}}\delta(x-\overset{*}{x}) = \frac{1}{(2\pi)^n} \int_{{\Bbb R}^n} e^{i(x-\overset{*}{x})\cdot\xi} \bigg(\frac{1}{2\pi i} \int_{\mathcal{C}} e^{-t\tau}\; \iota \big((\tau {I} -\mathcal{P})^{-1}\big) d\tau\bigg) d\xi\\
&& \qquad \qquad\;\;\; = \frac{1}{(2\pi)^n} \int_{{\Bbb R}^n} e^{i(x-\overset{*}{x})\cdot\xi} \bigg(\frac{1}{2\pi i} \int_{\mathcal{C}} e^{-t\tau}\; \iota \big((\tau {I} -{\mathcal{P}})^{-1}\big) d\tau\bigg) d\xi\nonumber\\
 && \qquad\;\;\;\qquad = \frac{1}{(2\pi)^n} \int_{{\mathbb{R}}^n} e^{i(x-\overset{*}{x})\cdot \xi} \Big( \frac{1}{2\pi i} \int_{\mathcal{C}} e^{-t\tau} \sum_{l\ge 0} q_{-2-l} (x, \xi, \tau) \,d\tau \Big) d\xi,\nonumber
 \end{eqnarray}
 where $\sum_{l\ge 0} {\mathbf{q}}_{-2-l} (x,\xi,\tau) $ is the full symbol of $(\tau I-P_g)^{-1}$.

 Firstly, from the discussion on p.$\,$10182 of \cite{Liu-21}, we know that
   \begin{eqnarray}\label{2022.11.9-1}
   && \;\;\;\;{\mathbf{q}}_{-2} (x,\xi,\tau) =\frac{1}{\tau- \mu \sum\limits_{l,m=1}^n g^{lm}\xi_l \xi_m }\,{\mathbf{I}}_n\\
 && \qquad \;\;\;\; \; \; +\frac{\mu+\lambda}{  \big(\tau- \mu \sum\limits_{l,m=1}^n g^{lm}\xi_l \xi_m \big)\big(\tau- (2\mu+\lambda)  \sum\limits_{l,m=1}^n g^{lm}\xi_l \xi_m \big)}\begin{bmatrix} \sum\limits_{r=1}^n g^{1r} \xi_r \xi_1 &\cdots &  \sum\limits_{r=1}^n g^{1r} \xi_r \xi_n\\
\vdots& {} &\vdots \\
 \sum\limits_{r=1}^n g^{nr} \xi_r \xi_1 &\cdots &  \sum\limits_{r=1}^n g^{nr} \xi_r \xi_n\end{bmatrix}\nonumber\end{eqnarray}
 and   \begin{eqnarray}\label{2020.6.6-1}  && \mbox{Tr}\, \big({\mathbf{q}}_{-2} (x,\xi,\tau)\big)=
   \frac{n}{\big(\tau \!-\!\mu \!\sum_{l,m=1}^n \!g^{lm} \xi_l\xi_m\big)} \\
 && \qquad \qquad \qquad \qquad \;\;  + \frac{(\mu+\lambda)\sum_{l,m=1}^n g^{lm} \xi_l\xi_m}{
  \big(\!\tau\! -\!\mu\! \sum_{l,m=1}^n \!g^{lm} \xi_l\xi_m\!\big)
 \big(\!\tau \!-\!(2\mu\!+\!\lambda)  \!\sum_{l,m=1}^n \!g^{lm} \xi_l\xi_m\!\big)}
.\nonumber \end{eqnarray}
       For each $x\in \Omega$, we use a geodesic normal coordinate system centered at this $x$. It follows from \S11 of Chap.1 in \cite{Ta-1}
        that in such a coordinate system, $g_{jk}(x)=\delta_{jk}$ and $\Gamma_{jk}^l(x)=0$. Then (\ref{2020.6.6-1}) reduces to
        \begin{eqnarray}\label{2020.7.6-3} \quad \quad \,\;\mbox{Tr} \,\big({\mathbf{q}}_{-2} (x, \xi,\tau)\big)=
     \frac{n}{(\tau -\mu |\xi|^2)}+
      \frac{(\mu+\lambda)|\xi|^2}{(\tau -\mu |\xi|^2) (\tau-(2\mu+\lambda)|\xi|^2)},\end{eqnarray}
  where $|\xi|=\sqrt{\sum_{k=1}^n \xi^2_k}$ for any $\xi\in {\mathbb{R}}^n$.
   By applying the residue theorem (see, for example, Chap.$\,$4, \S5 in \cite{Ahl}) we get
 \begin{eqnarray} \label{3.10} && \frac{1}{2\pi i}\! \int_{\mathcal{C}}\! e^{-t\tau} \bigg(\frac{n}{(\tau -\mu |\xi|^2)}\!+\!
      \frac{(\mu+\lambda)|\xi|^2}{(\tau -\mu |\xi|^2) (\tau-(2\mu\!+\!\lambda)|\xi|^2)} \bigg) d\tau\!=\! (n\!-\!1) e^{-t\mu|\xi|^2} \!+\!  e^{-t(2\mu+\lambda)|\xi|^2}.\end{eqnarray}
       It follows that
      \begin{eqnarray}\label{2020.7.10-1} \frac{1}{(2\pi)^n}\!\!\!\!&\!\!\!&\!\!\!\!\! \!\!\! \!\! \int_{{\mathbb{R}}^n}\Big( \frac{1}{2\pi i} \int_{\mathcal{C}} e^{-t\tau}\, \mbox{Tr}\,({\mathbf{q}}_{-2} (x, \xi,\tau) ) d\tau \Big) d\xi \\
      \!\!\! &=\!\!\!&
    \frac{1}{(2\pi)^n} \int_{{\Bbb R}^n} \bigg((n-1) e^{-t\mu|\xi|^2} +  e^{-t(2\mu+\lambda)|\xi|^2}\bigg)
         d\xi \nonumber\\
               \!\!\! &=\!\!\!&\frac{n-1}{(4\pi \mu t)^{n/2}} +  \frac{1}{(4\pi (2\mu+\lambda) t)^{n/2}},
                  \nonumber\end{eqnarray}
                  and hence
        \begin{eqnarray}\label{2020.7.12-3}  && \int_{\Omega}\! \left\{\! \frac{1}{(2\pi)^n}\! \int_{{\mathbb{R}}^n}\!\Big( \frac{1}{2\pi i} \int_{\mathcal{C}} e^{-t\tau}\, \mbox{Tr}\,({\mathbf{q}}_{-2} (x, \xi,\tau) ) d\tau\! \Big) d\xi\!\right\}\! dV\\
        && \qquad \;\;=\Big(\frac{n-1}{(4\pi \mu t)^{n/2}} \! +\!  \frac{1}{(4\pi (2\mu\!+\!\lambda) t)^{n/2}}\!\Big){\mbox{Vol}}(\Omega).\nonumber\end{eqnarray}

In the above discussion, if we replace $x\in \Omega$ by $\overset{*}{x} \in \Omega^*$, then (\ref{2022.11.9-1}) will become
\begin{eqnarray*}\label{2022.11.9-2}
  \!\! \!\!\!&& \!\!\!\!\!{\mathbf{q}}_{-2} (\overset{*}{x},\xi,\tau) =\frac{1}{\tau- \mu \sum\limits_{l,m=1}^n (g^{lm}(\overset{*}{x}))\xi_l \xi_m }\,{\mathbf{I}}_n +\frac{\mu\!+\!\lambda}{  \big(\tau\!-\! \mu \sum\limits_{l,m=1}^n\! (g^{lm}(\overset{*}{x}))\xi_l \xi_m \big)\big(\tau\!- \!(2\mu\!+\!\lambda)  \sum\limits_{l,m=1}^n \!(g^{lm}(\overset{*}{x}))\xi_l \xi_m \big)}\\
 && \;\;\;\;\times \begin{bmatrix}\! \sum\limits_{r=1}^n \!(g^{1r}(\overset{*}{x})) \xi_r \xi_1 \!&\! \cdots \!&\! \sum\limits_{r=1}^n( g^{1r}(\overset{*}{x}))\xi_r \xi_{n-1}\!& \!\sum\limits_{r=1}^n \!(-g^{1r}(\overset{*}{x})) \xi_r \xi_n\\
\vdots\!& {}\! &\vdots &\vdots \\
\sum\limits_{r=1}^n (g^{n\!-\!1,r}(\overset{*}{x})) \xi_r \xi_1 \!\!&\!\cdots \!&\! \sum\limits_{r=1}^n (g^{n-1,r}(\overset{*}{x})) \xi_r \xi_{n\!-\!1}\!& \! \sum\limits_{r=1}^n (-g^{n\!-\!1,r}(\overset{*}{x})) \xi_r \xi_n
\\
 \sum\limits_{r=1}^n\! (-g^{nr}(\overset{*}{x})) \xi_r \xi_1 \!&\!\cdots \!& \!\sum\limits_{r=1}^n \!(-g^{nr}(\overset{*}{x})) \xi_r \xi_{n-1}\!& \! \sum\limits_{r=1}^n (g^{nr} (\overset{*}{x}))\xi_r \xi_n\end{bmatrix}\nonumber\end{eqnarray*}
 and   \begin{eqnarray*}\label{2022.11.9-4}  \mbox{Tr}\, \big({\mathbf{q}}_{-2} (\overset{*}{x},\xi,\tau)\big)=
   \frac{n}{\big(\tau \!-\!\mu \!\sum_{l,m=1}^n \!(g^{lm}(\overset{*}{x})) \xi_l\xi_m\big)} \!+\! \frac{(\mu+\lambda)\sum_{l,m=1}^n (g^{lm}(\overset{*}{x})) \xi_l\xi_m}{
  \big(\!\tau\! -\!\mu\! \sum_{l,m=1}^n \!(g^{lm}(\overset{*}{x})) \xi_l\xi_m\!\big)
 \big(\!\tau \!-\!(2\mu\!+\!\lambda)  \!\sum_{l,m=1}^n \!(g^{lm}(\overset{*}{x})) \xi_l\xi_m\!\big)}
.\end{eqnarray*}
This implies that all expressions (\ref{2020.6.6-1})--(\ref{2020.7.12-3}) of the above  trace symbols  have the same form either in $\Omega$ or in $\Omega^*$.

      For given (small) $\epsilon>0$ , denote by $U_\epsilon(\partial \Omega)=\{z\in {\mathcal{M}}\big| \mbox{dist}\, (z, \partial \Omega)<\epsilon\}$ the $\epsilon$-neighborhood of $\partial \Omega$ in $\mathcal{M}$.
             When $x\in \Omega\setminus U_\epsilon (\partial \Omega)$,
we see by taking geodesic normal coordinate system at $x$ that (\ref{2020.7.6-3}) still holds at this $x$. According to (\ref{3.10}) we have
  that
     \begin{eqnarray} \label{3.11}  \mbox{Tr}\,({\mathbf{q}}_{-2}(t, x, \overset{*}{x}))\!\!\! &=\!\!\!&
     \frac{1}{(2\pi)^n} \int_{{\Bbb R}^n} e^{i(x-\overset{*}{x})\cdot\xi}
    \bigg((n-1) e^{-t\mu|\xi|^2} +  e^{-t(2\mu+\lambda)|\xi|^2}\bigg)
         d\xi \nonumber\\
               \!\!\! &=\!\!\!&\frac{n-1}{(4\pi \mu t)^{n/2}}e^{-\frac{|x-\overset{*}{x}|^2}{4t\mu}} +  \frac{1}{(4\pi (2\mu+\lambda) t)^{n/2}}e^{-\frac{|x-\overset{*}{x}|^2}{4t(2\mu+\lambda)}} \quad \, \mbox{for any}\;\, x\in \Omega \setminus U_\epsilon (\partial \Omega),
                  \nonumber\end{eqnarray}
     which exponentially tends to zero as $t\to 0^+$ because $|x-\overset{*}{x}|\ge \epsilon$.  Hence
     \begin{eqnarray} \label{2020.7.6-4}  \int_{\Omega \setminus U_\epsilon (\partial \Omega)} \left(\mbox{Tr}\,({\mathbf{q}}_{-2}(t, x, \overset{*}{x}))\right)\,dV\!\!\! &=\!\!\!&
     O(t^{1-\frac{n}{2}}) \quad \;\mbox{as}\;\; t\to 0^+.
                 \end{eqnarray}

                Secondly, for $l\ge 1$, it can be verified that
                  $\mbox{Tr}\, ({\mathbf{q}}_{-2-l} (x, \xi, \tau))$ is a sum of finitely many terms, each of which  has the following form:
                  $$\frac{r_k(x, \xi)}{(\tau-\mu \sum_{l,m=1}^n g^{lm} \xi_l \xi_m )^s (\tau -(2\mu +\lambda)\sum_{l,m=1}^n g^{lm} \xi_l\xi_m )^j },$$
where $k-2s-2j=-2-l$, and $r_k(x,\xi)$ is the symbol independent of $\tau$ and homogeneous of degree $k$.
 Again we take the geodesic normal coordinate systems center at $x$ (i.e., $g_{jk}(x)=\delta_{jk}$ and $\Gamma_{jk}^l(x)=0$), by applying residue theorem we see  that, for $l\ge 1$ , \begin{eqnarray*}\label{2020.7.12-1} && \frac{1}{(2\pi)^n} \int_{{\mathbb{R}}^n}\Big( \frac{1}{2\pi i} \int_{\mathcal{C}} e^{-t\tau}\, \mbox{Tr}\,({\mathbf{q}}_{-2-l} (x, \xi,\tau) ) d\tau \Big) d\xi=O(t^{l-\frac{n}{2}}) \;\;\mbox{as}\;\, t\to 0^+ \;\;\, \mbox{uniformly for} \;\, x\in \Omega,
 \end{eqnarray*} and
\begin{eqnarray} \label{2020.7.14-1} \\
 \frac{1}{(2\pi)^n} \!\int_{{\mathbb{R}}^n}\! e^{i(x-\overset{*}{x})\cdot \xi} \Big(\! \frac{1}{2\pi i} \!\int_{\mathcal{C}}\! e^{-t\tau}\, \mbox{Tr}\,({\mathbf{q}}_{-2-l} (x, \xi,\tau) ) d\tau \!\Big) d\xi\!=\!O(t^{l-\frac{n}{2}}\!) \;\,\mbox{as}\;\, t\to 0^+ \,\; \mbox{uniformly for} \;\, x\in \Omega.\nonumber
      \end{eqnarray}
Therefore
\begin{eqnarray}\label{2020.7.12-10} && \int_{\Omega} \bigg\{\frac{1}{(2\pi)^n} \int_{{\mathbb{R}}^n}\Big( \frac{1}{2\pi i} \int_{\mathcal{C}} e^{-t\tau}\, \sum_{l\ge 1}\mbox{Tr}\,({\mathbf{q}}_{-2-l} (x, \xi,\tau) ) d\tau \Big) d\xi\bigg\} dV =O(t^{1-\frac{n}{2}}) \;\;\mbox{as}\;\, t\to 0^+,
 \end{eqnarray} and
\begin{eqnarray}\label{2020.7.13-5}\;\;\;\;\;\;\;\quad\,\;
\int_{\Omega} \!\bigg\{\! \frac{1}{(2\pi)^n} \!\int_{{\mathbb{R}}^n}\!\!\! e^{i(x-\overset{*}{x})\cdot \xi} \Big( \frac{1}{2\pi i} \!\int_{\mathcal{C}} e^{-t\tau} \sum_{l\ge 1}\mbox{Tr}\,({\mathbf{q}}_{-2-l} (x, \xi,\tau) ) d\tau \!\Big) d\xi\bigg\}dV\!=\!O(t^{1-\frac{n}{2}})\;\;\mbox{as}\,\, t\to 0^+.
      \end{eqnarray}
Combining (\ref{2020.7.5-1}), (\ref{2020.7.12-3}) and (\ref{2020.7.12-10}), we have
       \begin{eqnarray} \label{3.11}&& \int_{\Omega}\mbox{Tr}\,({\mathbf{K}}(t, x, x)) \, dV =    \bigg[  \frac{n-1}{(4\pi \mu t)^{n/2}} +  \frac{1}{(4\pi (2\mu+\lambda) t)^{n/2}}\bigg]{\mbox{Vol}}(\Omega)+O(t^{1-\frac{n}{2}})\;\;\mbox{as}\;\; t\to 0^+. \end{eqnarray}

     Finally, we will consider the case of $\int_{\Omega\cap U_\epsilon(\partial \Omega)} \left\{\frac{1}{(2\pi)^n} \int_{{\mathbb{R}}^n} e^{i(x-\overset{*}{x})\cdot \xi} \Big(  \frac{1}{2\pi i} \int_{\mathcal{C}} e^{-t\tau}\,\mbox{Tr}\,({\mathbf{q}}_{-2} (x, \xi,\tau) ) d\tau \Big) d\xi\right\} dV$.
We pick a self-double patch $W$ of $\mathcal{M}$ (such that $W\subset U_\epsilon (\partial \Omega)$) covering a patch $W\cap \partial \Omega$ of $\partial \Omega$  endowed (see the diagram on p.$\,$54 of \cite{MS-67}) with local coordinates $x$ such that
\begin{figure}[h]
\centering
%\caption{\label{figure1}}
\begin{tikzpicture}[scale=1,line width=0.8]
\clip (-1,1.3) rectangle (6.5,7.4);

\path(70:7) coordinate(A);
\path(30:7) coordinate(B);
\path(70:4.5) coordinate(C);
\path(30:4.5) coordinate(D);

\draw (A) arc (70:30:7);
\draw (C) arc (70:30:4.5);
\draw (A)--(C);
\draw (B)--(D);

\draw[fill=black] (3.177,3.6) circle (1pt);
\draw[fill=black] (4.3,5) circle (1pt);
\draw[fill=black] (3.67,4.32) circle (1pt);

\draw[->,>=stealth] (4.3,5) .. controls (3.5,4.3) and (2.8,3) .. (2.6,2.5);

\node at (68:7.5) {$\Omega^{*}$};
\node at (2.2,2.65) {$\Omega$};
\node at (2.6,2) {$x_n>0$};
\node at (3.35,3.35) {$x$};
\node at (4.5,4.7) {$x^*$};
\node at (-0.5,6.7) {$\partial \Omega$};
\node at (-2,-2){};

\draw (-0.524,6.3709) arc (80.0001:25:8.5);
\end{tikzpicture}
\end{figure}
   $\epsilon>x_n>0$ in $W\cap \Omega$; $\,x_n=0$ on $W\cap \partial \Omega$;
  $\; x_n (\overset{*}{x})=-x_n(x)$; and the positive $x_n$-direction is perpendicular to $\partial \Omega$. This has the effect that (\ref{2021.2.6-3})--(\ref{2021.2.6-4})  and  \begin{eqnarray}
   \label{2021.2.6-5}  \sqrt{|g|/g_{nn}} \; dx_1\cdots dx_{n-1}\!\!\!&\!=\!&\!\!\! \mbox{the element of (Riemannian) surface area on} \,\, \partial \Omega.\end{eqnarray}
    We choose coordinates $x'=(x_1,\cdots, x_{n-1})$ on an open set in $\partial \Omega$ and then
 coordinates $(x', x_{n})$ on a neighborhood in $\bar\Omega$ such that
$x_{n}=0$ on $\partial \Omega$ and $|\nabla x_{n}|=1$ near $\partial \Omega$ while $x_{n}>0$ on $\Omega$ and such that $x'$ is constant
on each geodesic segment in $\bar\Omega$ normal to $\partial \Omega$.
   Then the metric tensor on $\bar \Omega$ has
the form (see \cite{LU} or p.$\,$532 of \cite{Ta-2})
\begin{eqnarray} \label{a-1}  \big(g_{jk} (x',x_{n}) \big)_{n\times n} =\begin{pmatrix} ( g_{jk} (x',x_{n}))_{(n-1)\times (n-1)}& 0\\
      0& 1 \end{pmatrix}. \end{eqnarray}
       Furthermore, we can take a geodesic normal coordinate system for $(\partial \Omega, g)$ centered at $x_0=0$, with respect to $e_1, \cdots, e_{n-1}$, where  $e_1, \cdots, e_{n-1}$ are the principal curvature vectors. As Riemann showed, one has (see p.$\,$555 of \cite{Ta-2})
            \begin{eqnarray} \label{7/14/1}& g_{jk}(x_0)= \delta_{jk}, \; \; \frac{\partial g_{jk}}{\partial x_l}(x_0)
 =0  \;\;  \mbox{for all} \;\; 1\le j,k,l \le n-1,\\
 & -\frac{1}{2}\frac{\partial g_{jk}}{\partial x_n} (x_0) =\kappa_k\delta_{jk}  \;\;  \mbox{for all} \;\; 1\le j,k \le n-1,\nonumber
 \end{eqnarray}
 where $\kappa_1\cdots, \kappa_{n-1}$ are the principal curvatures of $\partial \Omega$ at point $x_0=0$.
   Due to the special geometric normal coordinate system and (\ref{7/14/1})--(\ref{a-1}), we see that for any $x\in \{z\in \Omega\big|\mbox{dist}(z, \partial \Omega)< \epsilon\}$,
   \begin{eqnarray} \label{18-4-5-1} x-\overset{\ast} {x}=(0,\cdots, 0, x_n-(-x_n))=(0,\cdots, 0,2x_n).\end{eqnarray}
              By (3.17) of \cite{Liu-21}, (\ref{7/14/1}), (\ref{3.10}) and (\ref{18-4-5-1}), we find that
      \begin{eqnarray*} \label{18-4-1-1} &&\int_{W\cap \Omega} \bigg\{\frac{1}{(2\pi)^n} \int_{{\Bbb R}^{n}} e^{i\langle x-\overset{*}{x}, \xi\rangle} \Big( \frac{1}{2\pi i} \int_{\mathcal{C}} e^{-t\tau} \,\mbox{Tr}\,\big({\mathbf{q}}_{-2} (x, \xi, \tau)\big) d\tau \Big)  \, d\xi\bigg\} dV \\
           &&   =\! \int_0^\epsilon dx_n\! \int_{W\cap \partial \Omega}\!
 \frac{dx'}{(2\pi)^n}\!\int_{{\Bbb R}^{n}} \!e^{i\langle 0, \xi'\rangle + i2x_n \xi_n} \bigg[\!\frac{1}{2\pi i}\! \int_{\mathcal{C}}\! e^{-t\tau} \!\Big(  \frac{n}{(\tau \!-\!\mu |\xi|^2)}\!+\!
      \frac{(\mu\!+\!\lambda)|\xi|^2}{(\tau \!-\!\mu |\xi|^2) (\tau\!-\!(2\mu\!+\!\lambda)|\xi|^2)}\Big) d\tau\!\bigg] d\xi\\
  && \,=  \int_0^\epsilon dx_n \int_{W\cap \partial \Omega}
 \frac{dx'}{(2\pi)^n}\int_{{\Bbb R}^{n}} e^{i2x_n \xi_n} \bigg((n-1) e^{-t\mu|\xi|^2} +  e^{-t(2\mu+\lambda)|\xi|^2} d\tau\bigg) d\xi\\
 &&\,= \int_0^\epsilon dx_n \int_{W\cap \partial \Omega}  \frac{dx'}{(2\pi)^n}\int_{-\infty}^\infty e^{2ix_n \xi_n}\bigg[ \int_{{\Bbb R}^{n-1}} \bigg( (n-1) e^{-t\mu(|\xi'|^2+\xi_n^2)} +  e^{-t(2\mu+\lambda)(|\xi'|^2+\xi_n^2)} \bigg) d\xi' \bigg]d\xi_n  \nonumber\\
    &&\,= \int_0^\epsilon dx_n \int_{W\cap \partial \Omega}  \frac{1}{(2\pi)^n}\bigg[\int_{-\infty}^\infty e^{2ix_n \xi_n}e^{-t\mu \xi_n^2} \bigg( \int_{{\Bbb R}^{n-1}}  (n-1) e^{-t\mu \sum_{j=1}^{n-1}\xi_j^2}  d\xi' \bigg)d\xi_n\bigg] dx'\nonumber\\
  &&\quad\, + \int_0^\epsilon dx_n \int_{W\cap \partial \Omega}  \frac{1}{(2\pi)^n}\bigg[\int_{-\infty}^\infty e^{2ix_n \xi_n}e^{-t(2\mu+\lambda) \xi_n^2} \bigg( \int_{{\Bbb R}^{n-1}} e^{-t(2\mu+\lambda) \sum_{j=1}^{n-1}\xi_j^2} d\xi' \bigg)d\xi_n\bigg]dx'\nonumber, \end{eqnarray*}
 where  $\xi=(\xi', \xi_n)\in {\Bbb R}^n$, $\xi'=(\xi_1, \cdots, \xi_{n-1})$.
   A direct calculation shows that
   \begin{eqnarray*}
 &&  \frac{1}{(2\pi)^n}\bigg[\int_{-\infty}^\infty e^{2ix_n \xi_n}e^{-t\mu \xi_n^2} \bigg( \int_{{\Bbb R}^{n-1}}  (n-1) e^{-t\mu \sum_{j=1}^{n-1}\xi_j^2}  d\xi' \bigg)d\xi_n\bigg]=\frac{n-1}{(4\pi \mu  t)^{n/2}}\, e^{-\frac{(2x_n)^2}{ 4\mu t} },\\
  &&    \frac{1}{(2\pi)^n}\bigg[\int_{-\infty}^\infty e^{2ix_n \xi_n}e^{-t(2\mu+\lambda) \xi_n^2} \bigg( \int_{{\Bbb R}^{n-1}} e^{-t(2\mu+\lambda) \sum_{j=1}^{n-1}\xi_j^2} d\xi' \bigg)d\xi_n= \frac{1}{(4\pi (2 \mu+\lambda)  t)^{n/2}}\, e^{-\frac{(2x_n)^2}{ 4(2\mu+\lambda) t} }.\end{eqnarray*}
   Hence
    \begin{eqnarray} \label{18-4-1-10.} \\
 && \int_{W\cap \Omega} \bigg\{\frac{1}{(2\pi)^n} \int_{{\Bbb R}^{n}} e^{i\langle x-\overset{*}{x}, \xi\rangle} \Big( \frac{1}{2\pi i} \int_{\mathcal{C}} e^{-t\tau} \,\mbox{Tr}\,\big({\mathbf{q}}_{-2} (x, \xi, \tau)\big) d\tau \Big)  \, d\xi\bigg\} dV\nonumber \\
  &&\;\;\quad = \int_0^\epsilon dx_n \int_{W\cap \partial \Omega} \left[\frac{n-1}{(4\pi \mu  t)^{n/2}}\, e^{-\frac{(2x_n)^2}{ 4\mu t} }+\frac{1}{(4\pi (2 \mu+\lambda)  t)^{n/2}}\, e^{-\frac{(2x_n)^2}{ 4(2\mu+\lambda) t} }\right]dx'
    \nonumber\\
    &&\;\; \quad= \int_0^\infty dx_n \int_{W\cap \partial \Omega} \left[ \frac{n-1}{(4\pi \mu  t)^{n/2}}\, e^{-\frac{(2x_n)^2}{ 4\mu t} }+\frac{1}{(4\pi (2\mu+\lambda)  t)^{n/2}}\, e^{-\frac{(2x_n)^2}{ 4(2\mu+\lambda) t} }\right]dx'\nonumber\\
    && \;\; \quad\;\quad - \int_\epsilon^\infty dx_n \int_{W\cap \partial \Omega} \left[\frac{n-1}{(4\pi \mu  t)^{n/2}}\, e^{-\frac{(2x_n)^2}{ 4\mu t} }+ \frac{1}{(4\pi (2\mu +\lambda) t)^{n/2}}\, e^{-\frac{(2x_n)^2}{ 4(2\mu+\lambda) t} }\right]dx'\nonumber\\
   &&\;\;\quad= \frac{n-1}{4} \cdot \frac{\mbox{Vol}(W\cap \partial \Omega)}{(4\pi \mu t)^{(n-1)/2}} + \frac{1}{4} \cdot \frac{\mbox{Vol}(W\cap \partial \Omega)}{(4\pi (2\mu +\lambda) t)^{(n-1)/2}}\nonumber\\
   &&\;\;\quad \quad \, -   \int_{W\cap \partial \Omega} \left\{\int_\epsilon^\infty\bigg[\frac{n-1}{(4\pi \mu  t)^{n/2}}\, e^{-\frac{(2x_n)^2}{ 4\mu t} }+ \frac{1}{(4\pi (2\mu +\lambda) t)^{n/2}}\, e^{-\frac{(2x_n)^2}{ 4(2\mu+\lambda) t}} \bigg]dx_n\right\}dx'. \nonumber\end{eqnarray}
      It is easy to verify that for any  fixed $\epsilon>0$,
   \begin{eqnarray} \label{18-4-1-11.}\begin{aligned} && \int_\epsilon^\infty \frac{1}{(4\pi \lambda t)^{\frac{n}{2}}} e^{-\frac{(2x_n)^2}{4\mu t}}dx_n = O(t^{1-n/2})\quad \; \, \mbox{as} \;\, t\to 0^+,\qquad \;\; \quad \\
  && \int_\epsilon^\infty \frac{1}{(4\pi (2\mu+\lambda) t)^{\frac{n}{2}}} e^{-\frac{(2x_n)^2}{4(2\mu+\lambda) t}}dx_n = O(t^{1-n/2})
   \quad \; \, \mbox{as} \;\, t\to 0^+.\end{aligned}\end{eqnarray}
   From (\ref{18-4-1-10.}) and (\ref{18-4-1-11.}), we get that
     \begin{eqnarray}\label{18-4-1-6}      &&\;\, \int_{W\cap \Omega} \bigg\{\frac{1}{(2\pi)^n} \int_{{\Bbb R}^{n}} e^{i\langle x-\overset{*}{x}, \xi\rangle} \Big( \frac{1}{2\pi i} \int_{\mathcal{C}} e^{-t\tau} \,\mbox{Tr}\,\big({\mathbf{q}}_{-2} (x, \xi, \tau)\big) d\tau \Big)  \, d\xi\bigg\} dV=
     \frac{n\!-\!1}{4} \cdot \frac{\mbox{Vol}(W\cap \partial \Omega)}{(4\pi \mu t)^{(n-1)/2}}\\
  && \; \quad\,\quad \;\,\quad  + \frac{1}{4} \cdot \frac{\mbox{Vol}(W\cap \partial \Omega)}{(4\pi (2\mu +\lambda) t)^{(n-1)/2}}
 + O(t^{1-n/2}) \quad \; \mbox{as} \, \; t\to 0^+.\nonumber\end{eqnarray}
  For any $x\in \Omega\cap U_\epsilon (\partial \Omega)$, we have
  \begin{eqnarray} \label{18-4-3-1}\;\;\,\;\, \;\;\;\;\; \mbox{Tr}\,(K(t, x, \overset{*}{x}))
 \!\!\!&\!=&\!\!\!\! \frac{1}{(2\pi)^n}\int_{{\Bbb R}^n} \!  e^{i\langle x-\overset{*}{x}, \xi\rangle}\big( \frac{1}{2\pi i} \!\int_{\mathcal{C}}\! e^{-t\tau} \mbox{Tr}\,\big({\mathbf{q}}_{-2} (x, \xi, \tau)\big) d\tau\big) d\xi\\
  \!\!\!&\!&\!\!\!\!+ \frac{1}{(2\pi)^n}\int_{{\Bbb R}^n}   e^{i\langle x-\overset{*}{x}, \xi\rangle}\big(\sum_{l\ge 1} \frac{1}{2\pi i} \int_{\mathcal{C}} e^{-t\tau} \mbox{Tr}\,\big({\mathbf{q}}_{-2-l} (x, \xi, \tau)\big) d\tau\big) d\xi\nonumber\\
  \!\!\!&\!=&\!\!\!\! \!\frac{1}{(2\pi)^n}\!\!\int_{{\Bbb R}^n} \! \! e^{i\langle x-\overset{*}{x}, \xi\rangle}\big( \frac{1}{2\pi i} \!\int_{\mathcal{C}} e^{-t\tau} \mbox{Tr}\,\big({\mathbf{q}}_{-2} (x, \xi, \tau)\big) d\tau\big) d\xi\!+\!O(t^{1+\frac{n}{2}}) \;\; \mbox{as}\;\, t\to 0^+,\nonumber\end{eqnarray}  where the second equality used (\ref{2020.7.14-1}).
  Combining (\ref{18-4-1-6}) and (\ref{18-4-3-1}), we have
     \begin{eqnarray} \label{3..15}
     &&\int_{W\cap \Omega}\mbox{Tr}\big( \mathbf{K}(t, x, \overset{*}{x})\big)dx=  \frac{n-1}{4} \cdot \frac{\mbox{Vol}(W\cap \partial \Omega)}{(4\pi \mu t)^{(n-1)/2}} \\
  && \quad\,\quad \;\,\quad \;+ \frac{1}{4} \cdot \frac{\mbox{Vol}(W\cap \partial \Omega)}{(4\pi (2\mu +\lambda) t)^{(n-1)/2}}
 + O(t^{1-n/2}) \quad \; \mbox{as} \, \; t\to 0^+.\nonumber\end{eqnarray}
     It follows from  (\ref{c4-23}), (\ref{2020.7.5-2}), (\ref{2020.7.6-4}), (\ref{2020.7.13-5}), (\ref{3.11}) and (\ref{3..15}) that
 \begin{eqnarray} \label{a-4-1-3} && \int_{W\cap \Omega}\mbox{Tr}\big({\mathbf{K}}(t, x, x) \big) dx  \mp \int_{W\cap \Omega}\mbox{Tr}\big( {\mathbf{K}}(t, x, \overset{*}{x})\big)dx \\
 &&\qquad\;\; =\bigg[  \frac{n-1}{(4\pi \mu t)^{n/2}}  +  \frac{1}{(4\pi (2\mu+\lambda) t)^{n/2}}\bigg]\mbox{Vol}(W\cap\Omega) \nonumber\\
&&\qquad \;\, \;\,\;  \;\mp \frac{1}{4}\bigg[(n-1) \frac{\mbox{Vol}(W\cap \partial \Omega)}{(4\pi \mu t)^{(n-1)/2}}+ \frac{\mbox{Vol}(W\cap \partial \Omega)}{(4\pi (2\mu +\lambda) t)^{(n-1)/2}}\bigg]\nonumber \\
 &&\quad \qquad \; \; +O(t^{1-n/2})\quad \; \mbox{as} \, \; t\to 0^+.\nonumber\end{eqnarray}
  Note that \begin{eqnarray}  \label {2.12.10-6} \sum_{k=1}^\infty e^{- t \tau_k^{-} } = \int_\Omega \mbox{Tr}\, \big( \mathbf{K}^- (t, x,x) \big) \,dx  
 = \int_\Omega  \mbox{Tr}\, \Big( \mathbf{E} (t,x,x) -\mathbf{K} (t, x, \overset{*}{x})\Big) \,dx    \end{eqnarray} 
and   \begin{eqnarray} \label{23.12.10-7}
\sum_{k=1}^\infty e^{- t \tau_k^{+} } 
\!\!\!\!  &&  \!\!\!= \int_\Omega \mbox{Tr}\, \big( \mathbf{K}^+ (t, x,x) \big) \,dx  \\
 && \!\!\! = \int_\Omega  \mbox{Tr} \,\Big( \mathbf{K} (t,x,x) +\mathbf{K} (t, x, \overset{*}{x})\Big) \,dx  - \int_\Omega  \mbox{Tr}\, \Big( \mathbf{H}  (t, x, x) \Big) \,dx. \nonumber \end{eqnarray} 
From  (\ref{a-4-1-3})--(\ref{23.12.10-7}) and (\ref{23.12.9-1}),  we get  (\ref{1-7}).
   $\square$

\addvspace{16.9mm}

\centerline {\bf  Acknowledgments}

\vskip 0.58 true cm

 This research was supported by NNSF of China (12271031)  and  NNSF
of China  (11671033/A010802).

\addvspace{7mm}

\addvspace{12mm}

\addvspace{5mm}

\end{document}